\documentclass[11pt, reqno]{amsart}
\usepackage{tikz-cd}

\usepackage{amsfonts}
\usepackage{mathrsfs}
\usepackage{chngcntr}
\usepackage[all]{xy}
\usepackage{amsxtra,amsthm,amsmath,amssymb}

\theoremstyle{plain}

\newcommand{\cleqn}{\setcounter{equation}{0}}
\newcommand{\clth}{\setcounter{theorem}{0}}
\newcommand {\sectionnew}[1]{\section{#1}\cleqn\clth}

\let\oldtocsection=\tocsection

\let\oldtocsubsection=\tocsubsection

\let\oldtocsubsubsection=\tocsubsubsection

\renewcommand{\tocsection}[2]{\hspace{0em}\oldtocsection{#1}{#2}}
\renewcommand{\tocsubsection}[2]{\hspace{2em}\oldtocsubsection{#1}{#2}}
\renewcommand{\tocsubsubsection}[2]{\hspace{4.5em}\oldtocsubsubsection{#1}{#2}}

\makeatletter
\def\subsection{\@startsection{subsection}{2}
 \z@{.5\linespacing\@plus.7\linespacing}{-.5em}
 {\normalfont\bfseries}}
\def\subsubsection{\@startsection{subsubsection}{3}
 \z@{.5\linespacing\@plus.7\linespacing}{-.5em}
 {\normalfont\bfseries}}
\makeatother

\newtheorem{theorem}{Theorem}[section]
\newtheorem{lemma}[theorem]{Lemma}
\newtheorem{definition-theorem}[theorem]{Definition-Theorem}
\newtheorem{proposition}[theorem]{Proposition}
\newtheorem{corollary}[theorem]{Corollary}
\newtheorem{definition}[theorem]{Definition}
\newtheorem{example}[theorem]{Example}
\newtheorem{remark}[theorem]{Remark}
\newtheorem{notation}[theorem]{Notation}
\newtheorem{assumption}[theorem]{Assumption}
\newtheorem{assumption-notation}[theorem]{Assumption-Notation}
\newtheorem{lemma-definition}[theorem]{Lemma-Definition}
\newtheorem{lemma-notation}[theorem]{Lemma-Notation}
\newtheorem{question}[theorem]{Question}
\newtheorem{remark-definition}[theorem]{Remark-Definition}
\newcommand \bthm[1] { \begin{theorem}\label{t#1} }
\newcommand \ble[1] { \begin{lemma}\label{l#1} }

\newcommand \bpr[1] { \begin{proposition}\label{p#1} }
\newcommand \bco[1] { \begin{corollary}\label{c#1} }
\newcommand \bde[1] { \begin{definition}\label{d#1}\rm }
\newcommand \bex[1] { \begin{example}\label{e#1}\rm }
\newcommand \bre[1] { \begin{remark}\label{r#1}\rm }

\newcommand \bnota[1] {\begin{notation}\label{n#1}\rm }
\newcommand \bas[1] { \begin{assumption}\label{a#1}\rm }

\newcommand \bld[1] { \begin{lemma-definition}\label{ld#1} }

\newcommand \bqu[1] { \begin{question}\label{q#1}\rm }

\newcommand {\ethm} { \end{theorem} }
\newcommand {\ele} { \end{lemma} }

\newcommand {\epr} { \end{proposition} }
\newcommand {\eco} { \end{corollary} }
\newcommand {\ede} { \end{definition} }
\newcommand {\eex} { \end{example} }
\newcommand {\ere} { \end{remark} }
\newcommand {\enota} { \end{notation} }
\newcommand {\eas} {\end{assumption}}

\newcommand {\eld}{ \end{lemma-definition} }

\newcommand {\equ} {\end{question}}

\newcommand \thmref[1]{Theorem \ref{t#1}}
\newcommand \leref[1]{Lemma \ref{l#1}}
\newcommand \prref[1]{Proposition \ref{p#1}}
\newcommand \coref[1]{Corollary \ref{c#1}}
\newcommand \deref[1]{Definition \ref{d#1}}
\newcommand \exref[1]{Example \ref{e#1}}
\newcommand \reref[1]{Remark \ref{r#1}}
\newcommand \lb[1]{\label{#1}}

\def \CC {{\mathbb C}}
\def \ZZ {{\mathbb Z}}

\def \O  {{\mathcal{O}}}

\def \al {\alpha}

\def \la {\lambda}

\def \Map {\longmapsto}
\def \ra  {\rightarrow}           

\def \lrw {\longrightarrow}

\def \la {\langle}
\def \ra {\rangle}

\def \g  {\mathfrak{g}}   
\def \h  {\mathfrak{h}}

\def \n  {\mathfrak{n}}

\def \b  {\mathfrak{b}}

\def \a  {\mathfrak{a}}

\def \t  {\mathfrak{t}}

\def \sA {{\scriptscriptstyle A}}
\def \sB {{\scriptscriptstyle B}}
\def \sC {{\scriptscriptstyle C}}

\def \sF {{\scriptscriptstyle F}}
\def \sG {{\scriptscriptstyle G}}

\def \sM {{\scriptscriptstyle M}}

\def \sQ {{\scriptscriptstyle Q}}

\def \sX {{\scriptscriptstyle X}}
\def \sY {{\scriptscriptstyle Y}}
\def \sZ {{\scriptscriptstyle Z}}

\newcommand{\beqa}{\begin{eqnarray*}}                     
\newcommand{\eeqa}{\end{eqnarray*}}
\def \hs {\hspace{.2in}}
\def \lara {\la \, , \, \ra}

\def \bfv {\underline{v}}

\def \dw {\dot{w}}
\def \dv {{\dot{v}}}
\def \du {\dot{u}}

\def \bfu {{\bf u}}
\def \bfv {{\bf v}}
\def \bfw {{\bf w}}

\def \lam {\lambda}

\def \pist {\pi_{\rm st}}

\def \piG {{\pi_{\scriptscriptstyle G}}}

\def \piGs {{\pi_{{\scriptscriptstyle G}^*}}}
\def \piX {{\pi_{\scriptscriptstyle X}}}
\def \piY {{\pi_{\scriptscriptstyle Y}}}

\def \Lam {{\Lambda}}

\def \wF_mn {\wF_m \times \wF_n}
\def \wF_mnC {\wF_{m, n, \, \sC}}

\def \wF {\widetilde{F}}
\def \swF {{\scriptscriptstyle{\wF}}}

\def \Guv {G^{\bfu, \bfv}}

\def \Ouuinv {\O^{(\bfu, \bfu^{-1})}}
\def \Ouuinve {\O^{(\bfu, \bfu^{-1})}_e}

\def \Ou {\O^\bfu}
\def \sG {{\scriptscriptstyle G}}

\def \sQ {{\scriptscriptstyle Q}}
\def \tF {\widetilde{F}}
\def \stF {{\scriptscriptstyle \widetilde{F}}}

\def \calG {{\mathcal{G}}}
\def \calI {{\mathcal{I}}}

\def \calJ {{\mathcal{J}}}

\def \mX {{\mathfrak{X}}}

\def \oapi \overrightarrow{\pi}
\def \calA {{\mathcal{A}}}

\def \calB {{\mathcal{B}}}

\def \calC {{\mathcal{C}}}

\def \Maps {\longmapsto}

\def \dbfu {\dot{\bfu}}
\def \dbfv {\dot{\bfv}}
\def \dbfw {\dot{\bfw}}

\def \tpi {\widetilde{\pi}}
\def \opi {\overline{\pi}}
\def \hpi {\widehat{\pi}}

\def \hTheta {\widehat{\Theta}}

\def \sN {{\scriptscriptstyle N}}

\def \calQ {{\mathcal{Q}}}
\def \calS {{\mathcal{S}}}

\def \TT {\mathbb{T}}
\def \Gauu {\Gamma^{(\bfu, \bfu^{-1})}}

\def \dbu {\dot{\bfu}}

\def \Gdbu {\calG^{\dbu, \dbu}}

\def \Guu {G^{\bfu, \bfu}}

\def \wbfw {{\widetilde{\bfw}}}
\def \Tbfw {{\widetilde{T}^{\bfw}}}
\def \Tbfu {{\widetilde{T}^{\bfu}}}

\begin{document}

\setlength{\baselineskip}{1.2\baselineskip}
\title[Configuration Poisson groupoids of flags]{Configuration Poisson groupoids of flags}
\author{Jiang-Hua Lu}
\address{
Department of Mathematics   \\
The University of Hong Kong \\
Pokfulam Road               \\
Hong Kong}
\email{jhlu@maths.hku.hk}
\author{Victor Mouquin}
\address{}
\email{mouquinv@sjtu.edu.cn}
\author{Shizhuo Yu}
\address{
School of Mathematical Sciences\\
Nankai University}
\email{yusz@nankai.edu.cn}
\date{}
\begin{abstract} Let $G$ be a connected complex semi-simple Lie group and ${\mathcal{B}}$ its flag variety. For every positive integer $n$,
we introduce a Poisson groupoid over $\calB^n$, called the 
$n$th {\it total configuration Poisson groupoid of flags} of $G$,  which contains a family of Poisson sub-groupoids  whose total spaces
 are  {\it generalized double Bruhat cells} and bases {\it generalized Schubert cells} in ${\mathcal{B}}^n$.  
Certain symplectic leaves of these Poisson sub-groupoids are then shown to be symplectic groupoids over generalized Schubert cells.
We also give explicit descriptions of symplectic leaves in three series of Poisson varieties associated to $G$.
\end{abstract}
\maketitle
\tableofcontents
\addtocontents{toc}{\protect\setcounter{tocdepth}{1}}

\vspace{-.4in}
\sectionnew{Introduction and statements of results}\label{s:intro}
\subsection{Introduction}\label{ss:intro} Symplectic groupoids, and more generally Poisson groupoids, were introduced by M. Karasev \cite{Kara} and
A. Weinstein \cite{We:sym-groupoid, We:coiso} to study singular
foliations and quantizations of Poisson manifolds. 
A Poisson manifold
is said to be integrable if it is the base of a symplectic groupoid.
While not every Poisson manifold is 
integrable (see \cite{C-F:inte-Poi} for the obstructions), ``natural" Poisson manifolds are expected to have natural integrations to symplectic groupoids. 
When the Poisson manifold is algebraic, one would also want the symplectic groupoids to be algebraic. 
It is a fundamental problem of Poisson geometry to construct explicit and systematic examples of symplectic groupoids, especially in the category of 
algebraic Poisson manifolds. 

Lie theory provides a rich class of Poisson manifolds: 
every connected complex semi-simple Lie group $G$ carries a {\it standard} multiplicative Poisson structure $\pist$, 
defined using the choice of a pair $(B, B_-)$ of opposite Borel subgroups of $G$ (see $\S$\ref{ss:pist}), and many important manifolds in Lie theory 
carry Poisson structures closely related to the Poisson Lie group $(G, \pist)$.  Four series of such Poisson manifolds have been introduced and studied 
by the first two authors in \cite{Lu-Mou:mixed, Lu-Mou:flags}. Among them are the two series of quotients of $G^n$, denoted respectively as
\[
\wF_n = G \times_B  \times \cdots \times_B G \hs \mbox{and} \hs F_n = \wF_n/B = G \times_B  \times \cdots \times_B G/B, \hs n \geq 1.
\]
Here and for the rest of the paper, we consider the right action of $G^n$ on itself by 
\begin{equation}\label{eq:Gn-Bn}
(g_1, g_2, \ldots, g_n) \cdot (h_1, h_2, \ldots, h_n) = (g_1h_1, h_1^{-1}g_2 h_2, \ldots, h_{n-1}^{-1}g_nh_n), \;\; h_j, g_j \in G.
\end{equation}
Then $\wF_n$ is the quotient of $G^n$ by $B^{n-1} \times \{e\} \subset G^n$, while $F_n$ is the quotient of $G^n$ by $B^n \subset G^n$. 
By \cite[$\S$7.1]{Lu-Mou:mixed}, the product Poisson structure $(\pist)^n$ on $G^n$
projects to well-defined Poisson structures on $\wF_n$ and on $F_n$, 
respectively denoted as $\tpi_n$ and $\pi_n$. See $\S$\ref{ss:XnYn} for a more general construction.

For $Z = \wF_n$ or $F_n$, denote by $[g_1, \ldots, g_n]_\sZ$ the image in $Z$ of
$(g_1, \ldots, g_n) \in G^n$. Let $T = B\cap B_-$, a maximal torus of $G$. Then $T$ acts on $\tF_n$ and $F_n$ respectively by
\begin{align}\label{eq:T-Fn}
&t \cdot [g_1, g_2, \ldots, g_n]_{\sF_n} = [tg_1, g_2, \ldots, g_n]_{\sF_n},\\
\label{eq:T-tFn}
&t \cdot [g_1, g_2, \ldots, g_n]_{\stF_n} = [tg_1, g_2, \ldots, g_n]_{\stF_n},
\end{align}
 preserving the Poisson structures $\tpi_n$ and $\pi_n$. A systematic study of the $T$-orbits of symplectic leaves, or {\it $T$-leaves} for short (see
\deref{de:T-leaves}), of both $(\tF_n, \tpi_n)$ and $(F_n, \pi_n)$ are given in \cite{Lu-Mou:flags}. 
In particular, it is shown in \cite[Theorem 1.3 and Theorem 1.1]{Lu-Mou:flags}  that both $(\tF_n, \tpi_n)$ and $(F_n, \pi_n)$ have finitely many $T$-leaves.
Setting $\calB = G/B$, the flag variety of $G$, note that for each $n$ one has the isomorphism
\[
F_n \lrw \calB^n, \;\; [g_1, g_2, \ldots, g_n]_{\sF_n} \longmapsto ({g_1}_\cdot B, \, (g_1g_2)_\cdot B, \, \ldots, \, (g_1g_2\cdots g_n)_\cdot B).
\]
We thus also regard $F_n$ as a product of flag varieties. Similarly, one has the isomorphism
\[
\tF_n \lrw \calB^{n-1} \times G, \;\; [g_1, g_2, \ldots, g_n]_{\stF_n} \longmapsto ({g_1}_\cdot B, \, (g_1g_2)_\cdot B, \, \ldots, \, (g_1g_2\cdots g_{n-1})_\cdot B, \,
g_1g_2 \cdots g_n).
\]

As the main results of the paper, we show that for each $n \geq 1$ the Poisson manifold $(\tF_{2n}, \tpi_{2n})$ is a Poisson groupoid
over $(F_n, \pi_n)$, and that certain $T$-leaves of $(\tF_{2n}, \tpi_{2n})$ (resp. symplectic leaves therein) are Poisson (resp. symplectic) sub-groupoids.
We also give several isomorphic models for the Poisson groupoids, which shed different lights and put them in different contexts. 
The isomorphic models are established through 
Poisson isomorphisms between $T$-leaves in various $T$-Poisson varieties, the proofs of which, being technical, are presented in the appendices. 
In Appendix $\S$\ref{s:app-leaves} in particular, we determine symplectic leaves in all the $T$-leaves of three series of $T$-Poisson varieties 
including $(\tF_n, \tpi_n)$
for all $n$. We remark that while $T$-leaves of many $T$-Poisson varieties associated to the Poisson Lie group $(G, \pist)$ have been
determined (see, for example, \cite{e-l:Grothendieck, Lu:twisted, Lu-Mou:flags}), describing the symplectic leaves therein is a harder problem
and it has been done only for the case of $(G, \pist)$ itself by M. Kogan and A. Zelevinsky in \cite{k-z:leaves}. 
Results in Appendix $\S$\ref{s:app-leaves} thus constitute a big step towards a general theory of {\it leaves in $T$-leaves}
 and is thus of independent interest.

In the rest of the Introduction, we explain our motivation and give more details of the main results of the paper. See in particular $\S$\ref{ss:models}
on identifications of the total and base spaces of the Poisson groupoids in this paper with cluster varieties studied by L. Shen and D. Weng 
in \cite{Shen-Weng:BS}, and with augmentation varieties of Legendrian links by H. Gao, L. Shen and D. Weng in \cite{GSW:aug}. 
We also point out the recent work
\cite{Daniel:moduli} by D. Alvarez, which contains a construction of a Poisson groupoid over $F_n$ as the moduli space of flat $G$-bundles over the disk with 
decorated boundary. Precise relations between the Poisson groupoids in \cite{Daniel:moduli} and the ones in this paper, and further study on  their
 dual Poisson groupoids and double symplectic groupoids will be given elsewhere.

\subsection{Generalized Schubert cells and configuration Poisson groupoids of flags}\label{ss:Ou-intro}
Let $W$ be the Weyl group of $(G, T)$, and recall that the flag variety
$\calB = G/B$ has the decomposition into Schubert cells\footnote{In the literature, $BuB/B$, for $u \in W$, is sometimes 
called a Bruhat cell in $G/B$. 
In this paper, we use the term Schubert cell, reserving the term {\it Bruhat cell} for the sub-manifold $BuB$ in $G$ as suggested by A. Berenstein.}
$\O^u :=BuB/B$, where $u \in W$. Similarly, for $n \geq 1$ and $\bfu = (u_1, \ldots, u_n) \in W^n$, denote the image of 
$Bu_1 B \times \cdots \times Bu_nB$ in $F_n$ as
\[
\Ou =  (Bu_1 B) \times_B \cdots \times_B (Bu_nB)/B \subset F_n.
\]
One has the disjoint union decomposition 
\[
F_n = \bigsqcup_{\bfu \in W^n} \O^\bfu.
\]
For $\bfu \in W^n$, $\Ou \subset F_n$ is called a {\it generalized Bruhat cell} in \cite{Elek-Lu:BS, Lu-Mou:flags}. In this paper, we will
refer to them as {\it generalized Schubert cells} to be consistent with the case when $n = 1$. 
By \cite[Theorem 1.1]{Lu-Mou:flags}, each 
 generalized Schubert cell $\O^\bfu \subset F_n$,  being a union of (finitely many) $T$-leaves of $\pi_n$, 
is a Poisson sub-manifold of $(F_n, \pi_n)$.
Following \cite{Elek-Lu:BS, Lu-Mou:flags}, the
restriction of $\pi_n$ to $\O^\bfu$, again denoted as $\pi_n$, is called the {\it standard Poisson structure} on $\O^\bfu$.

In this paper, we first give in $\S$\ref{s:general}
a general construction of a series of Poisson groupoids associated to any Poisson Lie group and a closed Poisson Lie sub-group.
The natural Poisson groupoid structure on $(\tF_{2n}, \tpi_{2n})$
over $(F_n, \pi_n)$  is then a special case of the general construction (see \thmref{:tF2n-00}). We introduce the sub-manifold
\[
\Gamma_{2n} \stackrel{{\rm def}}{=} \{[g_1, g_2, \ldots, , g_{2n}]_{\stF_{2n}}: g_1g_2 \cdots g_{2n} \in B_-\}
\]
of $\tF_{2n}$, which, by \prref{:tFn-0} and \thmref{:tF2n-00}, is a union of $T$-leaves of $(\tF_{2n}, \tpi_{2n})$ and a Poisson sub-groupoid of 
the Poisson groupoid $(\tF_{2n}, \tpi_{2n}) \rightrightarrows (F_n, \pi_n)$. We call 
\[
\theta_{\pm}: \; (\Gamma_{2n}, \tpi_{2n}) \rightrightarrows (F_n, \pi_n), \hs n \geq 1,
\]
the {\it $n$th total configuration Poisson groupoid of flags of $G$} or simply the {\it a total configuration Poisson groupoid}, where
$\theta_+$ and $\theta_-$ are  the source and target maps. 
 For each $\bfu \in W^n$, let
\[
\Gauu \stackrel{\rm def}{=}\theta_+^{-1}(\Ou) \cap \theta_-^{-1}(\Ou) \rightrightarrows \Ou
\]
be the {\it full sub-groupoid}  of $\Gamma_{2n} \rightrightarrows F_n$ over $\Ou$ (see \deref{:wide-full}) . We show in \thmref{:Gauu} that 
$\Gauu$ is a single $T$-leaf of  $(\Gamma_{2n}, \tpi_{2n})$. Consequently 
\begin{equation}\label{eq:Gauu-intro}
\theta_+, \theta_-: \;\; (\Gauu, \tpi_{2n}) \rightrightarrows (\Ou, \pi_n)
\end{equation}
is an (algebraic) Poisson sub-groupoid of  $(\Gamma_{2n}, \tpi_{2n}) \rightrightarrows (F_n, \pi_n)$. Furthermore, 
we show that 
the symplectic leaf of 
$(\Gauu, \tpi_{2n})$ through the section of units of the groupoid in \eqref{eq:Gauu-intro} is an (algebraic) symplectic groupoid over $(\Ou, \pi_n)$. 
We call $(\Gauu, \tpi_{2n}) \rightrightarrows (\Ou, \pi_n)$ a {\it special configuration Poisson groupoid of flags of $G$}.

Our interest in generalized Schubert cells $\Ou$ with their standard Poisson structures $\pi_n$, for 
$\bfu \in W^n$, stems from some of their remarkable features
in relation to cluster algebras and to other Poisson manifolds related to the Poisson Lie group $(G, \pist)$.

First of all, it is shown in \cite{Elek-Lu:BS} that for any $\bfu = (u_1, \ldots, u_n) \in W^n$, one can use root subgroups of $G$ and
reduced decompositions for each $u_i$ to 
parametrize $\O^\bfu$ by $\CC^{l(\bfu)}$, thus obtaining the so-called 
{\it Bott-Samelson coordinates} $(z_1, \ldots, z_{l(\bfu)})$ on $\Ou$ and a Poisson bracket $\{\, ,\, \}_\bfu$ on 
$\CC[z_1, z_2, \ldots, z_{l(\bfu)}]$. 
Here $l(\bfu) = l(u_1) + \cdots + l(u_n)$ and $l(u_i)$ is the length of $u_i$.
Explicit formulas for $\{\, , \, \}_\bfu$ are given in \cite[Theorem 5.15]{Elek-Lu:BS} in terms of root strings and structure
constants of the Lie algebra of $G$. 
In particular, it is shown in \cite{Elek-Lu:BS} that the polynomial Poisson algebra
$(\CC[z_1, z_2, \ldots, z_{l(\bfu)}], \{\, ,\, \}_\bfu)$ is a {\it symmetric Poisson CGL extension} in the sense of K. Goodearl and M. Yakimov,
a special class of Poisson polynomial algebras
introduced and studied in \cite{GY:PNAS, GY:Poi-CGL} by the same authors in the context of cluster algebras.

Secondly, generalized Schubert cells with the standard Poisson structures
form basic building blocks for 
many of the Poisson manifolds associated to the Poisson Lie group $(G, \pist)$.  
For example, it is shown in \cite{Lu-Yu:BS-atlas} that a number of Poisson homogeneous spaces $(G/Q, \pi_{\sG/\sQ})$ of 
the Poisson Lie group $(G, \pist)$,
including $(G, \pist)$ itself and $(G/B, \pi_1)$, 
admit so-called {\it Bott-Samelson atlases} which are built out of generalized Schubert cells,  and  the Poisson structure $\pi_{\sG/\sQ}$
 is presented as symmetric Poisson CGL extensions in all of the coordinate charts of the Bott-Samelson atlas. 
We refer to \cite{Lu-Yu:BS-atlas} for more detail. 

The explicit and natural Poisson and symplectic groupoids over the generalized Schubert cells $(\O^u, \pi_n)$
constructed in this paper add another dimension to 
this distinguished class of Poisson manifolds. The Poisson groupoids $(\Gauu, \tpi_{2n}) \!\rightrightarrows \!(\Ou, \pi_n)$
are interesting on their own. Indeed, 
we give three additional isomorphic models of $(\Gauu, \tpi_{2n}) \!\rightrightarrows \!(\Ou, \pi_n)$, each having
advantages over the others and putting
these Poisson and symplectic groupoids in different perspectives. 
We now give more details on these models. 

\subsection{Two isomorphic models of $(\Gauu, \tpi_{2n}) \rightrightarrows (\Ou, \pi_n)$}\label{ss:models}
Let $\calB = G/B$ be the flag variety of $G$, and let $\calA = G/N$ be the
{\it decorated flag variety}  (also known as the basic affine space) of $G$,  where
 $N$ be the unipotent sub-group of $B$.  Let $\calA^o = B_-N/N$, an open subvariety of $\calA$. For $n \geq 1$, set
$\calC_{2n} = \calB^{2n-1} \times \calA^o$.
Under natural isomorphisms $\calC^{2n} \to \Gamma_{2n}$ and $F_n \to \calB^n$,
the Poisson groupoid
$(\Gamma_{2n}, \tpi_{2n}) \rightrightarrows (F_n, \pi_n)$ becomes the Poisson groupoid 
\begin{equation}\label{eq:C2n-intro}
(\calC_{2n}, \, \hpi_{2n}) \rightrightarrows (\calB^n, \, \opi_n),
\end{equation}
 see \thmref{:C2n} for detail.
Correspondingly, for each $\bfu \in W^n$, the special configuration Poisson groupoid 
$(\Gauu, \tpi_{2n}) \rightrightarrows (\Ou, \pi_n)$ is then isomorphic to a Poisson sub-groupoid 
\[
(\calC^{(\bfu, \bfu^{-1})}, \, \hpi_{2n}) \rightrightarrows (\calB^\bfu, \opi_n)
\]
 of the Poisson groupoid in \eqref{eq:C2n-intro}, where 
$\bfu^{-1} = (u_n^{-1}, \ldots, u_2^{-1}, u_1^{-1}) \in W^n$ if $\bfu = (u_1, u_2, \ldots, u_n)$ (see \coref{:Cuu-sub}).
The variety $\calC^{(\bfu, \bfu^{-1})} \subset \calC_{2n}$ consists of sequences of flags with relative Tits distances prescribed by $\bfu$ (see \leref{le:Cw})
and is an example of a {\it decorated double Bott-Samelson cell} introduced by L. Shen and D. Weng in 
\cite{Shen-Weng:BS}. 
We explain in $\S$\ref{ss:Fn-mixed} that $(\calB^n, \opi_n)$ is a {\it mixed product} of $n$ copies of the Poisson variety $(\calB, \pi_1)$ in the sense defined in 
\cite{Lu-Mou:mixed}, and a similar statement holds for the Poisson variety $(\calC_{2n}, \hpi_{2n})$. 
By \cite{Shen-Weng:BS}, each $\calC^{(\bfu, \bfu^{-1})}$ is a Poisson cluster variety when $G$ is of adjoint type.
Relations between the Poisson 
groupoid structure defined in this paper and the cluster structure on these varieties defined in \cite{Shen-Weng:BS} 
will be a very interesting topic to explore (see \cite[Remark 1.10]{Shen-Weng:BS}).

For the second isomorphic model, consider the open sub-manifold $F_n^o$ of $F_n$ given by 
\[
F_n^o = \{[g_1, g_2, \ldots, g_n]_{\sF_n}: \; g_1g_2 \cdots g_n \in B_-B/B\},
\]
and for 
$\bfw \in W^n$, let 
\begin{equation}\label{eq:Owe-intro}
\O^\bfw_e =\O^\bfw \cap F_n^o  = \{[g_1, g_2, \ldots, g_n]_{\sF_n} \in \O^{\bfw}: g_1g_2 \cdots g_n \in B_-B/B\}.
\end{equation}
By \cite[Theorem 1.1]{Lu-Mou:flags}, $\O^\bfw_e$ is the open $T$-leaf of $(\O^\bfw,\pi_n)$. The Poisson variety
$(F_n^o, \pi_n)$ has a natural {\it $T$-extension} $(F_n^o \times T, \pi_n \bowtie 0)$, whose $T$-leaves with respect to the
diagonal action of $T$ are precisely all the sub-varieties $\O^\bfw_e \times T$ for $\bfw \in W^n$.
Here $\pi_n \bowtie 0$, as a Poisson structure on $F_n^o \times T$, is the 
sum of the product Poisson structure $(\pi_n, 0)$ and a certain mixed term defined using the $T$-action on $F_n^o$ 
(see \eqref{eq:pix0}). We show in \coref{:Cuu-sub} that, via a Poisson isomorphism 
$J_{2n}: (\Gamma_{2n}, \tpi_{2n}) \to (F_{2n}^o \times T, \pi_{2n} \bowtie 0)$ and for each $\bfu \in W^n$, 
the Poisson groupoid $(\Gauu, \tpi_{2n}) \rightrightarrows (\Ou, \pi_n)$ is isomorphic to a Poisson groupoid
\begin{equation}\label{eq:OC0} 
(\Ouuinve \times T, \, \pi_{2n} \bowtie 0) \rightrightarrows (\Ou, \pi_n).
\end{equation}

The advantages of the isomorphic model in \eqref{eq:OC0} are at least three-fold.
First of all, as  
$(\Ouuinve \times T, \pi_{2n} \bowtie 0)$ is the $T$-extension of $(\O^{(\bfu, \bfu^{-1})}_e, \pi_{2n})$,
we can apply results in \cite{Elek-Lu:BS, Lu-Mi:Kostant} on arbitrary generalized Schubert cells and their $T$-extensions.
In particular, the model in \eqref{eq:OC0} allows us to describe all the symplectic leaves of $(\Gauu, \tpi_{2n})$, thereby proving that 
the symplectic leaf of $(\Gauu, \tpi_{2n})$ through the section of units of the groupoid $\Gauu \rightrightarrows \O^\bfu$ is a symplectic groupoid over
$(\O^{\bf u}, \pi_n)$. 

Secondly, each $(\Ouuinve \times T, \pi_{2n} \bowtie 0)$ is a (localization of a) 
Poisson symmetric CGL extension by \cite[Theorem 5.12]{Elek-Lu:BS} and 
\cite[$\S$4.2]{Lu-Yu:BS-atlas}, so
the Goodearl-Yakimov theory
in \cite{GY:Poi-CGL} gives a cluster variety structure on $\Ouuinve \times T$. The isomorphisms $\Ouuinve\times T \cong \Gauu \cong \calC^{(\bfu, \bfu^{-1})}$
thus provide tools for future research on comparing the cluster structure on $\Ouuinve \times T$ via the Goodearl-Yakimov theory
with that on $\calC^{(\bfu, \bfu^{-1})}$ established  by L. Shen and D. Weng in \cite{Shen-Weng:BS}. 

Thirdly, in their work \cite{GSW:aug} on Lagrangian fillings of Legendrian links, H. Gao, L. Shen, and D. Weng show that 
varieties of the form $\O^\bfw_e$ for $G = SL_n$ are isomorphic to  {\it augmentation varieties} of certain 
positive braid Legendrian links. It would be very interesting to explore connections between the Poisson groupoid structure on $\Ouuinve \times T$ 
in this paper and the results in \cite{GSW:aug}.

\subsection{The third isomorphic model via generalized double Bruhat cells}\label{ss:third}
For the third isomorphic model for the Poisson groupoid $(\Gauu, \tpi_{2n}) \rightrightarrows (\Ou, \pi_n)$, we first explain some background. 
Recall  from \cite{FZ:double} that 
associated to each pair $u,v \in W$ one has the double Bruhat cell $G^{u, v} =BuB \cap B_-vB_-$,  and that
(see, for example, \cite{k-z:leaves}) the decomposition 
$G = \bigsqcup_{u, v \in W} G^{u, v}$ is that of $(G, \pist)$ into $T$-leaves for the $T$-action on $G$ by left translation.
The first two authors proved in  \cite{Lu-Mou:groupoids}
that for any $u \in W$, $(G^{u, u}, \pist)$ is a Poisson groupoid over $(\O^u, \pi_1)$, but the groupoid structure
depends on a choice $\du$ of a 
representative of $u$ in the normalizer subgroup $N_G(T)$ of $T$ in $G$.  Using the description of
symplectic leaves of $\pist$ in $G^{u, u}$ given by M. Kogan and A. Zelevinsky in \cite{k-z:leaves}, it is then proved in \cite{Lu-Mou:groupoids} that
all symplectic leaves of $(G^{u, u}, \pist)$ are symplectic groupoids over $(\O^u, \pi_1)$.

 Generalizing the decomposition of $G$ into double Bruhat cells, the first two authors introduced in \cite{Lu-Mou:flags}
a Poisson manifold $(G_{n, n}, \tpi_{n, n})$, for each integer $n \geq 1$, and its decomposition 
\[
G_{n, n} = \bigsqcup_{\bfu, \bfv \in W^n} G^{\bfu, \bfv}
\]
into $T$-leaves, where each $G^{\bfu, \bfv}$ is called a {\it generalized double Bruhat cell} (one has $(G_{1, 1}, \pi_{1, 1}) \cong (G, \pist)$). 
The second author then showed in \cite{Mou:Guu} that 
 for every $\bfu \in W^n$, $(G^{\bfu, \bfu}, \tpi_{n, n})$ is a Poisson groupoid over $(\Ou, \pi_n)$, where again the 
groupoid structure depends on a  choice $\dbfu = (\du_1, \ldots, \du_n) \in N_G(T)^n$ representing $\bfu$. 
The question of
whether or not symplectic leaves of $(G^{\bfu, \bfu}, \tpi_{n, n})$ are symplectic groupoids over $(\Ou, \pi_n)$ was left unanswered in 
\cite{Mou:Guu}, due to the fact that description of symplectic leaves of $(G^{\bfu, \bfu}, \tpi_{n, n})$ was not available. 

In this paper, we show that for any $\bfu \in W^n$, each choice $\dbfu \in N_G(T)^n$ representing $\bfu$ 
gives a $T$-equivariant  Poisson embedding  
\begin{equation}\label{eq:Eu-intro}
E_{\dbfu}: \; (\Guu, \,\tpi_{n, n})\; {\hookrightarrow}\; (\Gamma_{2n}, \,\tpi_{2n})
\end{equation}
whose image is exactly $\Gauu$.
The Poisson groupoid $(\Guu, \tpi_{n, n}) \rightrightarrows (\Ou, \pi_n)$ defined in \cite{Mou:Guu} using the choice
$\dbfu$ is then shown to become precisely the groupoid $(\Gauu, \tpi_{2n}) \rightrightarrows (\Ou, \pi_n)$ via $E_{\dbfu}$. 
Through the ($\dbfu$-dependent) Poisson embedding
$E_{\dbfu}$, we have thus an intrinsic explanation on the origin of
the Poisson groupoid structures
on $\Guu$ as well as on their dependence on the representatives $\dbfu$.
At this connection, we point out that the construction of the Poisson groupoids $(\Guu, \tpi_{n, n}) \rightrightarrows (\Ou, \pi_n)$ 
in \cite{Mou:Guu} is based on a general theory 
on local Poisson groupoids over mixed product Poisson manifolds and actions by double symplectic groupoids,
 an approach completely different from what we use in this paper. 

When $n = 1$, 
$\Gamma_{2} =(G/B) \times B_-$ is the action groupoid 
$(G/B) \times B_- \rightrightarrows G/B$
for the action of $B_-$ on $G/B$ by left translation. The fact that $(\Gamma_2, \tpi_2)$ is a Poisson groupoid over
$(G/B, \pi_1)$, and that one has the Poisson embeddings in \eqref{eq:Eu-intro} for $n = 1$, are also proved in  \cite{Lu-Mou:groupoids}. 
Putting these 
results in \cite{Lu-Mou:groupoids} (for $n = 1$) and that in \cite{Mou:Guu} (on the Poisson groupoid structures on 
$\Guu$) in one unified framework was part of the motivation for this paper.

\subsection{Organization of the paper} After a general construction in $\S$\ref{s:general} of
a series of Poisson groupoids associated to any Poisson Lie group and a closed Poisson Lie sub-group, we turn to 
the Poisson Lie group $(G, \pist)$ and its Poisson Lie sub-group $B$ in $\S$\ref{s:full}, where we introduce the 
total configuration Poisson groupoids 
\[
(\Gamma_{2n}, \tpi_{2n}) \rightrightarrows (F_n, \pi_n), \hs n \geq 1,
\]
of flags. In $\S$\ref{s:Gauu}, we discuss the Poisson sub-groupoids
\[
(\Gauu, \tpi_{2n}) \rightrightarrows (\Ou, \pi_n), \hs \bfu \in W^n, \, n \geq 1,
\]
of $(\Gamma_{2n}, \tpi_{2n}) \rightrightarrows (F_n, \pi_n)$ and their three isomorphic models.
In $\S$\ref{s:leaves} we prove that the symplectic leaf of $(\Gauu, \tpi_{2n})$ passing through the section of units of the groupoid
$\Gauu \rightrightarrows \O^\bfu$ is a  symplectic groupoid over $(\Ou, \pi_n)$.

The paper contains several extensive appendices. 
The most technical parts of the paper are the proofs that 
various isomorphisms between varieties are in fact  isomorphisms of Poisson varieties, and we present these proofs in the appendices 
$\S$\ref{s:app-mix} - $\S$\ref{s:app-Guv}.  

In $\S$\ref{s:app-mix},  we 
show how the various Poisson varieties considered in this paper are {\it mixed product Poisson varieties} as defined in \cite{Lu-Mou:mixed}.

In $\S$\ref{s:app-Jn}, we show that the  $T$-equivariant isomorphism $J_n: \Gamma_n \to F_n^o \times T$ defined in \eqref{eq:Jn} is a Poisson 
isomorphism from $(\Gamma_n, \tpi_n)$ to $(F_n^o \times T, \pi_n \bowtie 0)$.

Generalizing the case of $m = n$ from \cite[$\S$1.4]{Lu-Mou:flags}, we introduce in $\S$\ref{s:app-Guv} 
a $T$-Poisson manifold $(G_{m, n}, \tpi_{m, n})$ for any pair of integers $m, n \geq 1$ whose $T$-leaves are shown to be 
 {\it generalized double Bruhat cells} $G^{\bfu, \bfv}$, where  $(\bfu, \bfv) \in W^m \times W^n$. 
The main results of $\S$\ref{s:app-Guv}  are certain explicit $T$-equivariant Poisson isomorphisms between the single $T$-leaves
\[
(G^{\bfu, \bfv}, \tpi_{m, n}) \stackrel{\sim}{\longrightarrow}
\left(\O^{(\bfu, \bfv^{-1})}_e \times T, \,\pi_{m+n} \bowtie 0\right)
\]
for all $(\bfu, \bfv) \in W^m \times W^n$, and similar Poisson isomorphisms between $T$-leaves of
 $(\tF_n, \tpi_n)$  to those of the form 
 $(\O^{\bfw}_e \times T, \pi_{n+1} \bowtie 0)$ for $\bfw \in W^{n+1}$.
These facts illustrate again the role of generalized Schubert cells (or their $T$-extensions for the examples in this paper) as building blocks for 
Poisson varieties associated to the Poisson Lie group $(G, \pist)$.

In $\S$\ref{s:app-leaves} we determine the symplectic leaves of
$(\O^\bfw_e \times T, \pi_n \bowtie 0)$ for any $n \geq 1$ and $\bfw \in W^n$. Although only the cases of $\bfw = (\bfu, \bfu^{-1})$ are needed in the main
text of the paper, the results for arbitrary $\bfw$ allow us to determine the symplectic leaves in all generalized double Bruhat cells $G^{\bfu, \bfv}$,
thereby extending the result of M. Kogan and A. Zelevinsky in \cite{k-z:leaves} for the case of $u, v \in W$. 
We in fact describe the symplectic leaves in all the three series
\[
(F_n^o \times T, \,\pi_n \bowtie 0), \hs (G_{m+n}, \,\tpi_{m, n}), \hs (\tF_n, \,\tpi_n), \hs m, n \geq 1,
\]
of $T$-Poisson varieties. Results in $\S$\ref{s:app-leaves}
provide test stone examples towards a general theory of symplectic leaves in $T$-leaves to be carried out elsewhere.

\subsection{Notation and basic definitions}\label{ss:nota-00}
For a manifold $X$ and $1 \leq k \leq \dim X$, let $\mX^k(X)$ be the space of all
$k$-vector fields on $X$, i.e., all smooth sections of $\wedge^k TX$, where $TX$ is the tangent bundle of $X$. When $X$ is a complex manifold,
$TX$ will stand for the holomorphic tangent bundle and $\mX^k(X)$ the space of all holomorphic $k$-vector fields on $X$.
 If $X$ and $Y$ are manifolds and
if  $V_\sX\in \mX^k(X)$ and $V_\sY \in \mX^k(Y)$,  let
$(V_\sX, 0), (0, V_\sY)\in \mX^k(X \times Y)$ be given by
\[
(V_\sX, \,0)(x, y) = i_y V_\sX(x) \hs \mbox{and} \hs (0, \,V_\sY)(x, y) = i_x V_\sY(y),\hs x \in X, y \in Y,
\]
where $i_y: X \to X \times Y, \,x' \mapsto (x', y)$
for $x'\in X$, and $i_x: Y \to X \times Y, \,y'\mapsto (x, y')$ for $y' \in Y$.  This convention also extends to multi-vector fields on
$n$-fold product manifolds $X_1 \times \cdots \times X_n$ for any integer $n \geq 2$. If $V_i \in \mX^k(X_i)$ for $i \in [1, n]$, we also write
\begin{align*}
(V_1, V_2, \ldots, V_n) &= V_1 \times V_2 \times \cdots \times V_n \\
&= (V_1, 0, \ldots, 0) + (0, V_2, \ldots, 0) + \cdots + (0, 0, \ldots, V_n)
\in \mX^k(X_1 \times X_2\times  \cdots \times X_n).
\end{align*}

Recall that a Poisson manifold is a pair $(X, \pi)$, where $X$ is a manifold and
$\pi \in \mX^2(X)$, called a Poisson structure, satisfies $[\pi, \pi] = 0$, where $[\,, \, ]$ is the
Schouten bracket on $\mX^\bullet(X) = \oplus_k \mX^k(X)$. For a Poisson structure $\pi$ on $X$, define
\[
\pi^\#: \;\; T^*X \lrw TX, \;\; (\pi^\#(\alpha), \; \beta) = \pi(\alpha, \beta),
\]
where $\alpha$ and $\beta$ are any $1$-forms on $X$. If $X_1$ is a Poisson sub-manifold of $(X, \pi)$, i.e., $X_1$ is a sub-manifold of $X$ such that
$\pi(x) \in \wedge^2T_xX_1$ for all $x \in X_1$, the restriction of $\pi$ to $X_1$ will still be denoted by $\pi$, so that $(X_1, \pi)$ is a Poisson manifold.

The identity element of a group is typically denoted as $e$. 
Let $A$ be any Lie group with Lie algebra $\a$. The left and right translations on $A$ by $a \in A$ will be
denoted by $L_a$ and $R_a$ respectively, and for any integer $k \geq 1$ and
$x \in \wedge^k\a$ and $\xi \in \wedge^k\a^*$, $x^L$ and $x^R$ (resp. $\xi^L$ and $\xi^R$)
will respectively
denote the left and right invariant $k$-vector fields (resp. $k$-forms) on $A$ with value $x$ (resp. $\xi$)
at $e \in A$.
If $A$ acts on a manifold $Y$ from the right with the action map
$\rho: Y \times A \to Y$, we will also denote by $\rho$ the Lie algebra homomorphism
\[
\rho: \; \a \lrw \mX^1(X), \;\; \rho(x)(y) = \frac{d}{dt}|_{t=0} \rho(y, \exp(tx)), \hs x \in \a, \, y \in Y.
\]
Similarly, if $\lambda: A \times Y \to Y$ is a left action by $A$, one has the Lie algebra anti-homomorphism
\[
\lam: \; \a \lrw \mX^1(X), \;\; \lam(x)(y) = \frac{d}{dt}|_{t=0} \lam(\exp(tx), y), \hs x \in \a, \, y \in Y.
\]

For clarity and for convenience of the reader, we recall some terminology on groupoids and Poisson groupoids and refer to 
\cite{MX, We:sym-groupoid, We:coiso} for the basics of the subject.

\bde{:wide-full}
{\rm
1) A  sub-groupoid of  a groupoid $\Gamma \rightrightarrows M$ is said to be {\it wide}
if it contains the set of units of $\Gamma \rightrightarrows M$.

2) The 
{\it full sub-groupoid} of a groupoid $\Gamma \rightrightarrows M$ over a subset $M' \subset M$ is
the intersection $\theta_+^{-1}(M') \cap \theta_-^{-1}(M')$ as a groupoid over $M'$ whose structure maps are the restrictions of 
those for $\Gamma \rightrightarrows M$, where $\theta_+$ and $\theta_-$ are respectively the source and target maps of $\Gamma \rightrightarrows M$. 
}
\ede

\bde{de:Poi-goid}
{\rm A {\it Poisson groupoid} 
is a Lie groupoid $\Gamma \rightrightarrows M$ together with a Poisson structure $\pi$ on $\Gamma$  such that
$\{(\gamma, \gamma', m(\gamma,\gamma')): \theta_-(\gamma) = \theta_+(\gamma')\} \subset \Gamma^3$
 is a coisotropic sub-manifold of $(\Gamma^3, \pi \times \pi \times (-\pi))$, where
$\theta_+, \theta_-: \Gamma \to M$ are the source and target maps, and
\[
m: \{(\gamma, \,\gamma') \in \Gamma^2:\; \theta_-(\gamma) = \theta_+(\gamma')\} \lrw \Gamma
\]
is the partially defined
multiplication on $\Gamma$.  In such a case, there is a unique Poisson structure $\pi_\sM$ on $M$ such that $\theta_+: (\Gamma, \pi) \to (M, \pi_\sM)$
is Poisson and $\theta_-: (\Gamma, \pi) \to (M, \pi_\sM)$
is anti-Poisson, and we also say that $(\Gamma, \pi)$ is a Poisson groupoid over $(M, \pi_\sM)$. If in addition $\pi$ is non-degenerate and
$\dim \Gamma = 2\dim M$,  one says that $(\Gamma, \pi) \rightrightarrows (M, \pi_\sM)$ is a symplectic groupoid  over $(M, \pi_\sM)$. 

A Poisson (resp. symplectic) groupoid  $(\Gamma, \pi) \rightrightarrows (M, \pi_\sM)$ is said to be {\it (complex) algebraic} if 
both $\Gamma$ and $M$ are smooth algebraic manifolds over $\CC$, all structure maps of $\Gamma \rightrightarrows M$ are smooth algebraic morphisms, and 
both $\pi$ and $\pi_\sM$ are algebraic Poisson (resp. symplectic) structures.
\hfill $\diamond$}
\ede

We also recall the notion of $\TT$-leaves that will be used throughout the paper.

\bde{de:T-leaves} 
{\rm If $\TT$ is a torus, by a $\TT$-Poisson manifold we mean a Poisson manifold $(X, \piX)$ with an action of
$\TT$ by Poisson isomorphisms. For a $\TT$-Poisson manifold  $(X, \piX)$, a {\it $\TT$-orbit of symplectic leaves}, or 
a {\it $\TT$-leaf} for short, of $(X, \piX)$  is a sub-manifold $L$ of $X$ of the form 
\[
L = \bigcup_{t \in \TT} t\Sigma,
\]
 where $\Sigma$ is a symplectic leaf of $(X, \piX)$, and the map $\TT \times \Sigma \to L, (t, x) \to tx,$ is a submersion.
}
\ede

We now recall a construction from \cite[$\S$2.6]{Lu-Mi:Kostant}. 
Let $\t$ be the Lie algebra of $\TT$ and assume that $\lara$ is a symmetric non-degenerate bilinear form on $\t$. Let $\{h_i: i = 1, \ldots, r = \dim \t\}$
be an orthonormal basis of $\t$ with respect to $\lara$. Given a $\TT$-Poisson manifold $(X, \piX)$ with the $\TT$-action 
$\TT \times X \to X$, we then have the Poisson structure $\piX \bowtie_\sigma 0$ on $X \times \TT$ given by
\begin{equation}\label{eq:pix0}
\piX \bowtie_\sigma 0 = (\piX, \, 0) +  \sum_{i=1}^{r} (\sigma(h_i), 0) \wedge (0, h_i^R),
\end{equation}
where $h_i^R$ denoted the right (left) invariant vector field on $\TT$ defined by $h_i$. 
We call $\piX \bowtie_\sigma 0$ a {\it $\TT$-extension of the Poisson structure} $\piX$, and we call
$(X \times\TT, \piX \bowtie_\sigma 0)$ a {\it $\TT$-extension of the $\TT$-Poisson manifold} $(X, \piX)$.
It is easy to see that $(X \times \TT, \piX \bowtie 0)$ is a $\TT$-Poisson manifold with respect to the diagonal 
$\TT$-action. The following fact is proved in \cite[Lemma 2.23]{Lu-Mi:Kostant}.

\ble{:T-leaves}
With respect to the diagonal $\TT$-action on $X \times \TT$,  the $\TT$-leaves of $(X \times \TT, \piX \bowtie_\sigma 0)$
are precisely all the sub-manifolds of the form $L \times \TT$, where $L$ is a $\TT$-leaf of $(X, \piX)$.
\ele

\subsection{Acknowledgments} We would like to thank L. Shen and D. Weng for explaining to us their work in \cite{Shen-Weng:BS}.
The research by the first and the third authors in this paper was partially supported by
the RGC of the Hong Kong SAR, China (GRF 17304415 and GRF 17307718), and the paper was largely written
while the second author was a faculty member at the School of Mathematical Sciences of Shanghai Jiaotong University. The third author acknowledges 
support from the NSF China (1210010907). 

\section{A series of Poisson groupoids associated to Poisson Lie groups}\label{s:general}
In this section, we give a construction of a series of Poisson groupoids associated to any Poisson Lie group and a closed Poisson Lie subgroup.

\subsection{Poisson Lie group actions and gauge Poisson groupoids}\label{ss:gauge-Poi}
We refer to \cite{chari-pressley, etingof-schiffmann} and especially to \cite[$\S$2]{Lu-Mou:mixed} for basic facts and sign conventions
on Poisson Lie groups
and Lie bialgebras.

Recall first that a Poisson Lie group is a pair $(G, \piG)$, where $G$ is a Lie group and
$\piG$ a Poisson bi-vector field on $G$ which is {\it multiplicative} in the sense that
the group multiplication $(G \times G, \, \piG \times \piG) \rightarrow (G, \piG)$ 
is a Poisson map. A right Poisson action of a Poisson Lie group $(G, \piG)$ on a
Poisson manifold $(X, \piX)$ is, by definition, a Poisson map
\begin{equation}\label{eq:rho0}
\rho: \;\;(X \times G, \,  \piX \times \piG)\lrw (X, \, \piX), \;\; (x, \,g) \longmapsto xg,
\end{equation}
which also defines a right Lie group action of $G$ on $X$. Left Poisson actions of $(G, \piG)$ are defined similarly.
Recall also that a {\it coisotropic subgroup} of a Poisson Lie group $(G, \piG)$ is a
Lie  subgroup of $G$ which is also a coisotropic sub-manifold with respect to
the Poisson structure $\piG$. The following fundamental fact is proved in  \cite{STS2}.

\ble{:quot1}
Suppose that $\rho$ is a Poisson Lie group action as in \eqref{eq:rho0}, and suppose that
 $Q$ is a coisotropic subgroup of $(G, \piG)$ such
that the restricted action of $Q$ on $X$ is free and the quotient $X/Q$ is a smooth manifold. Then the Poisson structure $\piX$ projects to
a well-defined Poisson structure on $X/Q$, which will be called the quotient Poisson structure of $\piX$.
\ele

\bex{:quot00}
{\rm
Suppose that $(X, \piX)$ and $(Y, \piY)$ are Poisson manifolds,
both with right Poisson actions by a Poisson Lie group $(G, \piG)$ and such that the diagonal action
of $G$ on $X \times Y$ is free and the quotient space $(X \times Y)/G$ is a smooth manifold. Then $(X \times Y, \piX \times (-\piY))$ has the right
product Poisson action by the product Poisson Lie group $(G \times G, \piG \times (-\piG))$.
Since the diagonal of $G \times G$ is a coisotropic subgroup with respect to $\piG \times (-\piG)$, the quotient space 
$(X \times Y)/G$ has the well-defined quotient Poisson structure of $\piX \times (-\piY)$.
\hfill $\diamond$
}
\eex

Suppose now that $X \to X/G$ is a principal $G$-bundle for a Lie group $G$.  Let $(X \times X)/G$ be the quotient space for the diagonal
$G$-action on $X \times X$, and denote elements in $X/G$ by $[x]$ and in $(X \times X)/G$ by $[x_1, x_2]$, where
$x, x_1, x_2 \in X$.  Recall that the {\it gauge groupoid} of  $X \to X/G$ is the manifold
$(X \times X)/G$ with the following Lie groupoid structure
over $X/G$:
\begin{align*}
&\mbox{source map}\; \theta_+: \;\; (X \times X)/G \lrw X/G: \;\; [x_1, \,x_2] \longmapsto [x_1],\\
 &\mbox{target map}\; \theta_-: \;\; (X \times X)/G \lrw X/G: \;\; [x_1, \,x_2]\longmapsto [x_2],\\
&\mbox{unit map}\; \epsilon: \;\; X/G  \lrw (X \times X)/G, \;\; [x] \longmapsto [x, \, x],\\
&\mbox{inverse map}\; {\rm inv}: \;\; (X \times X)/G \lrw (X \times X)/G:  \;\; [x_1, \, x_2] \longmapsto [x_2, x_1],\\
&\mbox{multiplication:}\; \mbox{for}\; \gamma = [x_1, \, x_2]\; \mbox{and} \;
\gamma' = [x_3, \, x_4]\; \mbox{with}\; \theta_-(\gamma) = \theta_+(\gamma'), \; \gamma \gamma' =
[x_1g, \, x_4],\\
&\hs \hs \hs \hs\hs \mbox{where}\; g \in G \; \mbox{is the unique element such that}\; x_2 g = x_3.
\end{align*}

Assume, in addition,  that $\piG$ is a multiplicative Poisson structure on $G$, $\piX$ is a Poisson structure
on $X$ such that the $G$-action on $X$ for the principal bundle $X \to X/G$ is a right Poisson action of the Poisson Lie group $(G, \pi_\sG)$
on $(X, \piX)$.
Let $\pi \in \mX^2((X \times X)/G)$ be the quotient Poisson structure of $\piX \times (-\piX)$ and
$\pi_{\sX/\sG} \in \mX^2(X/G)$ the quotient Poisson structure of $\piX$.
The following should be well-known but we have not been able to find a reference.

\ble{:gauge-Poi-0} With the Poisson structure $\pi$ on $(X \times X)/G$ and $\pi_{\sX/\sG}$ on $X/G$, the gauge groupoid
$(X \times X)/G\rightrightarrows X/G$ is a Poisson groupoid.
\ele

\begin{proof} Let $\Gamma = (X \times X)/G$ for notational simplicity.
Let $Y = \{(x_1, x_2, x_2 ,x_4,x_1,x_4) \in X^6\}$. It is clear that
$Y$ is a coisotropic sub-manifold of $(X^6, \pi_\sX^{(6)})$, where
\[
\pi_\sX^{(6)} = \piX \times (-\piX) \times \piX \times (-\piX) \times (-\piX) \times \piX \in \mX^2(X^6).
\]
Let $Y_1 = \{(\gamma, \gamma', \gamma \gamma') \in \Gamma^3: \theta_-(\gamma) = \theta_+(\gamma')\}$,
and let
\[
J: \; X^6 \lrw \Gamma^3, \;\; (x_1, x_2, x_3, x_4, x_5, x_6) \longmapsto ([x_1, x_2], [x_3, x_4], [x_5, x_6]).
\]
 Then
$J(Y) \subset Y_1$. Conversely, given $y_1 =(\gamma, \gamma', \gamma \gamma') \in Y_1$, let
$\gamma  = [x_1, \, x_2]$ and
$\gamma' = [x_3, \, x_4]$ for some $x_1, x_2, x_3, x_4 \in X$ such that
$[x_2] =[x_3]$, and let $g$ be the unique element in $G$ such that $x_2g = x_3$.
Then $y =(x_1g, x_2g, x_3, x_4,  x_1g, x_4) \in Y$ and $J(y) = y_1$.
This shows that  $J(Y) = Y_1$. Since $J: (X^6, \pi_\sX^{(6)}) \to (\Gamma^3, \pi \times \pi \times (-\pi))$ is Poisson,
it follows from the definition of coisotropic sub-manifolds that $Y_1=J(Y)$ is coisotropic in
$(\Gamma^3, \pi \times \pi \times (-\pi))$. Note also that both $\theta_+$ and $\theta_-$ are surjective submersions.
Thus $(\Gamma, \pi)$ is a Poisson groupoid. By the definition of
the Poisson structure $\pi_{\sX/\sG}$, $(\Gamma, \pi)$ is a Poisson groupoid over $(X/G, \pi_{\sX/\sG})$.
\end{proof}

We now adapt the constructions of quotient Poisson structures in \exref{:quot00} and the gauge Poisson groupoids in \leref{:gauge-Poi-0} to a setting
suitable for applications in this paper.

\bex{:quot}\cite[$\S$7.1]{Lu-Mou:mixed}
{\rm
Suppose that $(G, \piG)$ is a Poisson Lie group,  $(X, \piX)$ is  Poisson manifold with a free
{\it right} Poisson action $(x, g) \mapsto xg$ by $(G, \piG)$, and $(Y, \piY)$ is a Poisson manifold with a {\it left} Poisson action
$(g, y) \mapsto gy$
by $(G, \piG)$. Then one has the right Poisson action of the Poisson Lie group $(G \times G, \piG \times (-\piG))$ on
$(X \times Y, \, \piX \times \piY)$ given by
$(x, y) \cdot (g_1, g_2)= (xg_1, \;g_2^{-1}y)$. 
Denote by
$X \times_G Y$ the quotient of $X \times Y$ by the diagonal $G$-action
$(x, \, y) \cdot  g = (xg, \, g^{-1}y)$ for $x \in X$, $y \in Y$, and $g \in G$.
Assuming that $X \times_G Y$ is a smooth manifold, it then has the well-defined quotient Poisson structure of $\piX \times \piY$.
\hfill $\diamond$
}
\eex

Assume now that $(X, \piX)$ is a Poisson manifold with a right Poisson action $(x, g) \to xg$ by a Poisson Lie group $(G, \piG)$ and
assume that $X \to X/G$ is a principal bundle. We make the further assumption that $\kappa: X \to X$ is an {\it anti-Poisson involution} with
respect to the Poisson structure $\piX$.  One then has the unique left action $(g, x) \mapsto gx$ of $G$ on $X$ determined by
\[
\kappa(g^{-1} x) = \kappa(x)g, \hs x \in X, \, g \in G,
\]
which is a left Poisson action of $(G, \piG)$ on $(X, \piX)$. Applying the construction in \exref{:quot}, one has the quotient space
$X \times_G X$ of $X \times X$ by the right $G$-action
\[
(x_1, \, x_2) \cdot g = (x_1g, \, g^{-1}x_2), \hs x_1, x_2 \in X, \, g \in G,
\]
and the quotient Poisson structure $\pi$ on $X \times_G X$ of $\piX \times \piX$.  Denote again
elements in $X/G$ by $[x]$ and in $X \times_G X$ by $[x_1, x_2]$, where $x, x_1, x_2 \in X$, and let
again $\pi_{\sX/\sG} \in \mX^2(X/G)$ be the quotient Poisson structure
of $\piX$.

\ble{:gauge-Poi}
With the assumption and notation as above,  $(X \times_G X, \pi)$ is a Poisson groupoid over $(X/G, \pi_{\sX/\sG})$ with the following
groupoid structure:
\begin{align*}
&\mbox{source map}\; \theta_+: \;\; X \times_G X \lrw X/G: \;\; [x_1, \,x_2] \longmapsto [x_1],\\
 &\mbox{target map}\; \theta_-: \;\; X \times_G X \lrw X/G: \;\; [x_1, \,x_2]\longmapsto [\kappa(x_2)],\\
&\mbox{unit map}\; \epsilon: \;\; X/G  \lrw X \times_G X, \;\; [x] \longmapsto [x, \, \kappa(x)],\\
&\mbox{inverse map}\; {\rm inv}: \;\; X \times_G X \lrw X \times_G X:  \;\; [x_1, \, x_2] \longmapsto [\kappa(x_2), \, \kappa(x_1)],\\
&\mbox{multiplication:}\; \mbox{for}\; \gamma = [x_1, \, x_2]\; \mbox{and} \;
\gamma' = [x_3, \, x_4]\; \mbox{with}\; \theta_-(\gamma) = \theta_+(\gamma'), \; \gamma \gamma' =
[x_1g, \, x_4],\\
&\hs \hs \hs \hs \hs \mbox{where}\; g \in G \; \mbox{is the unique element such that}\; \kappa(x_2)g = x_3.
\end{align*}
\ele

\begin{proof} The map $X \times_G X \to (X \times X)/G, [x_1, x_2] \to [x_1, k(x_2)]$ is both an isomorphism of Lie groupoids and an isomorphism of
Poisson manifolds, where $(X \times X)/G$ has the gauge Poisson groupoid for the right
Poisson action $(x, g) \to xg$ as in \leref{:gauge-Poi-0}.
\end{proof}

\subsection{A series of Poisson groupoids}\label{ss:XnYn} 
Assume that $(G, \piG)$ is a Poisson Lie group and that $Q$ a closed Poisson Lie subgroup. For each integer $n \geq 1$, one then has
the  quotient spaces
\begin{equation}\label{eq:XnYn}
X_n = \overbrace{G \times_Q \cdots \times_Q G}^n \hs \mbox{and} \hs
Y_n = \overbrace{G \times_Q \cdots \times_Q G}^n/Q.
\end{equation}
The manifolds $X_n$ and $Y_n$ can be taken as successive quotient spaces as in \exref{:quot}. Consequently, one has the
well-defined quotient Poisson structures on $X_n$ and $Y_n$, i.e., well-defined projections of
the $n$-fold product Poisson structure $\pi_\sG^n$ on $G^n$, which we will  denote respectively as $\pi_{\sX_n}$ and $\pi_{\sY_n}$.
By the multiplicativity of $\piG$, one has the following Poisson
Lie group actions
\begin{align}\label{eq:Xn-G-left}
&(G, \piG) \times (X_n, \, \pi_{\sX_n}) \lrw (X_n, \, \pi_{\sX_n}), \; (g, \,[g_1, g_2, \ldots, g_n]_{\sX_n}) \longmapsto [gg_1, g_2, \ldots, g_n]_{\sX_n},\\
\label{eq:Xn-G-right}
&(X_n, \, \pi_{\sX_n}) \times (G, \piG) \lrw (X_n, \, \pi_{\sX_n}), \; ([g_1, g_2, \ldots, g_n]_{\sX_n}, \, g) \longmapsto [g_1, \ldots, g_{n-1}, g_ng]_{\sX_n},\\
\label{eq:Yn-G}
&(G, \piG) \times (Y_n, \, \pi_{\sY_n}) \lrw (Y_n, \, \pi_{\sY_n}), \; (g, \,[g_1, g_2, \ldots, g_n]_{\sY_n}) \longmapsto [gg_1, g_2, \ldots, g_n]_{\sY_n}.
\end{align}
As the inverse map on $G$ is an anti-Poisson involution for $(G, \piG)$, it follows that
\[
{\bf I}_{\sX_n}: \; (X_n, \, \pi_{\sX_n}) \lrw (X_n, \, \pi_{\sX_n}), \; [g_1, g_2, \ldots, g_n]_{\sX_n} \longmapsto [g_n^{-1}, \ldots, g_2^{-1}, g_1^{-1}]_{\sX_n}
\]
is anti-Poisson. As a direct application of \leref{:gauge-Poi} by taking $X = X_n$ with the right Poisson action by the Poisson Lie group
$(Q, \piG|_\sQ)$ given in \eqref{eq:Xn-G-right} and by taking $\kappa = {\bf I}_{\sX_n}$, one has

\bthm{:XnYn}
For any Poisson Lie group $(G, \piG)$, any closed Poisson Lie group $Q$ of $(G, \piG)$, and any positive integer $n$,
$(X_{2n}, \pi_{\sX_{2n}})$ is a Poisson groupoid over $(Y_n, \pi_{\sY_n})$ with the following groupoid structure:
\begin{align*}
&\mbox{source map}\; \theta_+: \;\; X_{2n} \to Y_n: \;\; [g_1, \ldots, g_{2n}]_{\sX_{2n}} \mapsto [g_1, \ldots, g_n]_{\sY_n};\\
 &\mbox{target map}\; \theta_-: \;\; X_{2n} \to Y_n: \;\; [g_1, \ldots, g_{2n}]_{\sX_{2n}} \mapsto [g_{2n}^{-1}, \ldots, g_{n+1}^{-1}]_{\sY_n};\\
&\mbox{unit map}\; \epsilon: \;\; Y_n \to X_{2n}, \;\; [g_1,  \ldots, g_{n}]_{\sY_{n}} \mapsto [g_1, \ldots, g_n, \,g_n^{-1}, \ldots, g_1^{-1}]_{\sX_{2n}};\\
&\mbox{inverse map}\; {\rm inv}: \;\; X_{2n} \to X_{2n}: \;\; [g_1, \ldots, g_{2n}]_{\sX_{2n}} \mapsto [g_{2n}^{-1}, \ldots, g_{1}^{-1}]_{\sX_{2n}};\\
&\mbox{multiplication:}\; \mbox{for}\; \gamma = [g_1, \ldots, g_{2n}]_{\sX_{2n}}\; \mbox{and} \;
\gamma' = [g_1^\prime,  \ldots, g_{2n}^\prime]_{\sX_{2n}} \; \mbox{with}\; \theta_-(\gamma) = \theta_+(\gamma'),\\
&\hs \hs \gamma \gamma' =
[g_1, \;\ldots, \;g_n, \;g_{n+1}\cdots g_{2n} g_1^\prime \cdots g_n^\prime g_{n+1}^\prime, \; g_{n+2}^\prime,  \ldots, g_{2n}^\prime]_{\sX_{2n}}.
\end{align*}
\ethm

\sectionnew{The total configuration Poisson groupoids of flags and isomorphic models}\label{s:full}
We now apply the construction in $\S$\ref{s:general} to the standard complex semisimple Lie group $(G, \pist)$.

\subsection{The complex semisimple Poisson Lie group $(G, \pist)$}\lb{ss:pist}
Let $G$ be a connected complex semi-simple  Lie group with Lie algebra $\g$. We recall
the so-called {\it standard multiplicative Poisson structure} on $G$ and refer the readers to
\cite{chari-pressley, etingof-schiffmann, Lu-Mou:mixed, Lu-Mou:flags} for more detail.

Fix again a pair
$(B, B_-)$ of opposite Borel subgroups of $G$ and let $T = B \cap B_-$, a maximal torus of $G$.  Let $N$ and $N-$ the respective
unipotent radicals of $B$ and $B_-$. Denote the Lie algebras of
$G, B, B_-, N, N_-$ and $T$ by  $\g, \b, \b_-,\n$ and $\h$ respectively.
Fix also a non-degenerate symmetric invariant bi-linear form
$\lara_\g$ on $\g$.
Let $\g = \h + \sum_{\alpha \in \Delta} \g_\alpha$
be the root decomposition of $\g$ with respect to $\h$, let $\Delta^+ \subset \g^*$ be the set of
positive roots with respect to $\b$, and for each $\al \in \Delta^+$, let $E_\al\in \g_\al$ and $E_{-\al}
\in \g_{-\al}$ be such that $\la E_\alpha, E_{-\alpha}\ra_\g = 1$. Let $\{h_i\}_{i = 1}^{\dim \h}$ be a basis of $\h$ which is orthonormal with respect to
the restriction of $\lara_\g$ to $\h$. Then the element
\begin{equation}\label{eq:rst}
r_{\rm st} =  \sum_{i = 1}^{\dim \h} h_i \otimes h_i + 2\sum_{\alpha \in \Delta^+} E_{-\alpha} \otimes E_\alpha \in \g \otimes \g
\end{equation}
is called the {\it standard $r$-matrix} on $\g$. Let
$\Lambda_{\rm st} \in \wedge^2\g$ be the skew-symmetric part of $r_{\rm st}$, i.e.,
\[
\Lambda_{\rm st} =  \sum_{\alpha \in \Delta^+} (E_{-\alpha} \wedge E_\alpha) =
 \sum_{\alpha \in \Delta^+} (E_{-\alpha} \otimes  E_\alpha - E_\alpha \otimes E_{-\alpha}),
\]
and let $\pist$ be the bi-vector field on $G$ defined by (see notation in $\S$\ref{ss:nota-00})
\begin{equation}\lb{eq:pist}
\pist = \Lam_{\rm st}^L - \Lam_{\rm st}^R.
\end{equation}
Then $\pist$ is a multiplicative Poisson structure on $G$, and $(G, \pist)$ is a {\it standard semi-simple Poisson
Lie group}.
It is well-known (see, for example, \cite{GY:GP}) that both $B$ and $B_-$
are Poisson Lie subgroups of $(G, \pist)$.
One thus has the Poisson Lie subgroups
$(B, \pist)$ and $(B_-, \pist)$ of $(G, \pist)$. 

\subsection{The total configuration Poisson groupoids of flags}\label{ss:t-F2n}
Continuing with the set-up in $\S$\ref{ss:pist}, we can now apply the constructions in $\S$\ref{ss:XnYn} to the Poisson Lie group
$(G, \pist)$ and its closed Poisson Lie subgroup $B$. In this particular case,  as we have already done in $\S$\ref{ss:intro}, we denote 
the Poisson spaces $X_n$ and $Y_n$
in \eqref{eq:XnYn} respectively as
\[
\tF_n =\overbrace{G \times_B \cdots \times_B G}^n \hs \mbox{and} \hs F_n =\overbrace{G \times_B \cdots \times_B G}^n/B,
\]
and we denote the quotient Poisson structures on $\tF_n$ and $F_n$ respectively as $\tpi_n$ and $\pi_n$. 
Note that we have Poisson Lie group actions
\begin{align}\label{eq:G-tFn}
&\tilde{\lambda}_n:\;\; (G, \pist) \times (\tF_n, \tpi_n) \lrw (\tF_n, \tpi_n),\;\; (g, [g_1, g_2, \ldots, g_n]_{\stF_n}) \longmapsto [gg_1, g_2, \ldots, g_n]_{\stF_n},\\
\label{eq:lam-n}
&\lambda_n: \;\; (G, \pist) \times (F_n, \pi_n) \lrw (F_n, \pi_n),\;\; (g, [g_1, g_2, \ldots, g_n]_{\sF_n}) \longmapsto [gg_1, g_2, \ldots, g_n]_{\sF_n}.
\end{align}
Recall that $T = B \cap B_-$, a maximal torus of $G$. Since $\pist|_T = 0$,  the restrictions of the actions to $T$ (see
\eqref{eq:T-Fn} and \eqref{eq:T-tFn}) make both $(\tF_n, \tpi_n)$ and $(F_n, \pi_n)$ into $T$-Poisson manifolds.  
The multiplicativity of $\pist$ also implies that we have the 
well-defined Poisson map
\begin{equation}\label{eq:mu-tFn}
\mu_{\stF_n}: \; (\tF_n, \tpi_n) \lrw (G, \pist), \;\; [g_1, g_2, \ldots, g_n]_{\stF_n} \Maps g_1g_2 \cdots g_n.
\end{equation}

For $\bfw = (w_1, \ldots, w_n) \in W^n$, set
\begin{equation}\label{eq:Ga-wn}
B \bfw B = (Bw_1B) \times_B \cdots \times_B (Bw_nB) \subset \tF_n,
\end{equation}
the image of $(Bw_1B) \times \cdots \times (Bw_nB)$ in $\tF_n$. 
One thus has 
\begin{equation}\label{eq:tFn-dcom}
\tF_n = \bigsqcup_{\bfw \in W^m, v \in W} (B\bfw B)\cap \mu_{\stF_n}^{-1}(B_-vB_-) \hs (\mbox{disjoint union}).
\end{equation}
For $w \in W$, let $l(w)$ be the length of $w$. For $\bfw = (w_1, \ldots, w_n) \in W^n$, let $l(\bfw) = l(w_1) + \cdots + l(w_n)$.
The following description of $T$-leaves of $(\tF_n, \tpi_n)$ is proved  in \cite[Theorem 1.3]{Lu-Mou:flags}. 

\bpr{:tFn-0} 1) For any $\bfw \in W^n$ and
$v \in W$, the intersection $(B\bfw B)\cap \mu_{\stF_n}^{-1}(B_-vB_-)$ is a non-empty smooth sub-manifold of $\tF_n$ of dimension $l(\bfw) + l(v) + \dim T$;

2) The decomposition in \eqref{eq:tFn-dcom}
is that of $\tF_n$ into the $T$-leaves of $\tpi_n$.
\epr

Introduce now  
\begin{equation}\label{eq:Gan}
\Gamma_n = \{[g_1, g_2, \ldots, g_n]_{\stF_n}: \, g_1g_2 \cdots g_n \in B_-\} = \mu_{\stF_n}^{-1}(B_-) \subset \tF_n.
\end{equation}
By \prref{:tFn-0}, $\Gamma_n$ is a 
union of $T$-leaves of $\tpi_n$ and thus a Poisson sub-manifold of $(\tF_n, \tpi_n)$.
We now apply \thmref{:XnYn} to the Poisson Lie group $(G, \pist)$ and its
Poisson Lie subgroup $B$.

\bthm{:tF2n-00}
For any integer $n \geq 1$, $(\tF_{2n}, \tpi_{2n})$ is a Poisson groupoid over $(F_n, \pi_n)$ with the following groupoid structure:
\begin{align*}
&\mbox{source map}\; \theta_+:  \tF_{2n} \to F_n,\;  [g_1, \ldots, g_{2n}]_{\stF_{2n}} \mapsto [g_1, \ldots, g_n]_{\sF_n};\\
 &\mbox{target map}\; \theta_-: \tF_{2n} \to F_n,\;  [g_1, \ldots, g_{2n}]_{\stF_{2n}} \mapsto [g_{2n}^{-1}, \ldots, g_{n+1}^{-1}]_{\sF_n};\\
&\mbox{unit map}\; \epsilon:  F_n \to \tF_{2n}, \; [g_1,  \ldots, g_{n}]_{\sF_{n}} \mapsto [g_1, \ldots, g_n, g_n^{-1}, \ldots, g_1^{-1}]_{\stF_{2n}};\\
&\mbox{inverse map}\; {\rm inv}: \tF_{2n} \to \tF_{2n},\;  [g_1, \ldots, g_{2n}]_{\stF_{2n}} \mapsto [g_{2n}^{-1}, \ldots, g_{1}^{-1}]_{\stF_{2n}};\\
&\mbox{multiplication:}\; \mbox{for}\; \gamma = [g_1, \ldots, g_{2n}]_{\stF_{2n}}\; \mbox{and} \;
\gamma' = [g_1^\prime,  \ldots, g_{2n}^\prime]_{\stF_{2n}} \; \mbox{with}\; \theta_-(\gamma) = \theta_+(\gamma'),\\
&\hs \hs \hs \hs \hs \gamma \gamma' =
[g_1, \;\ldots, \;g_n, \;g_{n+1}\cdots g_{2n} g_1^\prime \cdots g_n^\prime g_{n+1}^\prime, \; g_{n+2}^\prime, \, \ldots, \,g_{2n}^\prime]_{\stF_{2n}}.
\end{align*}
Furthermore, $\Gamma_{2n} \subset \tF_{2n}$ is a wide Poisson
sub-groupoid of $(\tF_{2n}, \tpi_{2n}) \rightrightarrows (F_n, \pi_n)$.
\ethm

\begin{proof}
Applying \thmref{:XnYn} directly to the Poisson Lie group $(G, \pist)$ and its
Poisson Lie subgroup $B$, we see that $(\tF_{2n}, \tpi_{2n}) \rightrightarrows (F_n, \pi_n)$ as described in
\thmref{:tF2n-00} is a Poisson groupoid.
One also checks directly from the definitions
of the structure maps of $\tF_{2n} \rightrightarrows F_n$
that $\Gamma_{2n}$
is a wide (set-theoretical) sub-groupoid of $\tF_{2n}\rightrightarrows F_n$ (\deref{:wide-full}).
As  $\Gamma_{2n}$ is a Poisson sub-manifold of
$(\tF_{2n}, \tpi_{2n})$, it is a Poisson sub-groupoid of  $(\tF_{2n}, \tpi_{2n}) \rightrightarrows (F_n, \pi_n)$ as long as we prove that $\theta_+:
\tF_{2n} \to F_n$ restricts to a
submersion from $\Gamma_{2n}$ to $F_n$. To prove this latter statement, we note the isomorphisms
\begin{align*}
&\tF_{2n} \longrightarrow (G/B)^{2n-1} \times G, \;\; [g_1, g_2, \ldots, g_{2n}]_{\stF_{2n}} 
\longmapsto (g_1 {}_\cdot B, g_1g_2 {}_\cdot B, \ldots, g_1g_2 \cdots 
g_{2n-1}{}_\cdot B, g_1g_2 \cdots g_{2n}),\\
&F_n \longrightarrow (G/B)^n, \;\; [g_1, g_2, \ldots, g_{n}]_{\sF_{n}} \longmapsto (g_1 {}_\cdot B, g_1g_2 {}_\cdot B, \ldots, 
g_1g_2 \cdots g_{n}{}_\cdot B),
\end{align*}
under which $\Gamma_{2n}$ is mapped to $(G/B)^{2n-1} \times B_-$, and $\theta_+|_{\Gamma_{2n}}: \Gamma_{2n} \to F_n$ becomes the projection of
$(G/B)^{2n-1} \times B_-$ from  the product $(G/B)^n$ of the first $n$ factors and is thus a submersion.
\end{proof}

\bde{:Ga2n}
{\rm
For $n \geq 1$
 the Poisson groupoid $(\Gamma_{2n}, \widetilde{\pi}_{2n})) \rightrightarrows (F_n, \pi_n)$ in \thmref{:tF2n-00} is called 
the {\it $n$th total configuration Poisson groupoid of flags of $G$}.
\hfill $\diamond$
}
\ede

In the next $\S$\ref{ss:C2n} and $\S$\ref{ss:F2no},  we introduce two isomorphic models of
 $(\Gamma_{2n}, \widetilde{\pi}_{2n})) \rightrightarrows (F_n, \pi_n)$.

\subsection{The Poisson groupoid $(\calC_{2n}, \hpi_{2n}) \rightrightarrows (\calB^n, \opi_n)$}\label{ss:C2n}
Recall again that  $\calB = G/B$ is the (full) flag variety of $G$, and $\calA = G/N$ is 
the decorated flag variety of $G$, where $N$ is the unipotent radical of $B$. For an integer $n \geq 1$, recall that we have set
\begin{equation}\label{eq:Cn}
\calC_n =\calB^{n-1} \times \calA^o, 
\end{equation}
where $\calA^o = B_-N/N$ is the open $B_-$-orbit in $\calA = G/N$.
Referring to an element
$f \in \calB$ as a {\it flag} and an element $\hat{f}$ 
in $\calA$ as a {\it decorated flag}, the space $\calC_n$ then consists of all $n$-tuples 
$(f_1, \ldots, f_{n-1}, \hat{f}_n)$, where $f_1, \ldots, f_{n-1}$ are flags and $\hat{f}_n$ is a decorated flag that is 
in general position with the flag $f_-$ represented by $B_-$.

The spaces $\calC_n$ appeared in \cite{Shen-Weng:BS}, where the authors consider the special cases when
$G = G_{\rm sc}$ is simply connected and when $G = G_{\rm ad}$ is of adjoint type. 
For any pair of words $(b, d)$ of length $n$ in the simple reflections in $W$ (called positive braids in \cite{Shen-Weng:BS}), 
Shen and Weng introduced certain configuration spaces of flags, denoted as   
${\rm Conf}^b_d(\calA_{\rm sc})$ and ${\rm Conf}^b_d(\calA_{\rm ad})$ and called
{\it decorated double Bott-Samelson cells} (see \cite[$\S$2.2]{Shen-Weng:BS}), which can 
be embedded in $\calC_n$ for $G = G_{\rm sc}$ and $G = G_{\rm ad}$ respectively \cite[$\S$2.3]{Shen-Weng:BS}
(see \reref{:Cw}).
As one of their main results, Shen and Weng prove in \cite[Theorem 1.1 and Theorem 1.2]{Shen-Weng:BS}  that
both ${\rm Conf}^b_d(\calA_{\rm sc})$ and ${\rm Conf}^b_d(\calA_{\rm ad})$
are smooth affine varieties; the coordinate rings $\O({\rm Conf}^b_d(\calA_{\rm sc}))$ is an upper cluster algebra, 
$\O({\rm Conf}^b_d(\calA_{\rm ad}))$ is an upper cluster Poisson algebra, and the pair $({\rm Conf}^b_d(\calA_{\rm sc}), {\rm Conf}^b_d(\calA_{\rm ad}))$
form a cluster ensemble for which the Fock-Goncharov cluster duality conjecture holds. See \cite{Shen-Weng:BS} for detail.

For $n \geq 1$, consider now the isomorphisms 
\begin{align}\label{eq:hThetan}
&\hTheta_n: \; \Gamma_{n} \lrw \calC_{n}, \; [g_1, \cdots, g_{n-1}, g_{n}]_{\stF_{n}} \longmapsto
(g_1{}_\cdot B, \ldots, g_1\cdots g_{n-1}{}_\cdot B, \, g_1\cdots g_{n-1}g_{n}{}_\cdot N),\\
\label{eq:Thetan}
&\Theta_n: \; F_n \lrw \calB^n, \;  [g_1, \cdots, g_{n-1}, g_{n}]_{\sF_{n}} \longmapsto
(g_1{}_\cdot B, \ldots, g_1\cdots g_{n-1}{}_\cdot B, \,g_1\cdots g_{n-1} g_{n}{}_\cdot B).
\end{align}
Under the isomorphisms $\hTheta_{2n}: \Gamma_{2n} \to \calC_{2n}$ and $\Theta_n: F_n \to \calB^n$,  one checks directly that
the Lie groupoid $\Gamma_{2n} \rightrightarrows F_n$ in \thmref{:tF2n-00} becomes the
following Lie groupoid $\calC_{2n}\rightrightarrows \calB^n$, where we denote an element in $\calA^o = B_-N/N$ as $b_-{}_\cdot N$ for a unique 
$b_- \in B_-$:
\begin{align*}
&\mbox{source map}\; \theta_+:  \calC_{2n}\to \calB^n:  (f_1, \ldots, f_{2n-1}, b_-{}_\cdot N) \mapsto (f_1, \ldots, f_n);\\
&\mbox{target map}\; \theta_-: \calC_{2n} \to \calB^n:  (f_1, \ldots, f_{2n-1}, b_-{}_\cdot N) \mapsto (b_-^{-1}f_{2n-1}, \ldots, b_-^{-1} f_n);\\
&\mbox{unit map} \;\epsilon: \calB^n \to \calC_{2n}, (f_1, \ldots, f_n) \to (f_1, \ldots, f_n, f_{n-1}, \ldots, f_1, e{}_\cdot N);\\
&\mbox{inverse map}\; {\rm inv}: \calC_{2n} \to \calC_{2n}:  (f_1, \ldots, f_{2n-1}, b_-{}_\cdot N)
\mapsto (b_-^{-1}f_{2n-1}, \ldots, b_-^{-1}f_2, b_-^{-1}f_1, \, b_-^{-1}{}_\cdot N);\\
&\mbox{multiplication:}\; \mbox{for}\; \gamma = (f_1, \ldots, f_{2n-1}, b_-{}_\cdot N) \; \mbox{and} \;
\gamma' = (f_1^\prime, \ldots, f_{2n-1}^\prime, b_-^\prime{}_\cdot N)  \; \mbox{with}\; \theta_-(\gamma) = \theta_+(\gamma'),\\
&\hs \hs \hs \hs \gamma \gamma' = (f_1, \; \ldots, \; f_n, \; b_-f_{n+1}^\prime, \; \ldots \; b_- f_{2n-1}^\prime, \; b_- b_-^\prime{}_\cdot N).
\end{align*}
Define Poisson structures $\hpi_n$ on $\calC_n$ and $\opi_n$ on $\calB^n$ respectively by
\begin{equation}\label{eq:hpi-opi}
\hpi_n = \hTheta_n(\tpi_n) \hs \mbox{and} \hs \opi_n = \Theta_n(\pi_n).
\end{equation}
We now have the following direct consequence of \thmref{:tF2n-00}.

\bthm{:C2n}
For any positive integer $n$, $(\calC_{2n}, \hpi_{2n})$, with the groupoid structure as above, is a algebraic Poisson groupoid over $(\calB^n, \opi_n)$.
\ethm

For each $n \geq 1$, we call $\hpi_n$ (resp. $\opi_n$) the {\it standard Poisson structure} on $\calC_n$ (resp. $\calB^n$). 
We prove in \prref{:piQ} (see also \coref{:Jmn}) 
that $(\calB^n, \opi_n)$ is a {\it mixed product} of $n$ copies of the Poisson variety $(\calB, \pi_1)$ in the sense defined in 
\cite{Lu-Mou:mixed}. A similar statement for
$(\calC_n, \hpi_n)$ is given in \reref{:pin-mix}.

\subsection{The Poisson groupoid $(F_{2n}^o \times T, \pi_{2n} \bowtie 0) \rightrightarrows (F_n, \pi_n)$}\label{ss:F2no}
For $n \geq 1$, let again 
\begin{equation}\label{eq:Fno}
F_n^o = \{[g_1, g_2, \ldots, g_n]_{\sF_n}: \; g_1g_2 \cdots g_n \in B_-B/B\},
\end{equation}
an open sub-manifold of $F_n$. 
Recall that for $g \in B_-B$, we write $g = [g]_-[g]_0[g]_+$ with $[g]_- \in N_-, [g]_0 \in T$, and $[g]_+ \in N$. 
Set also $[g]_{\geq 0} = [g]_0[g]_+$ for $g \in B_-B$. Define
\begin{equation}\label{eq:Jn}
J_n:\;\; \Gamma_n \longrightarrow F_n^o \times T,\;\; J_n([g_1, g_2, \ldots, g_n]_{\stF_n}) = ([g_1, g_2, \ldots, g_n]_{\sF_n}, \, [g_1g_2 \cdots g_n]_0).
\end{equation}
Recall that $T$ acts on $F_n$ and $\tF_n$  by \eqref{eq:T-Fn}  and \eqref{eq:T-tFn}
respectively. Let $T$ act on itself by translation and on $F_n^o \times T$ diagonally. 
Then $J_n$ is $T$-equivariant and that the inverse of $J_n$ is  given by 
\begin{equation}\label{eq:Jn-inv}
 J_n^{-1}([g_1, \ldots, g_{n-1}, g_n]_{\sF_n}, t) = [g_1, \ldots, g_{n-1}, g_n [g_1 \cdots g_{n-1}g_n]_{\geq 0}^{-1}t]_{\stF_n}.
\end{equation}
On the other hand, we have the Poisson structure $\pi_n \bowtie 0$ on $F_n \times T$ which is the 
$T$-extension of $\pi_n$ with respect to the $T$-action $\lambda_n$ on $F_n$ (see $\S$\ref{ss:nota-00} and \eqref{eq:lam-n}), namely, 
\begin{equation}\label{eq:pin-bow0}
\pi_n \bowtie 0  := (\pi_n, 0) + \sum_{i = 1}^{\dim \h} (\lam_n(h_i), 0) \wedge (0, h_i^R) \in \mX^2(F_n \times T),
\end{equation}
where 
$\{h_i\}$ is a basis of the Lie algebra 
$\h$ of $T$ orthonormal with respect to $\lara_\g|_\h$,  and for $x \in \h$,
$x^R$ is again the right (and left) invariant vector field on $T$ with value $x$ at the
identity element. 

\bthm{:Jn}
For any $n \geq 1$, $J_n:  (\Gamma_n, \, \tpi_n) \to (F_n^o \times T, \, \pi_n \bowtie 0)$ is a 
 Poisson isomorphism.
\ethm

\begin{proof}
This is \thmref{:Jn-apdx} in the $\S$\ref{s:app-Jn}.
\end{proof}

Using the isomorphism $J_{2n}: \Gamma_{2n} \to F_{2n}^o \times T$, we now transfer the groupoid structure on $\Gamma_{2n}$ to one on
$F_{2n}^o \times T$. 
The following \leref{:oid-FoT} is straightforward to check.

\ble{:oid-FoT}
Under the isomorphism $J_{2n}: \Gamma_{2n} \to F_{2n}^o \times T$, the Lie groupoid $\Gamma_{2n} \rightrightarrows F_n$ becomes the following 
Lie groupoid $F_{2n}^o \times T \rightrightarrows F_n$:
\begin{align*}
&\mbox{source map} \;\theta_+: \; F_{2n}^o \times T\to F_n,\; ([g_1, \ldots, g_{2n}]_{\sF_{2n}}, t) \mapsto [g_1, \ldots, g_n]_{\sF_n};\\
 &\mbox{target map}\; \theta_-: \; F_{2n}^o \times T\to F_n,\; \\
&\hs \hs \hs \hs \hs \;\;\;\;  ([g_1, \ldots, g_{2n}]_{\sF_{2n}}, t) \mapsto 
[t^{-1}[g_1 g_2 \cdots g_{2n}]_{\geq 0}g_{2n}^{-1}, g_{2n-1}^{-1}, \ldots, g_{n+1}^{-1}]_{\sF_n};\\
&\mbox{unit map}\; \epsilon: \; F_n \to F_{2n}^o \times T, \; [g_1,  \ldots, g_{n}]_{\sF_{n}} \mapsto 
([g_1, \ldots, g_n, g_n^{-1}, \ldots, g_1^{-1}]_{\sF_{2n}}, e);\\
&\mbox{inverse map}\; {\rm inv}: F_{2n}^o \times T \to F_{2n}^o \times T: \;\; \\
&\hs \hs \hs \hs \hs \;\;\;\;
([g_1, \ldots, g_{2n}]_{\sF_{2n}}, t) \mapsto (t^{-1}[g_1 g_2 \cdots g_{2n}]_{\geq 0} g_{2n}^{-1}, g_{2n-1}^{-1}, \ldots, g_{1}^{-1}]_{\sF_{2n}}, \; t^{-1});\\
&\mbox{multiplication:}\; \mbox{for}\; \gamma = ([g_1, \ldots, g_{2n}]_{\sF_{2n}}, t)\; \mbox{and} \;
\gamma' = ([g_1^\prime,  \ldots, g_{2n}^\prime]_{\sF_{2n}}, t') \; \mbox{with}\; \theta_-(\gamma) = \theta_+(\gamma'),\\
&\hs \hs  \gamma \gamma' =
([g_1, \;\ldots, \;g_n, \;g_{n+1}\cdots g_{2n} [g_1g_2 \cdots g_{2n}]_{\geq 0}^{-1}g_1^\prime \cdots g_n^\prime g_{n+1}^\prime, 
\; g_{n+2}^\prime, \, \ldots, \,g_{2n}^\prime]_{\sF_{2n}}, \, tt').
\end{align*}
\ele

By construction, we have now made  $(F_{2n}^o \times T, \pi_{2n} \bowtie 0)$ into a Poisson groupoid over $(F_n, \pi_n)$, as stated 
in the following  companion of \thmref{:tF2n-00} and \thmref{:C2n}.

\bthm{:F2nT}
For any integer $n \geq 1$, $(F_{2n}^o \times T, \pi_{2n} \bowtie 0)$ with the groupoid structure as above  is a Poisson groupoid over $(F_n, \pi_n)$.
\ethm

\section{Special configuration Poisson groupoids of flags and isomorphic models}\label{s:Gauu} 
\subsection{$T$-leaves}\label{ss:T-leaves0}
In view of the three isomorphic models of the total configuration Poisson groupoids given in $\S$\ref{ss:t-F2n} -
$\S$\ref{ss:F2no}, 
we now look at 
the $T$-leaves in  $(\Gamma_{n}, \tpi_{n})$, $(\calC_{n}, \tpi_{n})$, and $(F_{n}^o \times T, \pi_{n} \bowtie 0)$ for any integer $n \geq 1$. 

Recall again that $T$ acts on $F_n$ and $\tF_n$ via \eqref{eq:T-Fn}  and \eqref{eq:T-tFn} respectively 
and on $\calC_n = \calB^{n-1} \times \calA^o$ and on $F_n^o \times T$ diagonally.
Recall also the $T$-equivariant Poisson isomorphisms 
\[
\hTheta_n: \; (\Gamma_n, \tpi_n) \longrightarrow (\calC_n, \hpi_n) \hs \mbox{and} \hs
J_n: \; (\Gamma_n, \tpi_n) \longrightarrow (F_n^o \times T, \pi_n \bowtie 0),
\]
respectively given in \eqref{eq:hThetan} and \eqref{eq:Jn}.
For $\bfw = (w_1, \ldots, w_n) \in W^n$, introduce
\begin{equation}\label{eq:Ga-w}
\Gamma^{\bfw} = \Gamma_n \cap (B\bfw B) \subset \Gamma_n, \hs 
 \calC^{\bfw} = \hTheta_n(\Gamma^{\bfw}) \subset \calC_n,  \hs \mbox{and} \hs \O^{\bfw}_e \times T \subset F_n^o \times T,
\end{equation}
where recall that 
$B \bfw B$ is the image of $Bw_1B \times \cdots \times Bw_n B$ in $\tF_n$ and  that (see \eqref{eq:Owe-intro}) 
\begin{equation}\label{eq:Owe}
\O^\bfw_e = \{[g_1, g_2, \ldots, g_n]_{\sF_n} \in \O^\bfw: \; g_1g_2 \cdots g_n \in B_-B/B\} \subset \O^\bfw.
\end{equation}
By the definition of  $J_n$, one has $J_n(\Gamma^\bfw) = \O_e^\bfw \times T$ for each $\bfw \in W^n$.
It follows from the decomposition $\tF_n = \bigsqcup_{\bfw \in W^n} B \bfw B$ that
one has the disjoint unions
\begin{equation}\label{eq:T-leaves-all}
\Gamma_n = \bigsqcup_{\bfw \in W^n} \Gamma^\bfw, \hs 
\calC_n = \bigsqcup_{\bfw \in W^n} \calC^{\bfw}, \hs
F_n^o \times T = \bigsqcup_{\bfw \in W^n} (\O^{\bfw}_e \times T).
\end{equation}

\bpr{:T-leaves-all}
For any integer $n \geq 1$, the decompositions in \eqref{eq:T-leaves-all} are that of the 
$T$-leaves  of $(\Gamma_{n}, \tpi_{n})$, $(\calC_{n}, \tpi_{n})$, and $(F_{n}^o \times T, \pi_{n} \bowtie 0)$. In particular,
for any $\bfw \in W^n$, 
\[
\hTheta_n: \;  (\Gamma^\bfw, \tpi_n) \longrightarrow (\calC^\bfw, \hpi_n) \hs \mbox{and} \hs
J_n: \;(\Gamma^\bfw, \tpi_n) \longrightarrow (\O^{\bfw}_e \times T, \pi_n \bowtie 0)
\]
are  $T$-equivariant isomorphisms of single $T$-leaves.
\epr

\begin{proof}
The statement for the $T$-leaf decomposition for $(\Gamma_n, \tpi_n)$ follows \prref{:tFn-0}. The rest of \prref{:T-leaves-all}
follows from the fact that both  $\hTheta_n: (\Gamma_n, \tpi_n) \rightarrow (\calC_n, \hpi_n)$ and 
 $J_n: (\Gamma_n, \tpi_n) \rightarrow (F_n^o \times T, \pi_n \bowtie 0)$ are
$T$-equivariant Poisson isomorphisms. We also note that 
by \cite[Theorem 1.1]{Lu-Mou:flags}, each $\O^\bfw_e$ is the unique open $T$-leaf of $(\O^\bfw, \pi_n)$, so the fact that the
$T$-leaves of $(F_n^o \times T, \pi_n \bowtie 0)$ are precisely of the form $\O^\bfw_e \times T$ for $\bfw \in W^n$ also follows from \leref{:T-leaves}.
\end{proof}

We give a description of $\calC^{\bfw}$. 
Let $G_{\rm diag}$ be the diagonal of $G \times G$.
Recall that for $(f_1, f_2) \in \calB^2$, the Tits distance from $f_1$ to $f_2$ is defined to be 
the unique 
$w \in W$ such that $(f_1, f_2) \in G_{\rm diag} (e_\cdot B, \, w_\cdot B)$,
and we write $f_1 \stackrel{w}{\lrw} f_2$. 
In particular, $f_1$ and $f_2$ are said to be in general position if $f_1 \stackrel{w_0}{\lrw} f_2$, where $w_0$ is the longest element in $W$.
Set
\[
f_0 = e_\cdot B \in \calB \hs \mbox{and} \hs f_- = w_0{}_\cdot B \in \calB.
\]
For $f = b_-{}_\cdot N \in \calA^o$, where $b_- \in B_-$, let $\bar{f} = b_-{}_\cdot B \in \calB$.
 The following statement is now clear from the definition of the isomorphism $\hTheta: \Gamma_n \to \calC_n$.

\ble{le:Cw}
For $\bfw = (w_1, w_2, \ldots, w_n) \in W^n$, the sub-variety $\calC^{\bfw}$ of $\calC_n  = \calB^{n-1} \times \calA^o$ 
consists of all $(f_1, f_2, \ldots, f_{n-1}, f_n) \in \calB^{n-1} \times \calA^o$ such that
\[
f_0 \stackrel{w_1}{\lrw} f_1  \stackrel{w_2}{\lrw} f_2 \lrw
\cdots \stackrel{w_n}{\lrw} \bar{f}_n \stackrel{w_0}{\lrw} f_-.
\]
\ele

\bre{:Cw}
{\rm
For ${\bf w} = (w_1, \ldots, w_n) \in W^n$, choose
 any a reduced decomposition $w_i = s_{\alpha_{i, 1}} \cdots s_{\alpha_{i,l_i}}$ for each $1 \leq i \leq n$, where $l_i = l(w_i)$, and let
\[
b = (s_{\alpha_{1, 1}}, \,  \ldots, \, s_{\alpha_{1,l_1}}, \, s_{\alpha_{2, 1}}, \,  \ldots, \, s_{\alpha_{2,l_2}},\; \ldots, \;
s_{\alpha_{n, 1}}, \,  \ldots, \, s_{\alpha_{n,l_n}}) \in W^{l(\bfw)}.
\]
One then has $\calC^\bfw \cong \calC^{b} \subset \calC_{l(\bfw)}$. In the notation of \cite{Shen-Weng:BS},
$\calC^{b} = {\rm Conf}_e^{{b}}(\calA)$.
\hfill $\diamond$
}
\ere

\subsection{The Poisson groupoid $(\Gamma^{(\bfu, \bfu^{-1})}, \tpi_{2n}) \rightrightarrows (\O^\bfu, \pi_n)$ and two isomorphic models}\label{ss:Gauu-2}
Let now $\bfu = (u_1, \ldots, u_n) \in W^n$ and let $\bfu^{-1} = (u_n^{-1}, \ldots, u_1^{-1})$. 
Consider
\[
\Gauu = \Gamma_{2n} \cap (B (\bfu, \bfu^{-1}) B).
\]
By \prref{:T-leaves-all}, $\Gauu$ is a single $T$-leaf of $(\Gamma_{2n}, \tpi_{2n})$. 
One checks directly from the definitions that $\Gauu$ is the full sub-groupoid (see \deref{:wide-full}) of 
$\Gamma_{2n} \rightrightarrows F_n$ over $\Ou$, i.e., 
\[
\Gauu = \theta_{+}^{-1}(\Ou) \cap \theta_-^{-1}(\Ou),
\]
where $\theta_+, \theta_-: \Gamma_{2n} \to F_n$ are the source and target maps
of the groupoid $(\Gamma_{2n}, \tpi_{2n}) \rightrightarrows (F_n, \pi_n)$.

\bthm{:Gauu}
For any $\bfu = (u_1, u_2, \ldots, u_n) 
 \in W^n$, $(\Gauu, \tpi_{2n}) \rightrightarrows (\Ou, \pi_n)$ is a Poisson sub-groupoid of $(\Gamma_{2n}, \tpi_{2n}) \rightrightarrows (F_n, \pi_n)$.
\ethm

\begin{proof}
Being a $T$-leaf of $(\Gamma_{2n}, \tpi_{2n})$, $\Gauu$ is a Poisson sub-manifold of $(\Gamma_{2n}, \tpi_{2n})$. 
 It remains to show that $\theta_+|_{\Gauu}: \Gauu \to \Ou$ is a surjective submersion.
 Since $\Gauu$ contains the image of $\Ou \subset F_n$ under the unit map of $\Gamma_{2n} \rightrightarrows F_n$,
$\theta_+|_{\Gauu}: \Gauu \to \Ou$ is surjective. Define 
\[
\theta:\;\;\Ouuinv  \lrw \Ou, \;\; [g_1, g_2, \ldots, g_{2n}]_{\sF_{2n}} \longmapsto [g_1, g_2, \ldots, g_n]_{\sF_n}.
\]
Clearly $\theta$ is a submersion. 
Under the isomorphism $J_{2n}: \Gauu \to \Ouuinve \times T$, the map $\theta_+|_{\Gauu}: \Gauu \to \Ou$ becomes the projection
$\theta_+^\prime: \Ouuinve \times T \to \Ou, (q, t) \to \theta(q)$. 
Since $\Ouuinve$ is open in $\Ouuinv$, 
$\theta_+^\prime$ is a submersion. Thus $\theta_+|_{\Gauu}$ is a submersion.
\end{proof}
 
Recall from \eqref{eq:Thetan} the Poisson isomorphism $\Theta_n: (F_n, \pi_n) \to (\calB^n, \opi_n)$. Set
\[
\calB^\bfu = \Theta_n(\Ou) = \{(f_1, \ldots, f_n) \in \calB^n: f_0 \stackrel{u_1}{\lrw} f_1  \stackrel{u_2}{\lrw} f_2 \stackrel{u_3}{\lrw}
\cdots \stackrel{u_n}{\lrw} f_n\},
\]
so $\calB^\bfu$ is a Poisson sub-manifold of $(\calB^n, \opi_n)$. We have the following immediate consequence of \thmref{:Gauu}.

\bco{:Cuu-sub}
For any $n \geq 1$ and $\bfu \in W^n$, $(\calC^{(\bfu, \bfu^{-1})}, \hpi_{2n}) \rightrightarrows (\calB^\bfu, \opi_n)$ is a Poisson sub-groupoid
of $(\calC_{2n}, \hpi_{2n}) \rightrightarrows (\calB^n, \opi_n)$ in 
\thmref{:C2n} and is isomorphic, via the isomorphisms $\hTheta_{2n}: \Gauu \to \calC^{(\bfu, \bfu^{-1})}$ and $\Theta_n: \Ou \to \calB^\bfu$, to the Poisson groupoid 
$(\Gauu, \tpi_{2n}) \rightrightarrows (\Ou, \pi_n)$.
\eco

\bco{:Fuu-sub}
For any $n \geq 1$ and $\bfu \in W^n$, 
$(\O^{(\bfu, \bfu^{-1})}_e \times T, \, \pi_{2n} \bowtie 0) \rightrightarrows (\Ou, \pi_n)$ is a Poisson sub-groupoid of
$(F^o_{2n} \times T, \pi_{2n}\bowtie 0) \rightrightarrows (F^n, \pi_n)$ in 
\thmref{:F2nT} and is isomorphic, via the isomorphism $J_{2n}: \Gauu \to \O_e^{(\bfu, \bfu^{-1})}$, to the Poisson groupoid 
$(\Gauu, \tpi_{2n}) \rightrightarrows (\Ou, \pi_n)$.
\eco

\bde{:special}
{\rm
For $\bfu \in W^n$, we refer to either of the three isomorphic Poisson groupoids
\[
(\Gauu, \tpi_{2n}) \rightrightarrows (\Ou, \pi_n), \;\;\;\; (\calC^{(\bfu, \bfu^{-1})}, \hpi_{2n}) \rightrightarrows (\calB^\bfu, \opi_n),  \;\;\;\; 
(\O^{(\bfu, \bfu^{-1})}_e \times T, \, \pi_{2n} \bowtie 0) \rightrightarrows (\Ou, \pi_n)
\]
as a {\it special configuration Poisson groupoid} (of flags of $G$).
\hfill $\diamond$
}
\ede

The special configuration Poisson groupoids have a simple set-theoretical  description. We first set up some notation.
Recall that $N_G(T) \subset G$ is the normalizer  of $T$ in $G$.  
For $u \in W$, let $uT \subset N_G(T)$ be the set of representatives of $u$ in $N_G(T)$.
For $\du \in uT$, set
\begin{equation}\label{eq:Cu-Du}
C_{\du} = N \du \cap \du N_- 
\end{equation}
It is well-known (see, for example, \cite[Proposition 2.9]{FZ:double}) that the maps
\begin{equation}\label{eq:CBBC}
C_{\du} \times B \lrw BuB, \, (c, b) \Maps cb\hs \mbox{and} \hs B \times C_{\du} \lrw Bu^{-1}B, \, (b, c) \Maps bc^{-1},
\end{equation}
are both isomorphisms.
For $\bfu = (u_1, \ldots, u_n) \in W^n$, let $\bfu T^n = u_1T \times  \cdots \times  u_nT \subset N_G(T)^n$ and call any
$\dbfu = (\du_1, \ldots, \du_n) \in \bfu T^n$ a representative of $\bfu$. For $\dbfu = (\du_1, \ldots, \du_n) \in \bfu T^n$, set
\begin{equation}\label{eq:CD-bfu}
C_{\dbfu} = C_{\du_1} \times C_{\du_2} \times \cdots \times C_{\du_n}.
\end{equation}
One then has the isomorphisms
\begin{align}\label{eq:CuB}
&C_{\dbfu} \times B \lrw  B\bfu B,\;\; (c_1, \,c_2, \,\ldots, \,c_n, \, b) \Maps [c_1,  \cdots c_{n-1}, c_n b]_{\stF_n},\\
\label{eq:C-Ou}
&C_{\dbfu} \lrw \Ou = B \bfu B/B, \;\; (c_1, c_2, \ldots, c_n) \longmapsto [c_1, c_2, \ldots, c_n]_{\sF_n},\\
\label{eq:BDu} 
&B \times C_{\dbfu} \lrw B\bfu^{-1} B, \;\; (b, c_1^\prime, \ldots, c_{n-1}^\prime, c_n^\prime) \Maps [b(c_n^\prime)^{-1}, \,
(c_{n-1}^\prime)^{-1}, \, \ldots,  \, (c_1^\prime)^{-1}]_{\stF_n}.
\end{align}
For $c = (c_1, c_2, \ldots, c_n) \in C_{\dbfu}$, set $[c]_{\sF_n} = [c_1, \ldots, c_n]_{\sF_n}$ and $\underline{c} = c_1c_2 \cdots c_n \in G$. 
Define
\begin{equation}\label{eq:Gdbu}
\Gdbu =\{(c, b, b_-, c'):  c, c' \in C_{\dbfu},\, b \in B, \, b_- \in B_-, \,  \underline{c} \, b = b_- \underline{c'}\} \subset
 C_{\dbfu} \times B \times B_- \times C_{\dbfu}.
\end{equation}
We then have the isomorphism ${\calI}_{\dbu}:\Gdbu \to \Gauu$ given by 
\[
{\calI}_{\dbu}((c_1, \ldots, c_n), \,b, \,b_-, (c_1^\prime, \ldots, c_n^\prime)) =
[c_1, \ldots, c_{n-1}, c_nb, \,(c_n^\prime)^{-1}, (c_{n-1}^\prime)^{-1}, \ldots, (c_1^\prime)^{-1}]_{\stF_{2n}}.
\]

\bpr{:GaSig-uu}
For any $\bfu \in W^n$ and $\dbu \in \bfu T^n$,  ${\calI}_{\dbu}$ is an isomorphism of Lie groupoids over $\O^\bfu$, 
where the groupoid $\Gauu \rightrightarrows \Ou$ is given as in \thmref{:Gauu}, and 
$\Gdbu \rightrightarrows \Ou$ is the groupoid defined in \cite{Mou:Guu} as follows:
\begin{align*}
&\mbox{source map} \; \theta_+: \; \Gdbu\longrightarrow \Ou, 
\;\; (c, \,b,\, b_-, \,c') \longmapsto [c]_{\sF_n},\\
&\mbox{target map} \;\theta_-:\;  \Gdbu\longrightarrow \Ou, 
\; \; (c, \,b,\, b_-, \,c') \longmapsto [c']_{\sF_n},\\
&\mbox{unit map} \;\epsilon: \; \Ou \longrightarrow \Gdbu, \;\;
[c]_{\sF_n} \longmapsto (c,\, e, \,e,\, c),\\
&\mbox{inverse map} \; {\rm inv}: \; \Gdbu\longrightarrow \Gdbu,
\;\; (c, \,b,\, b_-, \,c') \longmapsto (c', \,b^{-1},\, b_-^{-1}, \,c) ,\\
&\mbox{multiplication} :  (c, \,b,\, b_-, \,c')  (c', \,b',\, b_-^\prime, \,c^{\prime\prime}) =
(c, \,bb',\, b_-b_-^\prime, \,c^{\prime\prime}).
\end{align*}
\epr

\begin{proof} By \eqref{eq:CuB} and \eqref{eq:BDu}, 
${\calI}_{\dbu}$ is an isomorphism of varieties. The fact that ${\calI}_{\dbu}$ is an isomorphism from 
the groupoid $\Gdbu \rightrightarrows \Ou$ to $\Gauu \rightrightarrows \Ou$ follows directly from the
definitions of the two groupoids.
\end{proof}

\bre{:sub-prod}
Under the isomorphism $C_{\dbfu} \to \Ou, c \mapsto [c]_{\sF_n}$, one can regard $\Gdbu \rightrightarrows \Ou$ as a
groupoid $\Gdbu \rightrightarrows C_{\dbfu}$, and as a such, it is a sub-groupoid of
the direct product
\[
C_{\dbfu}  \times B \times B_- \times C_{\dbfu}  \rightrightarrows C_{\dbfu} 
\]
 of two groupoids: the pair groupoid
$C_{\dbfu}  \times C_{\dbfu}  \rightrightarrows C_{\dbfu} $ and the direct product group $B \times B_-$
as a groupoid over the one point space, where we identify
\[
 C_{\dbfu}  \times B \times B_- \times C_{\dbfu}  \longrightarrow C_{\dbfu}  \times C_{\dbfu}  \times B \times B_-, \;\;
(c, \,b,\, b_-, \, c') \longmapsto 
(c, \, c',  \,b,\, b_-).
\]
\ere

\subsection{Generalized double Bruhat cells as Poisson groupoids}\label{ss:Guu} We now give the third isomorphic model of the special
configuration Poisson groupoids using {\it generalized double Bruhat cells}.
For $n \geq 1$,  recall the right action of $G^n$ on itself given in \eqref{eq:Gn-Bn}.
Let 
\[
\tF_{-n} = G \times_{B_-} \cdots \times_{B_-} G,
\]
the quotient of 
$G^n$ by $B_-^{n-1} \times \{e\} \subset G^n$, and let $\tpi_{-n}$ be the Poisson structure on $\tF_{-n}$
that is the (well-defined) projection of the Poisson structure $(\pist)^n$ on $G^n$. For 
$(g_1, g_2, \ldots, g_n) \in G^n$, denote its image  in $\tF_{-n}$ as
$[g_1, \ldots, g_n]_{\stF_{-n}}$. For any integers $m, n \geq 1$, we introduce in $\S$\ref{ss:Guv-app} a Poisson 
structure $\tpi_{m, n}$ on $\tF_m \times \tF_{-n}$, which is the sum of the product Poisson structure
$(\tpi_m, \tpi_{-n})$ with a certain {\it mixed term}. See \deref{:tpimn}. The Poisson structure $\tpi_{m, n}$ is $T$-invariant 
under the $T$-action on $\tF_m \times \tF_{-n}$ 
given by
\begin{equation}\label{eq:T-tFmn}
t \cdot ([g_1, g_2, \ldots, g_m]_{\stF_m}, [k_1, k_2, \ldots, k_n]_{\stF_{-n}}) =
([tg_1, g_2, \ldots, g_m]_{\stF_m}, [tk_1, k_2, \ldots, k_n]_{\stF_{-n}}).
\end{equation}
Introduce the sub-manifold 
\begin{equation}\label{eq:Gmn-0}
G_{m, n}= \{([g_1, g_2, \ldots, g_m]_{\stF_m}, [k_1, k_2, \ldots, k_n]_{\stF_{-n}}): g_1g_2 \cdots g_n = k_1k_2 \cdots k_n\} \subset \tF_m \times \tF_{-n}.
\end{equation}
For $\bfu =(u_1, \ldots, u_m)\in W^m$ and $\bfv =(v_1, \ldots, v_n) \in W^n$, let again $B \bfu B$ be the image of 
$(Bu_1B) \times \cdots \times (Bu_m B)$ in $\tF_m$,  let 
$B_- \bfv B_- \subset \tF_{-n}$ be
 the image of $(B_-v_1B_-) \times \cdots \times (B_-v_n B_-)$ in $\tF_{-n}$, and set
\begin{equation}\label{eq:Guv-0}
\Guv = G_{m, n} \cap (B \bfu B \times B_- \bfv B_-) \subset G_{m, n}.
\end{equation}
It follows from the Bruhat decomposition of $G$ that one has the disjoint union
\begin{equation}\label{eq:Gmn-Guv}
G_{m,n} = \bigsqcup_{(\bfu, \bfv) \in W^m \times W^n} \Guv.
\end{equation}
We prove in \coref{:Tleaf-Gmn} that $G_{m, n}$ is a $T$-invariant Poisson sub-manifold of $(\tF_m \times \tF_{-n}, \tpi_{m, n})$, and that 
\eqref{eq:Gmn-Guv} is the decomposition of $(G_{m, n}, \tpi_{m, n})$ into its $T$-leaves.  Generalizing the case of 
\[
(G_{1, 1}, \tpi_{1, 1}) \cong (G, \pist),
\]
we call $\Guv$, for any $\bfu \in W^m$ and $\bfv \in W^n$, a {\it generalized double Bruhat cell}.
When $m = n$, the Poisson manifold $(\tF_n \times \tF_{-n}, \tpi_{n, n})$, as well as the generalized double Bruhat cells $\Guv$ for $\bfu, \bfv \in W^n$,
were introduced in \cite[$\S$1.4]{Lu-Mou:flags}.  
We also note that for any $m, n$, the projections $(\tF_m \times \tF_{-n},  \tpi_{m, n}) \to (\tF_m, \tpi_m)$ and
$(\tF_m \times \tF_{-n},  \tpi_{m, n}) \to (\tF_{-n}, \tpi_{-n})$ to the factors are
Poisson. It follows in particular that for any $m \geq 1$, the map
\begin{equation}\label{eq:Gm1-tFm}
(G_{m, 1}, \; \tpi_{m, 1}) \longrightarrow (\tF_m, \tpi_m), \; ([g_1, g_1, \ldots, g_m]_{\stF_m}, g_1g_2 \cdots g_m) \longmapsto 
[g_1, g_1, \ldots, g_m]_{\stF_m},
\end{equation}
is a Poisson isomorphism.

Let now $\bfu = (u_1, \ldots, u_n) \in W^n$, and consider the generalized double Bruhat cell $\Guu$.  Fix any 
$\dbfu = (\du_1, \ldots, \du_n) \in \bfu T^n$. As $B_-uB_- = B_-  C_{\du}$ is a direct product decomposition 
for any $u \in W$ and any $\du \in uT$, one has 
the isomorphism
\begin{equation}\label{eq:BBv0}
B_- \times C_{\dbfu} \lrw B_- \bfu B_-, \;\; (b_-, \,c_1^\prime, \,c_2^\prime, \,\ldots, \,c_n^\prime) \Maps [b_-c_1^\prime, \, c_2^\prime, \, 
\ldots, \, c_n^\prime]_{\stF_{-n}}.
\end{equation}
Recalling the definition of $\Gdbu \subset C_{\dbu} \times B \times B_- \times C_{\dbu}$, 
one then has the isomorphism
\begin{equation}\label{eq:calJu}
\calJ_{\dbu}: \; \Gdbu \longrightarrow \Guu, \;
(c, b, b_-, c') \longmapsto ([c_1, \ldots, c_{n-1}, c_nb]_{\stF_n}, [b_-c_1^\prime, c_2^\prime, \ldots, c_n^\prime]_{\stF_{-n}}),
\end{equation}
where $c = (c_1, \ldots, c_{n-1}, c_n)$, $c' = (c_1^\prime, c_2^\prime, \ldots, c_n^\prime)$, and $(c, b, b_-, c') \in \Gdbu$. 

\bde{:Guu-oid} \cite{Mou:Guu} For each choice $\dbfu \in \bfu T^n$, define 
the groupoid $\Guu\rightrightarrows \Ou$ such that 
\[
\calJ_{\dbu}: \; \Gdbu \longrightarrow \Guu
\]
 is a groupoid isomorphism 
from the groupoid $\Gdbu \rightrightarrows \Ou$ in \prref{:GaSig-uu} 
to $\Guu \rightrightarrows \Ou$.
\ede

Using the parametrization $\calJ_{\dbu}: \Gdbu \to \Guu$, we can now define
\begin{align}\label{eq:Edu}
E_{\dbfu}: \;\Guu \longrightarrow  \Gamma_{2n}, \;\;&([c_1, \ldots, c_{n-1}, c_nb]_{\stF_n}, [b_- c_1^\prime, c_2^\prime, \ldots, c_n^\prime]_{\stF_{-n}})\\
\nonumber
& \longmapsto 
[c_1, \ldots, c_{n-1}, c_nb, \; (c_n^\prime)^{-1}, \ldots, (c_2^\prime)^{-1}, (c_1^\prime)^{-1}]_{\tF_{2n}},
\end{align}
where $((c_1, \ldots, c_n),  b, b_-, (c_1^\prime, \ldots, c_n^\prime)) \in \Gdbu$.

\bthm{:Guu-gpoid} For each $\dbfu \in \bfu T^n$, 
the map $E_{\dbfu}$ is a $T$-equivariant Poisson embedding of $(\Guu, \tpi_{n,n})$ into $(\Gamma_{2n}, \tpi_{2n})$
and gives a Poisson isomorphism
\[
E_{\dbfu}:\;\; (\Guu, \;\tpi_{n, n}) \longrightarrow (\Gauu, \;\tpi_{2n}) \subset (\Gamma_{2n}, \tpi_{2n}). 
\]
Consequently, with the groupoid structure defined as in \deref{:Guu-oid}, 
$(\Guu, \tpi_{n, n}) \rightrightarrows (\Ou, \pi_n)$ is a Poisson groupoid,
and $E_{\dbfu}$ is an isomorphism from the Poisson groupoid $(\Guu, \tpi_{n, n}) \rightrightarrows (\Ou, \pi_n)$ to the 
Poisson 
groupoid
$(\Gauu, \tpi_{2n}) \rightrightarrows (\Ou, \pi_n)$.
\ethm

\begin{proof}
The fact that $E_{\dbfu}$ is a Poisson isomorphism onto $(\Gauu, \tpi_{2n})$ is a special case of \coref{:Guv-emb}
by taking $\bfu = \bfv$ and $E_{\dbfu} = E_{\bfu, \dbfu}$. That $E_{\dbfu}$ is a groupoid isomorphism
follows from the definition of the groupoid $\Guu \rightrightarrows \Ou$ and \prref{:GaSig-uu}.
\end{proof}

\bre{re:Guu-du}
{\rm 
 We emphasize that the
groupoid structure on $\Guu$ over $\Ou$ depends on the choice of the representative $\dbfu \in \bfu T^n$. If $\hat{\bfu} =(\hat{u}_1, \ldots, 
\hat{u}_n) \in \bfu T^n$ is another
such choice, and if $t \in T$ is such that $\hat{u}_1 \cdots \hat{u}_n = t \dot{u}_1 \cdots \dot{u}_n$, then 
$E_{\dot{\bfu}} = E_{\hat{\bfu}} \circ r_t$, where $r_t: \Guu \to \Guu$ is given by 
\[
r_t ([g_1, \ldots, g_{n-1}, g_n]_{\stF_m}, [k_1, \ldots, k_{n-1}, k_n]_{\stF_{-n}}) =
([g_1, \ldots, g_{n-1}, g_n t]_{\stF_m}, [k_1, \ldots, k_{n-1}, k_n t]_{\stF_{-n}}) ,
\]
Thus $r_t: \Guu \to \Guu$ defines a groupoid isomorphism from the groupoid structure on $\Guu$ defined by $\dot{\bfu}$ to that defined by 
$\hat{\bfu}$.
\hfill $\diamond$
}
\ere

\sectionnew{Configuration symplectic groupoids of flags}\label{s:leaves}
In this section, we assume that $G$ is connected and simply connected. 
Symplectic leaves of $(\O^\bfw_e \times T, \pi_n \bowtie 0)$, for any $n\geq 1$ and $\bfw \in W^n$, are determined in 
$\S$\ref{s:app-leaves}, where we also give a complete description of
the symplectic leaves of all the three series
\[
(F_n^o \times T, \,\pi_n \bowtie 0), \hs (G_{m+n}, \,\tpi_{m, n}), \hs (\tF_n, \,\tpi_n), \hs m, n \geq 1,
\]
 of $T$-Poisson varieties. See $\S$\ref{s:app-leaves} for details. 
For any $n \geq 1$ and $\bfu \in W^n$, we show in this section that all the units of 
the groupoid  $\Gauu \rightrightarrows \O^\bfu$ are contained in a single symplectic leaf, denoted as 
$\Lambda^{(\bfu, \bfu^{-1})}$, of 
$(\Gauu, \tpi_{2n})$, and that $\Lambda^{(\bfu, \bfu^{-1})}$ 
 is a Lie sub-groupoid
of $(\Gauu, \tpi_{2n}) \rightrightarrows (\O^\bfu, \pi_n)$, obtaining thus a symplectic groupoid 
\[
(\Lambda^{(\bfu, \bfu^{-1})}, \tpi_{2n}) \rightrightarrows (\O^\bfu, \pi_n). 
\]
Using \thmref{:Guu-gpoid}, we then show that every symplectic 
leaf of $(\Guu, \tpi_{n, n})$ is symplectic groupoid over $(\O^\bfu, \pi_n)$. 

\subsection{The symplectic groupoid $(\Lambda^{(\bfu, \bfu^{-1})}, \tpi_{2n}) \rightrightarrows (\O^\bfu, \pi_n)$}\label{ss:Lambda-poid}
Assume that $G$ is connected and simply connected. Let $X^*(T)$ be the character lattice of $T$. For $\lambda \in X^*(T)$ and 
$t \in T$, write $t^\lambda$ for the value of $\lambda$ at $t$. Let $\Phi_0 \subset
X^*(T)$ be the set of simple roots determined by the choice of $B$, and let 
$\{\omega_\alpha: \alpha \in \Phi_0\} \subset X^*(T)$ be the corresponding set of fundamental weights.
For $u \in W$ and $t \in T$, we also set 
$t^u = \dot{u}^{-1} t \dot{u}$ using any $\dot{u} \in u T$. 

Let now $\bfu = (u_1, \ldots, u_n) \in W^n$ be arbitrary. Define 
\[
\Tbfu = \{t \in T: t^{\omega_\alpha} = 1, \, \forall \, \alpha \in {\rm supp}^o(\bfu)\},
\]
where ${\rm supp}^o(\bfu) = \{\alpha \in \Phi_0: u_i \omega_\alpha = \omega_\alpha, \, \forall \, i \in [1, n]\}$. 
Choose any $\dbfu =(\du_1, \ldots, \du_n)\in \bfu T^n$, and let again
$C_{\dbfu} = C_{\du_1} \times \cdots \times C_{\du_n}$.  Let $u = u_1u_2 \cdots u_n \in W$.
Using the isomorphism ${\calI}_{\dbu}:\Gdbu \to \Gauu$ to write an element in 
$\Gauu$ uniquely as 
\begin{equation}\label{eq:gamma-1}
\gamma = [c_1, \ldots, c_{n-1}, c_nb, \,(c_n^\prime)^{-1}, (c_{n-1}^\prime)^{-1}, \ldots, (c_1^\prime)^{-1}]_{\stF_{2n}},
\end{equation}
where $ (c_1, \ldots, c_n),  
(c_1^\prime, \ldots,  c_n^\prime) \in C_{\dbfu}$, $b \in B, b_- \in B_-$, and $c_1 \cdots c_{n-1} c_n \, b = b_- c_1^\prime c_2^\prime \cdots c_n^\prime$,
we define the sub-variety $\Lambda^{(\bfu, \bfu^{-1})}$ of $\Gauu$ by 
\[
\Lambda^{(\bfu, \bfu^{-1})}=\{\gamma \in \Gauu\, \mbox{as in \eqref{eq:gamma-1}}: \; [b]_0 [b_-]_0^u = e, \, [b]_0 \in \widetilde{T}^{\bfu}\}.
\]
Note that $\Lambda^{(\bfu, \bfu^{-1})}$ contains all the units of the groupoid $\Gauu \rightrightarrows \O^\bfu$, and that
\begin{equation}\label{eq:Suu}
\calI_{\dbu}^{-1}(\Lambda^{(\bfu, \bfu^{-1})}) = \calS^{\dot{\bfu}, \dot{\bfu}} \stackrel{\rm def}{=}
\{(c, b, b_-, c') \in \Gdbu: \, [b]_0 [b_-]_0^u = e, \, [b]_0 \in \widetilde{T}^{\bfu}\}.
\end{equation}
Note that $\calS^{\dot{\bfu}, \dot{\bfu}}$ is a wide sub-groupoid of $\Gdbu \rightrightarrows \O^\bfu$.
It thus follows from \prref{:GaSig-uu} that $\Lambda^{(\bfu, \bfu^{-1})}$ is a wide sub-groupoid of $\Gauu \rightrightarrows \O^\bfu$.

\bthm{:Gauu-leaf}
For any $\bfu \in W^n$, $\Lambda^{(\bfu, \bfu^{-1})}$ is a symplectic leaf of $(\Gauu, \tpi_{2n})$
and is a Lie sub-groupoid of $\Gauu \rightrightarrows \O^\bfu$. 
Consequently, $(\Lambda^{(\bfu, \bfu^{-1})}, \tpi_{2n}) \rightrightarrows
(\O^\bfu, \pi_n)$ is a symplectic groupoid.
\ethm

\begin{proof} Symplectic leaves of $(\Gamma^\bfw, \tpi_n)$ for any $n \geq 1$ and any $\bfw \in W^n$
are described in \thmref{:Gaw-leaves}, and the case of $\bfw = (\bfu, \bfu^{-1}) \in W^{2n}$ is given in \exref{:uu-leaf}.
More concretely, for $\gamma \in \Gauu$ in \eqref{eq:gamma-1}, let
$(c_{n+1}, \ldots, c_{2n}) \in 
C_{\dot{u}_n^{-1}} \times \cdots \times C_{\dot{u}_1^{-1}}$ and $b_n, \ldots, b_1 \in B$ be such that
\[
b (c_n^\prime)^{-1} = c_{n+1} b_n, \;\; b_n (c_{n-1}^\prime)^{-1} = c_{n+2} b_{n-1}, \;\ldots, \;\;
b_2 (c_1^\prime)^{-1} = c_{2n}b_1.
\]
Then $\gamma =[c_1, \ldots, c_{n-1}, c_n, \, c_{n+1}, \ldots, c_{2n-1}, c_{2n}b_1]_{\stF_{2n}}$ and 
$[b_1]_0 = [b]_0^{u^{-1}}$. Let
\[
\dbfw = (\dot{u}_1, \ldots, \dot{u}_n, \dot{u}_n^{-1}, \ldots, \dot{u}_1^{-1}) \in \bfw T^{2n}.
\]
We then
have the alternative description of $ \Lambda^{(\bfu, \bfu^{-1})} \subset \Gauu$ as consisting of all 
\[
\gamma =[c_1, \ldots, c_{n-1}, c_n, \, c_{n+1}, \ldots, c_{2n-1}, c_{2n}b_1]_{\stF_{2n}} \in \Gauu,
\]
where $(c_1, \ldots, c_n, c_{n+1}, \ldots, c_{2n}) \in C_{\dbfw}$,  $b_1 \in B$, and $b_- =c_1, \cdots c_n c_{n+1} \cdots c_{2n}b_1 \in B_-$
 are such that
$[b_-]_0 \in \Tbfu$ and $[b_1]_0[b_-]_0 = e$.
By \exref{:uu-leaf}, $\Lambda^{(\bfu, \bfu^{-1})}$ is the symplectic leaf of $(\Gauu, \tpi_{2n})$ passing through the point
\[
[\bfu, \bfu^{-1}]_{\stF_{2n}} \, \stackrel{\rm def}{=}\, [\dot{u}_1, \ldots, \dot{u}_n, \, \dot{u}_n^{-1}, \ldots, 
\dot{u}_1^{-1}]_{\stF_{2n}} \in \Gauu.
\]

We already know that 
$\Lambda^{(\bfu, \bfu^{-1})}$ is a sub-groupoid of $\Gauu \rightrightarrows \O^\bfu$. To show
that $\Lambda^{(\bfu, \bfu^{-1})}$ is a Lie sub-groupoid of $\Gauu \rightrightarrows \O^\bfu$, it remains to show that 
the source map $\theta_+: \Gauu \to \O^\bfu$ restricts to a surjective submersion $\theta:=\theta_+|_{\Lambda^{(\bfu, \bfu^{-1})}}: 
\Lambda^{(\bfu, \bfu^{-1})} \to \O^\bfu$. To this end, let 
\[
\Sigma^{(\bfu, \bfu^{-1})} = J_{2n}(\Lambda^{(\bfu, \bfu^{-1})}) \subset \O^{(\bfu, \bfu^{-1})}_e \times T.
\]
 Then 
$\Sigma^{(\bfu, \bfu^{-1})}$ is the symplectic leaf of $\pi_{2n} \bowtie 0$ in $\O^{(\bfu, \bfu^{-1})}_e \times T$ through the point
\[
([\dot{u}_1, \ldots, \dot{u}_n, \, \dot{u}_n^{-1}, \ldots, 
\dot{u}_1^{-1}]_{\sF_{2n}}, \; e) \in\O^{(\bfu, \bfu^{-1})}_e \times T.
\]
By \exref{:w-e},  $\Sigma^{(\bfu, \bfu^{-1})} =\Sigma^{\dbfw}$, consisting of all 
$([c_1, \ldots, c_{2n}]_{\sF_{2n}}, t) \in \O^{(\bfu, \bfu^{-1})}_e \times T$, where $t \in \Tbfu$ and 
$(c_1, \ldots,  c_{2n}) \in C_{\dbfw}$ such that 
\[
c_1c_2  \cdots c_{2n} \in B_-B \hs \mbox{and} \hs t^2 =[c_1, c_2, \ldots, c_{2n}]_0.
\]
We also know from \exref{:w-e} that the projection 
$P: \Sigma^{(\bfu, \bfu^{-1})}\to \O^{(\bfu, \bfu^{-1})}_e ,  (q, t) \mapsto q,$
is a $2^{|{\rm supp}(\bfu)|}$-to-$1$ covering map, where
${\rm supp}(\bfu) = \Phi_0\backslash {\rm supp}^o(\bfu)$.
Let 
\[
p: \; \O^{(\bfu, \bfu^{-1})}_e \longrightarrow \O^\bfu, \; [g_1, \ldots, g_n, g_{n+1}, \ldots, g_{2n}]_{\sF_{2n}} \longmapsto 
[g_1, \ldots, g_n]_{\sF_{n}}.
\]
Since $p$ is a submersion, one sees that 
$\theta  = p \circ P \circ  J_{2n}:  \Lambda^{(\bfu, \bfu^{-1})} \to \O^\bfu$ is a submersion. 
\end{proof}

\subsection{Symplectic leaves in $\Guu$ as symplectic groupoids}\label{ss:Guu-sym-poids}
Let $\bfu  \in W^n$ and consider the Poisson manifold $(\Guu, \tpi_{n, n})$. As $(\Guu, \tpi_{n, n})$ is a single $T$-leaf,
every symplectic leaf of $(\Guu, \tpi_{n, n})$ passes through the point $(\dbfu, \dbfu) \in \Guu$ for some 
$\dbfu  \in \bfu T^n$. Let $\dbfu \in \bfu T^n$, and let $S^{\dbfu, \dbfu}$ be the symplectic leaf
of $(\Guu, \tpi_{n, n})$ through $(\dbfu, \dbfu)$. 
Using the parametrization $\calJ_{\dbu}:  \Gdbu \rightarrow \Guu$ in \eqref{eq:calJu}, we write every ${\bf g} \in \Guu$ uniquely as
\begin{equation}\label{eq:bfg-0}
{\bf g} = ([c_1, \ldots, c_{n-1}, c_nb]_{\stF_n}, [b_- c_1^\prime, c_2^\prime, \ldots, c_n^\prime]_{\stF_{-n}}).
\end{equation}
where $((c_1, \ldots, c_n),  b, b_-, (c_1^\prime, \ldots, c_n^\prime)) \in \Gdbu$. Recall from \eqref{eq:Suu} the wide sub-groupoid
 $\calS^{\dbfu, \dbfu} \rightrightarrows \O^\bfu$ of the groupoid $\Gdbu \rightrightarrows \O^\bfu$. 
Recall also from \thmref{:Guu-gpoid} the Poisson isomorphism $E_{\dbfu}: (\Guu, \tpi_{n, n}) \to (\Gauu, \tpi_{2n})$.  
The following \thmref{:Guu-leaf} is now a direct consequence of \thmref{:Guu-gpoid} and 
\thmref{:Gauu-leaf}.

\bthm{:Guu-leaf}
The symplectic leaf $S^{\dbfu, \dbfu}$ of $(\Guu, \tpi_{n, n})$ through $(\dbfu, \dbfu)$ is given by 
\begin{align*}
S^{\dbfu, \dbfu} = E_{\dbfu}^{-1} (\Lambda^{\bfu, \bfu^{-1}}) 
= \{{\bf g} \in \Guu \, \mbox{as in \eqref{eq:bfg-0}}: \; [b]_0 [b_-]_0^u = e, \, [b]_0 \in \widetilde{T}^{\bfu}\}.
\end{align*}
Equip $S^{\dbfu, \dbfu}$ with the structure of a Lie groupoid over $\O^\bfu$ via the isomorphism
$\calJ_{\dbfu}: \calS^{\dbfu, \dbfu} \to S^{\dbfu, \dbfu}$. Then 
$(S^{\dbfu, \dbfu}, \tpi_{n, n}) \rightrightarrows (\O^\bfu, \pi_n)$ is a symplectic groupoid.
\ethm

\appendix
\markboth{Appendix}{Appendix}
\numberwithin{equation}{section}

\section{Mixed product Poisson structures}\label{s:app-mix}
In $\S$\ref{ss:mixed}, we recall from \cite{Lu-Mou:mixed} the construction of mixed product Poisson structures
using Poisson Lie group actions. We then show in $\S$\ref{ss:Fn-mixed} that the Poisson varieties
$(\calC_n, \hpi_n)$ and $(\calB^n, \opi_n)$ defined in $\S$\ref{ss:C2n} are all mixed product Poisson varieties.

\subsection{Mixed product Poisson structures}\label{ss:mixed}
The following definition was introduced in \cite{Lu-Mou:mixed}.

\bde{:mixed-df}
{\rm
Given two manifolds $Y_1$ and $Y_2$, by a {\it mixed product Poisson structure} on the product manifold $Y_1 \times Y_2$ we mean a Poisson bivector field $\pi$ on $Y_1 \times 
Y_2$ that projects to well-defined Poisson structures on $Y_1$ and $Y_2$. Given $Y_i$, $1 \leq i \leq n$, where $n \geq 2$, a Poisson structure $\pi$ on the product manifold 
$Y =Y_1 \times \cdots \times Y_n$ is said to be a {\it mixed product} if the projection of $\pi$ to $Y_i \times Y_j$ is a well-defined mixed product Poisson structure on $Y_i 
\times Y_j$ for any $1 \leq i < j \leq n$, and in this case, we also call the pair
$(Y_1 \times \cdots \times Y_n, \pi)$ a {\it mixed product Poisson manifold}.
\hfill $\diamond$
}
\ede

Assume now that $((G, \piG), (G^*, \piGs))$ is any
pair of dual Poisson Lie groups, i.e. the Lie bialgebras $(\g, \delta_\g)$ and $(\g^*, \delta_{\g^*})$
are dual to each other, where recall that  $\delta_\g: \g \to \wedge^2 \g$ and $\delta_{\g^*}: \g^* \to \wedge^2 \g^*$
are respectively the linearizations of $\piG$ and $\piGs$ at the identity elements of $G$ and $G^*$.
Suppose that $(X, \piX)$ and $(Y, \piY)$ are two Poisson manifolds, with respective {\it right} and {\it left}
Poisson actions
\begin{equation}\label{eq:rho-lam}
\rho: \; (X, \piX) \times (G^*, \piGs) \lrw (X, \piX) \hs \mbox{and} \hs
\lambda: \; (G, \piG) \times (Y, \piY) \lrw (Y, \piY)
\end{equation}
by the Poisson Lie groups $(G^*, \pi_{\sG^*})$ and $(G, \pi_\sG)$. Let
$\rho: \g^*\to \mX^1(X)$ and $\lam: \g \to \mX^1(Y)$ be the corresponding Lie algebra homomorphism
and Lie algebra anti-homomorphism
 (see notation in $\S$\ref{ss:nota-00}). Then $\rho$ and $\lam$ are respectively right and left Lie bialgebra actions (see
\cite[$\S$2.5]{Lu-Mou:mixed}).
 Define the bi-vector field
\begin{equation}\label{eq:mixed}
\piX \times_{(\rho, \lam)} \piY = (\piX, 0) + (0, \piY) -\sum_{i=1}^n (\rho(\xi_i^*), \, 0) \wedge (0, \; \lam(\xi_i))
\end{equation}
on $X \times Y$,
where $\{\xi_i\}_{i=1, \ldots, n}$ is any basis for $\g$ and
$\{\xi_i^*\}_{i=1, \ldots, n}$
its dual basis for  $\g^*$.  Note that the definition of $\piX \times_{(\rho, \lam)} \piY$ only uses the Lie algebra actions $\rho$ and $\lam$.

\bld{:mixed} \cite[Proposition 4.2]{Lu-Mou:mixed} For any pair $(\rho, \lam)$ of Poisson Lie group actions in \eqref{eq:rho-lam},
the bi-vector field $\piX \times_{(\rho, \lam)} \piY$ on $X \times Y$ is  Poisson, and it
is called the {\it mixed product} of $\piX$ and $\piY$ associated to the pair $(\rho, \lam)$.
\eld

\bex{:XGY} Let $(\rho, \lam)$ be a pair of Poisson Lie group actions as in \eqref{eq:rho-lam}. Equipping $G^* \times Y$ with the
Poisson structure $\pi = \pi_{\sG^*} \times_{(\rho_{\sG^*}, \lam)} \piY$, where $\rho_{\sG^*}$ is the right action action of $G^*$ on itself by right translation,
one then has the left Poisson Lie group action
\[
\lam': \;\; (G^*, \;\pi_{\sG^*}) \times (G^* \times Y, \;\pi) \lrw( G^* \times Y, \; \pi), \; (a, \, (a', y)) \Maps (aa', y), \hs a, a' \in G^*, \, y \in Y.
\]
Let $X \times_{G^*} (G^* \times Y)$  is the quotient of $X \times (G^* \times Y)$ by the right diagonal action of $G^*$.
The Poisson structure $\piX \times \pi$ on 
 $X \times (G^* \times Y)$ then projects to a well-defined Poisson structure on  $X \times_{G^*} (G^* \times Y)$, 
which we denote as $\pi'$.
On the other hand, one has the
isomorphism
\[
\psi: \;\; X \times_{G^*} (G^* \times Y) \lrw X \times Y, \; [x, \, (a, y)] \Map (xa, \, y), \; \hs x \in X, \, a \in G^*, \, y \in Y.
\]
By an argument similar to that used in the proof of \cite[Lemma A.1]{Lu-Mi:Kostant}, one sees that
\[
\psi: (X \times_{G^*} (G^* \times Y), \; \pi') \lrw (X \times Y, \; \piX \times_{(\rho, \lam)} \piY)
\]
 is a Poisson map.
\hfill $\diamond$
\eex

Consider now the standard complex semi-simple Lie group $(G, \pist)$, and we keep the notation set up in $\S$\ref{ss:pist}. 
Define the non-degenerate bilinear pairing $\lara_{(\b, \b_-)}$ between $\b$ and $\b_-$ by
\begin{equation}\lb{eq:bb-pair}
\la x_- + x_0, \;\, y_+ + y_0\ra_{(\b_-, \b)} = \frac{1}{2}\la x_-,\; y_+\ra_\g + \la x_0, \;y_0\ra_\g, \hs x_- \in \n_-,\,
\, x_0, \,y_0 \in \h, \, y_+ \in \n_.
\end{equation}
For a basis $\{x_i\}$ for $\b_-$ and dual basis $\{x^i\}$ for $\b$ under $\lara_{(\b_-, \b)}$, we choose (see $\S$\ref{ss:pist})
\begin{align}\label{eq:x-ii}
&\{x_i\}_{i=1}^{\dim \b_-} = \{h_i: 1 \leq i \leq \dim \h\} \cup \{2E_{-\alpha}: \alpha \in \Delta^+\},\\
\label{eq:x-ii-1}
&\{x^i\}_{i = 1}^{\dim \b} = \{h_i: 1 \leq i \leq \dim \h\} \cup \{E_{\alpha}: \alpha \in \Delta^+\}.
\end{align}
We now give three pairs of Poisson Lie groups to be used to form mixed Poisson products.

\bex{:pair-AB}
{\rm 
With the non-degenerate pairing $\lara_{(\b_-, \b)}$ between Lie algebras $\b_-$ of $B_-$ and $\b$ of $B$ given in \eqref{eq:bb-pair},
we  have the first pair of dual Poisson Lie groups 
$((B_-, \pist), \; (B, -\pist))$.
Under the same pairing $\lara_{(\b_-, \b)}$, the Poisson Lie groups
$(B_-^{\rm op}, \pist)$ and $(B, \pist)$ are also dual to each other,
where $B_-^{\rm op}$ is the manifold $B_-$ with the group structure opposite to that on $B_-$.
A third pair of dual Poisson Lie groups is $((B_-, -\pist), (B, \pist))$, with the pairing between the Lie algebras $\b_-$ and $\b$ given by
$-\lara_{(\b_-, \b)}$. Finally we also have the 
 pair of dual (direct product) Poisson Lie groups
\[
(A, \pi_{\sA})=(B_-, \pist) \times (B_-^{\rm op}, \pist), \hs (A^*, \pi_{\sA^*}) =(B, -\pist) \times (B, \pist),
\]
where the non-degenerate pairing between the Lie algebras $\a = \b_- \oplus \b_-$ and $\a^* = \b \oplus \b$ is the direct sum
of the pairing $\lara_{(\b_-, \b)}$ with itself.
\hfill $\diamond$
}
\eex

\bex{ex:T-ext}
{\rm
Consider the maximal torus $T$ of $G$ with the zero Poisson structure as a Poisson Lie group.  By identifying
the Lie algebra $\h$ of $T$ with $\h^*$ using the bilinear form $\lara_\g|_{\h}$, $T$ is also a dual Poisson Lie group of itself.
Let $(X, \piX)$ be a $T$-Poisson manifold  with  $T$-action $\sigma: T \times X\to X$. The Poisson structure $\piX \bowtie_\sigma 0$ on $X \times T$
given in \eqref{eq:pix0} is then the mixed product Poisson structure defined using $(T, 0)$ as a Poisson Lie group.
\hfill $\diamond$
}
\eex

\subsection{The Poisson varieties $(\calC_n, \hpi_n)$ and $(\calB^n, \opi_n)$ as mixed products}\label{ss:Fn-mixed}
Consider now the complex semi-simple Poisson Lie group $(G, \pist)$, and let the notation be as in $\S$\ref{ss:pist}. 
Let $Q$ be any closed coisotropic subgroup of $(G, \pist)$, i.e., $Q$ is a closed Lie
subgroup of $G$ which is also a coisotropic with respect to the Poisson structure $\pist$. For an integer $n \geq 1$, let
\[
F_{n, \sQ} =  \overbrace{G \times_B \cdots \times_B G}^n/Q = \tF_n/Q,
\]
and denoted by $\pi_{n, \sQ}$ the Poisson structure on $F_{n, \sQ}$ that is the quotient of the Poisson structure $(\pist)^n$ on $G^n$. 
When $n = 1$, the quotient Poisson structure on $G/Q$ will be denoted as $\pi_{\sG/\sQ}$. Note that when $Q = \{e\}$ and when $Q = B$, 
we respectively have $\pi_{n, \sQ}=\tpi_n$ and $\pi_{n, \sQ}=\pi_n$. In particular,
the flag variety $\calB :=G/B$ has the Poisson structure $\pi_1 = \pi_{\sG/\sB}$, and the decorated flag variety $\calA = G/N$
has the Poisson structure $\hpi_1:=\pi_{\sG/\sN}$.

Returning to the case of an arbitrary closed coisotropic Lie subgroup $Q$ of $(G, \pist)$, set
\[
\calQ_n = \calB^{n-1} \times (G/Q),\hs n \geq 1.
\]
One then has the isomorphism
\[
\Theta_{n, \sQ}: \; F_{n, \sQ} \to \calQ_n, \; [g_1, \cdots, g_{n-1}, g_{n}]_{\stF_{n}} \longmapsto
(g_1{}_\cdot B, \ldots, g_1\cdots g_{n-1}{}_\cdot B,  g_1\cdots g_{n-1}g_{n}{}_\cdot Q).
\]
Define the Poisson structure $\pi_{{\scriptscriptstyle  \calQ_n}}$ on $\calQ_n$  by
\[
\pi_{{\scriptscriptstyle  \calQ_n}} \, =\, \Theta_{n, \sQ}(\pi_{n, \sQ}).
\]
Our goal now is to express the Poisson
manifold $(\calQ_n, \pi_{{\scriptscriptstyle  \calQ_n}})$  as a mixed product of $n-1$ copies of $(\calB, \pi_1)$ with $(G/Q, \pi_{\sG/\sQ})$.

Let $\lam_{\sG/\sQ}$ be the left Poisson action  of $(G, \pist)$ on $(G/Q, \pi_{\sG/\sQ})$ given by
translation, i.e,
\[
\lam_{\sG/\sQ}: \; (G, \pist)  \times (G/Q, \, \pi_{\sG/\sQ}) \lrw (G/Q, \, \pi_{\sG/\sQ}), \; (g, \, g_1{}_\cdot Q) \Maps gg_1{}_\cdot Q.
\]
Recall the bases $\{x_i\}_{i = 1}^{\dim \b_-}$ of $\b_-$ and $\{x^i\}_{i = 1}^{\dim \b_-}$ of $\b$ in \eqref{eq:x-ii}.
For $1 \leq j < k \leq n-1$, define $\mu_{j,k}^{(n-1)} \in \mX^2(\calB^{n-1})$ by (see notation in $\S$\ref{ss:nota-00})
\[
\mu_{j,k}^{(n-1)} =
\sum_{i=1}^{\dim \b_-} (0, \; \ldots, \; 0, \; \underset{j^{\text{th}}\,\text{entry}}{\lambda_{\sG/\sB}(x^i)}, \; 0, \; \ldots, \; 0)
\wedge (0, \; \ldots, \; 0, \; \underset{k^{\text{th}}\,\text{entry}}{\lambda_{\sG/\sB}(x_i)}, \; 0, \; \ldots, \; 0),
\]
and for $1 \leq j \leq n-1$, define $\mu_j^{(n,\sQ)}\in \mX^2(\calQ_n)$ by
\[
\mu_j^{(n,\sQ)} = \sum_{i=1}^{\dim \b_-} (0, \; \ldots, \; 0, \; \underset{j^{\text{th}}\,\text{entry}}{\lambda_{\sG/\sB}(x^i)}, \; 0, \; \ldots, \; 0)
\wedge (0, \; \ldots, \; 0, \; \lambda_{\sG/\sQ}(x_i)).
\]

\bpr{:piQ}
For any closed coisotropic subgroup $Q$ of $(G, \pist)$,  one has
\begin{equation}\label{eq:pi-Qn}
\pi_{{\scriptscriptstyle \calQ_n}} = (\pi_1,  \, \ldots, \, \pi_1, \, \pi_{\sG/\sQ}) + \sum_{1 \leq j < k \leq n-1} (\mu_{j,k}^{(n-1)}, \; 0)
+\sum_{j=1}^{n-1} \mu_j^{(n, \sQ)}.
\end{equation}
\epr

\begin{proof}
Note first that the standard $r$-matrix $r_{\rm st}$ on $\g$ given in \eqref{eq:rst} is given by
\[
r_{\rm st} = \sum_{i = 1}^{\dim \b_-} x_i \otimes x^i \in  \in \g \otimes \g.
\]
The case of $Q = \{e\}$ follows from the proof of \cite[Proposition 8.1]{Lu-Mou:mixed} by setting, in the notation of \cite[Proposition 8.1]{Lu-Mou:mixed},
$Q_+ = B$ and $r = r_{\rm st}$. 
The case for arbitrary $Q$ then follows from the following commutative diagram of Poisson manifolds
\[
\begin{tikzcd}[column sep = huge]
(F_{n, \{e\}}, \, \pi_{n, \{e\}}) \arrow[d] \arrow[r, "\Theta_{n, \{e\}}"]& (\calB^{n-1} \times G, \, \pi_{{\scriptscriptstyle \calG_n}}) \arrow[d]\\
 (F_{n, \sQ}, \, \pi_{n, \sQ})  \;\arrow[r, "\Theta_{n, \sQ}" '] & \;\; (\calQ_n, \, \pi_{{\scriptscriptstyle \calQ_n}}),
\end{tikzcd}
\]
where the second vertical arrow is ${\rm Id}_{{\scriptscriptstyle \calB^{n-1}}} \times p_{\sG/\sQ}$, with $p_{\sG/\sQ}: G \to G/Q$
the projection.
\end{proof}

For the special case of $Q = B$, let $\opi_n = \Theta_{n, \sB}(\pi_n)$. 

\bco{:Jmn} For any $n \geq 1$, the Poisson structure $\opi_n$ on $\calB^n$  is given by
\begin{equation}\label{eq:Jn-pin}
\opi_n = (\pi_1, \pi_1, \ldots, \pi_1) + \sum_{1 \leq j < k \leq n} \mu_{jk}^{(n)},
\end{equation}
where for $1 \leq j < k \leq n$, $\mu_{j,k}^{(n)} \in \mX^2(\calB^{n})$ is given by
\[
\mu_{j,k}^{(n)} =
\sum_{i=1}^{\dim \b_-} (0, \; \ldots, \; 0, \; \underset{j^{\text{th}}\,\text{entry}}{\lambda_{\sG/\sB}(x^i)}, \; 0, \; \ldots, \; 0)
\wedge (0, \; \ldots, \; 0, \; \underset{k^{\text{th}}\,\text{entry}}{\lambda_{\sG/\sB}(x_i)}, \; 0, \; \ldots, \; 0).
\]
\eco

\bre{:pin-mix}
{\rm By \eqref{eq:pi-Qn},   $\pi_{{\scriptscriptstyle \calQ_n}}$
 is a mixed product Poisson structure on $\calQ_n$ in the sense of \deref{:mixed-df}. 
The Poisson structure $\pi_{{\scriptscriptstyle \calQ_n}}$ on $\calQ_n = \calB^{n-1} \times G/Q$ can also be written as two fold mixed products:
consider again the pair of dual Poisson Lie groups
$((B_-, \pist), (B, -\pist))$ from \exref{:pair-AB}.
For $1 \leq k \leq n-1$, define the right Poisson action $\tilde{\rho}_{k}$ of $(B, -\pist)$ on $\calB^{k}$ by
\begin{align*}
\tilde{\rho}_{k}: &\; (\calB^{k}, \opi_{k}) \times (B, -\pist)\lrw (\calB^{k}, \opi_{k}),\\
&\; ((f_1, f_2, \ldots, f_{k}), \,b) \Maps (b^{-1}f_1, \, b^{-1}f_2, \ldots, b^{-1}f_{k}),
\end{align*}
and let $\lam_{n-k, \sQ}$ be the left Poisson action of $(B_-, \pist)$ on $(\calQ_{n-k}, \pi_{{\scriptscriptstyle \calQ_{n-k}}})$ given by
\begin{align*}
\lam_{n-k, \sQ}: &\; (B_-, \pist) \times (\calQ_{n-k}, \, \pi_{{\scriptscriptstyle \calQ_{n-k}}}) \lrw (\calQ_{n-k}, \, \pi_{{\scriptscriptstyle \calQ_{n-k}}}),\\
&\; (b_-, \, (f_1,  \, \ldots, \, f_{n-k-1}, \,g_\cdot Q)) \Maps  (b_-f_1,  \, \ldots, \, b_-f_{n-k-1}, \,b_-g_\cdot Q)).
\end{align*}
It then follows from \prref{:piQ} that
\[
\pi_{{\scriptscriptstyle \calQ_n}}= \opi_{k} \times_{(\tilde{\rho}_k, \lam_{n-k, \sQ})} \pi_{{\scriptscriptstyle \calQ_{n-k}}}.
\]
See \cite[Remark 6.10]{Lu-Mou:mixed}
for the general setting. In particular,
setting $\hpi_n = \pi_{n, \sQ}$ for  $Q = N$, the Poisson structure $\hpi_n$ on $\calB^{n-1} \times \calA$ then satisfies
\begin{equation}\label{eq:hpi-mix}
\hpi_n =\opi_{k} \times_{(\tilde{\rho}_k, \lam_{n-k, \sN})} \hpi_{n-k}, \hs 1 \leq k \leq n-1.
\end{equation}
\hfill $\diamond$
}
\ere

\section{The Poisson isomorphism $J_n$}\label{s:app-Jn}
\subsection{The isomorphism $J_n$}\label{ss:Jn}
For $n \geq 1$, recall the Poisson action 
\[
\lambda_n:\; (G, \pist) \times (F_n, \pi_n) \lrw (F_n, \pi_n),\;\; \lambda_n(g, [g_1, g_2, \ldots, g_n]_{\sF_n}) = [gg_1, g_2, \ldots, g_n]_{\sF_n},
\]
which in particular makes $(F_n, \pi_n)$ into a $T$-Poisson manifold. By \eqref{eq:pix0}, 
one has  the $T$-extension 
$(F_n \times T, \, \pi_n \bowtie_{\lam_n} 0)$ of $(F_n, \pi_n)$.
For notational simplicity, we set
\[
\pi_n \bowtie 0 = \pi_n \bowtie_{\lam_n} 0  = (\pi_n, 0) + \sum_{i = 1}^{\dim \h} (\lam_n(h_i), 0) \wedge (0, h_i^R) \in \mX^2(F_n \times T).
\]
Recall that the open sub-manifold $F_n^o$ of $F_n$ is defined as
\[
F_n^o = \{[g_1, g_2, \ldots, g_n]_{\sF_n}: g_1g_2 \cdots g_n \in B_-\},
\]
and that we have  the $T$-equivariant isomorphism
 $J_n: \Gamma_n \to F_n^o \times T$ (see \eqref{eq:Jn}) given by 
\[
 J_n([g_1, g_2, \ldots, g_n]_{\stF_n}) = ([g_1, g_2, \ldots, g_n]_{\sF_n}, \, [g_1g_2\cdots g_n]_0),
\]
where $T$ acts on $F_n^o \times T$ diagonally.
In this section we prove the following fact also stated as \thmref{:Jn}.

\bthm{:Jn-apdx}
For any $n \geq 1$, $J_n:  (\Gamma_n, \, \tpi_n) \to (F_n^o \times T, \, \pi_n \bowtie 0)$ is a Poisson isomorphism.
\ethm

\subsection{Some auxiliary lemmas}
Consider the pair of dual Poisson Lie groups $((B_-, \pist), (B, -\pist))$ from \exref{:pair-AB}.
 The left Poisson action $\lam_n$ of $(G, \pist)$ on $(F_n, \pi_n)$
restricts to a left Poisson action of $(B_-, \pist)$ on $(F_n, \pi_n)$, still denoted by $\lam_n$. One also has the induced
{\it right} Poisson action of $(B, -\pist)$ on $(F_n, \pi_n)$ by
\begin{align}\label{eq:rho-n}
\rho_n: \;&(F_n, \, \pi_n) \times (B, -\pist) \lrw (F_n, \pi_n), \\
\nonumber
& ([g_1, g_2, \ldots, g_n]_{\sF_n}, \, b)\longmapsto [b^{-1}g_1, g_2, \ldots, g_n]_{\sF_n}.
\end{align}
Recall that $\lam_n$ also denotes the induced Lie algebra anti-homomorphism $\lambda_n: \g \to \mX^1(F_n)$, and recall that
for $x \in \g$, $x^R$ is the right invariant vector field on $G$ with value $x$
at the identity element $e$. For every  sub-manifold $S$ of $G$ that is invariant under left translation by elements in $B_-$, we have 
the Lie algebra anti-homomorphism
\[
\lambda_{\b_-}: \b_- \lrw \mX^1(S), \;\; \lambda_{\b_-}(x) = x^R|_{{\scriptscriptstyle S}}.
\]
We further simplify notation as follows.

\bnota{:bowtie}
{\rm
If $(X, \piX)$ is a Poisson manifold with a right Poisson action
\begin{equation}\label{eq:rho-X}
\rho: \; (X, \piX) \times (B, -\pist) \lrw (X, \piX), \;\; (x, b) \Map xb,
\end{equation}
of the Poisson Lie group $(B, -\pist)$, we set
\begin{align*}
\piX \bowtie_\rho \pist& = \piX \times_{(\rho, \lambda_{\b_-})} \pist = (\piX, 0) + (0, \pist) - \sum_{i=1}^{\dim \b_-} (\rho(x^i), 0)
\wedge (0, x_i^R) \in \mX^2(X \times G),\\
\piX \bowtie_\rho \pi_n& = \piX \times_{(\rho, \lam_n)} \pist = (\piX, 0) + (0, \pi_n) - \sum_{i=1}^{\dim \b_-} (\rho(x^i), 0)
\wedge (0, \lam_n(x_i)) \in \mX^2(X \times F_n)
\end{align*}
for $n \geq 1$, where $\{x^i\}$ and $\{x_i\}$ are given in \eqref{eq:x-ii}. In the special case of $(X, \piX) = (F_m, \pi_m)$
and $\rho = \rho_m$ for $m \geq 1$ as in
\eqref{eq:rho-n}, we further simplify the notation to set
\begin{align*}
\pi_m \bowtie \pist  &= \pi_m \bowtie_{\rho_m} \pist  = (\pi_m, 0) + (0, \pist) + \sum_{i=1}^{\dim \b_-} (\lam_m(x^i), 0)
\wedge (0, x_i^R) \in  \mX^2(F_m \times G), \hs\\
\pi_m \bowtie \pi_n &= \pi_m \bowtie_{\rho_m} \pi_n =
(\pi_m, 0) + (0, \pi_n) + \sum_{i=1}^{\dim \b_-} (\lam_m(x^i), 0)
\wedge (0, \lam_n(x_i))    \in \mX^2(F_m \times F_n).
\end{align*}
\hfill $\diamond$
}
\enota

For $m, n \geq 1$,  consider the isomorphism $\Theta_{m, n}: F_{m+n} \to F_m \times F_n$ given by
\begin{equation}\label{eq:Jmn}
 \Theta_{m, n}([g_1, \ldots, g_{m+n}]_{\sF_{m+n}}) =
([g_1, \ldots, g_m]_{\sF_m}, \, [g_1 \cdots g_m g_{m+1}, \,g_{m+2}, \,\ldots, \,g_{m+n}]_{\sF_n}.
\end{equation}
Consider also the isomorphism
\begin{equation}\label{eq:tJn-1}
\widetilde{\Theta}_n: \; \tF_n \lrw F_{n-1} \times G, \; [g_1, g_2, \ldots, g_n]_{\stF_n} \Maps ([g_1, g_2, \ldots, g_{n-1}]_{\sF_n}, \; g_1g_2 \cdots g_n).
\end{equation}
The following lemma follows direct from \reref{:pin-mix}.

\ble{:Theta-mn}
For $m, n \geq 1$, one has
\[
\Theta_{m, n}(\pi_{m+n}) = \pi_m \bowtie \pi_n \hs \mbox{and} \hs
\widetilde{\Theta}_n(\tpi_n) = \pi_{n-1} \bowtie \pist.
\]
\ele

Assume that $(X, \piX)$ is a Poisson manifold with a right Poisson action $\rho$ of the Poisson Lie group $(B, -\pist)$ as in \eqref{eq:rho-X},
and consider  the isomorphism
\begin{equation}\label{eq:IX}
I_\sX :\;\; X \times B_- \lrw (X \times F_1^o) \times T, \;\; (x, \, b_-) \Maps (x, \; b_-{}_\cdot B, \; [b_-]_0).
\end{equation}
Equip $X \times B_-$ and $X \times F_1^o$ with the Poisson structures $\piX \bowtie_\rho \pist$ and $\piX \bowtie_\rho \pi_1$,
respectively.
Equip $X \times F_1^o$ with the action of $T$ by
\begin{equation}\label{eq:TXF}
\sigma_1: \; T \times (X \times F_1^o) \lrw X \times F_1^o, \;\; (t, \, (x, \, g_\cdot B))\Maps  (xt^{-1}, \; tg_\cdot B),
\end{equation}
which preserves the Poisson structure $\piX \bowtie_\rho \pi_1$, so one has the $T$-extension Poisson structure
$(\piX \bowtie_\rho \pi_1) \bowtie_{\sigma_1} 0$ on $(X \times F_1) \times T$.

\ble{:XB}
The map $I_\sX: (X \times B_-, \, \piX \bowtie_\rho \pist) \to ((X \times F_1^o) \times T, \, (\piX \bowtie_\rho \pi_1) \bowtie_{\sigma_1} 0)$
is a Poisson isomorphism.
\ele

\begin{proof}  Consider the isomorphism
\[
I_0: \; B_-\lrw F_1^o\times T, \; \;b_- \longmapsto (b_-{}_\cdot B, \, [b_-]_0), \hs b_- \in B_-.
\]
By \cite[Proposition A.7]{Lu-Mi:Kostant}, $I_0:  (B_-, \; \pist) \to (F_1^o\times T, \, \pi_1 \bowtie 0)$ is Poisson, where
\[
\pi_1 \bowtie 0 = (\pi_1, 0) + \sum_{i=1}^{\dim \h} (\lam_1(h_i), 0) \wedge (0, h_i^R).
\]
Note that $I_\sX = {\rm Id}_\sX \times I_0$.  Using the bases $\{x^i\}$ of $\b$ and  $\{x_i\}$ of $\b_-$ in \eqref{eq:x-ii},
\[
\piX \bowtie_\rho \pist = (\piX, 0) + (0, \pist) -2\sum_{\alpha \in \Delta^+} (\rho(E_\alpha), 0) \wedge (0, E_{-\alpha}^R)
-\sum_{i=1}^{\dim \h} (\rho(h_i), 0) \wedge (0, h_i^R).
\]
It follows that as bivector fields on $X \times F_1^o \times T$,
\begin{align*}
I_\sX(\piX \bowtie_\rho\pist) & = (\piX, 0, 0)  +(0, \pi_1, 0) + \sum_{i=1}^{\dim \h}(0, \, \lam_1(h_i), \, 0) \wedge (0, \, 0, \, h_i^R)\\
& \;\;\;\;-2\sum_{\alpha \in \Delta^+} (\rho(E_\alpha), 0, 0) \wedge
(0, \lam_1(E_{-\alpha}), 0) -\sum_{i=1}^{\dim \h} (\rho(h_i), 0, 0) \wedge (0, \lam_1(h_i), h_i^R).
\end{align*}
On the other hand, by \eqref{eq:TXF},
\begin{align*}
(\piX \bowtie_\rho\pi_1) \bowtie_\sigma 0& =
(\piX, 0, 0) + (0, \pi_1, 0)-2\sum_{\alpha \in \Delta^+} (\rho(E_\alpha), 0, 0) \wedge (0, \lam_1(E_{_\alpha}), 0)\\
&\; \;\;\;-\sum_{i=1}^{\dim \h} (\rho(h_i), 0, 0) \wedge (0, \lam_1(h_i), 0) +\sum_{i=1}^{\dim \h} (-\rho(h_i), \lam_1(h_i), 0) \wedge (0, 0, h_i^R).
\end{align*}
Comparing terms, one sees that $I_\sX(\piX \bowtie_\rho \pist) =(\piX \bowtie_\rho \pi_1) \bowtie_{\sigma_1} 0$.
\end{proof}

\subsection{Proof of \thmref{:Jn-apdx}}  Since 
$J_n =\Theta_{n-1, 2}\circ  I_{\sF_{n-1}} \circ \widetilde{\Theta}_n$, and
\begin{align*}
(\Gamma, \, \tpi_n) \;& \stackrel{\widetilde{\Theta}_n}{\longrightarrow} \; (F_{n-1} \times B_-, \; \pi_{n-1} \bowtie \pist) \\
& \stackrel{I_{\sF_{n-1}}}{-\!\!\!-\!\!\!-\!\!\!-\!\!\!\longrightarrow} \; ((F_{n-1} \times F_1^o) \times T, \; 
(\pi_{n-1} \bowtie \pi_1) \bowtie_{\sigma_1} 0)\\
&\stackrel{\Theta_{n-1, 1}^{-1} \times {\rm Id}_T}{-\!\!\!-\!\!\!-\!\!\!-\!\!\!\longrightarrow} \;(F_{n}^o \times T, \; \pi_{n} \bowtie 0),
\end{align*}
are Poisson maps, respectively by \leref{:Theta-mn}, \leref{:XB}, and \leref{:Theta-mn} again.
we see that  $J_n: (\Gamma_n, \tpi_n) \to (F_{n}^o \times T, \, \pi_{n} \bowtie 0)$ is Poisson. 

This finishes the proof of \thmref{:Jn-apdx}.

\section{Generalized double Bruhat cells}\label{s:app-Guv}
Generalizing the case of $m = n$ from \cite[$\S$1.4]{Lu-Mou:flags},
we introduce in this appendix the $T$-Poisson manifold $(G_{m, n}, \tpi_{m, n})$, for any pair of integers $m, n \geq 1$,
and the {\it generalized double Bruhat cells} $G^{\bfu, \bfv}$ as its $T$-leaves, where  $(\bfu, \bfv) \in W^m \times W^n$.
The main results of this appendix are presented in $\S$\ref{ss:piece}, where we establish, using any 
representative $\dbfv \in \bfv T^n$ for each $\bfv \in W^n$,  a piece-wise 
Poisson isomorphism from 
$(G_{m, n}, \tpi_{m, n})$ to $(\Gamma_{m+n}, \tpi_{m+n})$ and thus also to $(F_{m+n}^o \times T, \pi_{m+n} \bowtie 0)$. These
 piece-wise Poisson isomorphisms carry $T$-leaves to $T$-leaves, giving rise in particular to
Poisson isomorphisms (see \coref{:JEv-Poi})
\[
K_{\bfu, \dbfv}:\; (G^{\bfu, \bfv}, \tpi_{m, n}) \stackrel{\sim}{\longrightarrow}
\left(\O^{(\bfu, \bfv^{-1})}_e \times T, \,\pi_{m+n} \bowtie 0\right).
\]
As  a special case, we show in \coref{:tFn-1} that each $(\tF_n, \tpi_n)$ is also piece-wise Poisson isomorphic to 
$(F_{n+1}^o \times T, \pi_{n+1} \bowtie 0)$. Consequently every $T$-leaf of 
$(\tF_n, \tpi_n)$  is also Poisson isomorphic to $(\O^{\bfw}_e \times T, \pi_{n+1} \bowtie 0)$ for some  $\bfw \in W^{n+1}$.

\subsection{Generalized double Bruhat cells associated to conjugacy classes}\label{ss:Guv-app}
Recall from $\S$\ref{ss:Guu} that for an integer $n \geq 1$, we have the quotient  space
\[
\wF_{-n}=G \times_{B_-}  \cdots \times_{B_-} G
\]
with the well-defined quotient Poisson structure  $\tpi_{-n}$, and that the image of 
$(g_1, g_2, \ldots, g_n) \in G^n$ in $\tF_{-n}$ is denoted as  $[g_1, \ldots, g_n]_{\stF_{-n}}$. Similar to the case of $(\tF_n, \tpi_n)$, 
$(\tF_{-n}, \tpi_{-n})$ is a $T$-Poisson manifold with the $T$-action given by
\begin{equation}\label{eq:T-tn}
t \cdot [g_1, g_2, \ldots, g_n]_{\stF_{-n}} =  [tg_1, g_2, \ldots, g_n]_{\stF_{-n}}.
\end{equation}

Recall from \exref{:pair-AB} that we have the
pair of dual Poisson Lie groups
\[
(A, \pi_{\sA})=(B_-, \pist) \times (B_-^{\rm op}, \pist), \hs (A^*, \pi_{\sA^*}) =(B, -\pist) \times (B, \pist).
\]
As the Poisson structure $\pist$ on $G$ is multiplicative,  one has the right and left Poisson actions
\begin{align*}
&\tilde{\rho}_n:\;(\wF_n, \pi_{\swF_n}) \times (A^*, \, \pi_{\sA^*})  \longrightarrow (\wF_n, \pi_{\swF_n}),\\
&\tilde{\lambda}_{-n}: \;(A, \, \pi_{\sA})\times (\wF_{-n}, \pi_{\swF_{-n}}) \longrightarrow (\wF_{-n}, \pi_{\swF_{-n}})
\end{align*}
of the Poisson Lie groups $(A, \piG)$ and $(A^*, \pi_{\sA^*})$, 
respectively given by 
\begin{align}\lb{eq:trho-n}
&\tilde{\rho}_n([g_1, \,\ldots, \, g_n]_{\swF_n}, \; (b_1, b_2))= [b_1^{-1}g_1, \, g_2, \ldots, \, g_{n-1}, \, g_nb_2]_{\swF_n},\\ 
\lb{eq:tlam-n}
&\tilde{\lam}_{-n}((b_{-1}, b_{-2}), \;[g_1, \, \ldots, \, g_n]_{\swF_{-n}}) = [b_{-1} g_1, \, g_2, \ldots, \, g_{n-1}, \,
g_n b_{-2}]_{\swF_{-n}},
\end{align}
where $g_j \in G$ for each $j \in [1, n]$, and $b_1, b_2 \in B, \, b_{-1}, b_{-2} \in B_-$.

\bde{:tpimn}
For integers $m, n \geq 1$, 
the Poisson structure $\tpi_{m, n}$ on $\wF_m \times \wF_{-n}$ is defined to be the
mixed product 
\[
\tpi_{m, n} = \tpi_m \times_{(\tilde{\rho}_m, \tilde{\lam}_{-n})} \tpi_{-n}.
\]
\ede

Note that the Poisson structure $\tpi_{m,n}$ is invariant under the diagonal $T$-action on $\tF_m \times \tF_{-n}$.

\bnota{:F-C}
For a conjugacy class $C$ in $G$, let
\[
G_{m, n, \sC} = \{([g_1, \ldots, g_m]_{\swF_m},\, [k_1, \ldots, k_n]_{\swF_{-n}}): g_1\cdots g_m
(k_1 \cdots k_n)^{-1} \in C\} \subset \wF_m \times \wF_{-n},
\]
and for $\bfu = (u_1, \ldots, u_m) \in W^m$ and $\bfv = (v_1, \ldots, v_n) \in W^n$, let
\[
G^{\bfu, \bfv}_{\sC} = G_{m, n, \sC} \cap \left(
(B \bfu B) \times B_- \bfv B_-)\right) \subset \wF_m \times \wF_{-n}.
\]
We will refer to $G^{\bfu, \bfv}_{\sC}$ as the {\it generalized double Bruhat cell associated to the conjugacy
class $C$ and the sequences $\bfu$ and $\bfv$}.
One thus has the decomposition of $\wF_m \times \wF_{-n}$ into the disjoint union
\begin{equation}\lb{eq:Fmn-C}
\wF_m \times \wF_{-n} = \bigsqcup_{\bfu, \bfv, \sC} G^{\bfu, \bfv}_{\sC},
\end{equation}
where $\bfu \in W^m, \bfv \in W^n$, and $C$ runs over the set of all conjugacy classes in $G$.
\enota

\bthm{:Fmn-app}
1) For any conjugacy class $C$ in $G$ and any $\bfu \in W^m, \bfv \in W^n$, $G^{\bfu, \bfv}_{\sC}$ is
a non-empty connected smooth sub-manifold of $\wF_m \times \wF_{-n}$ of dimension $l(\bfu) + l(\bfv) + \dim C + \dim T$;

2) The decomposition \eqref{eq:Fmn-C} is that of $\wF_m \times \wF_{-n}$ into the
$T$-leaves of $\tpi_{m, n}$ for the $T$-action given in \eqref{eq:T-tFmn}.
\ethm

\begin{proof}
The case when $m = n$ is proved in \cite[$\S$1.5]{Lu-Mou:flags}. Let $m$ and $n$ be arbitrary,
and define $Z_m \subset \wF_{m+n}$ and $Z_{-n} \subset \wF_{-(m+n)}$ by
\[
Z_m=
\overbrace{{G \times_B} \cdots \times_B G}^{m \; \mbox{copies of}\; G}
\times_B \overbrace{{B \times_B \cdots
\times_B B}}^{n \; \mbox{copies of}\; B}, \hs
Z_{-n} =
\overbrace{B_- \times_{B_-} \cdots \times_{B_-} B_-}^{m \; \mbox{copies of}\; B_-}
\times_{B_-} \overbrace{{G\times_{B_-} \cdots
\times_{B_-} G}}^{n \; \mbox{copies of}\; G}.
\]
By the definition of the Poisson structure  $\tpi_{m+n, m+n}$, 
$Z_m \times Z_{-n}$ is a Poisson sub-manifold of $(\tF_{m+n} \times \tF_{-(m+n)}, \tpi_{m+n, m+n})$.
Define
$\mu: \tF_{m+n} \to \tF_m$ and $\mu_{-}: \tF_{-(m+n)} \to \tF_{-n}$ by
\begin{align*}
& \mu([g_1, \ldots,\, g_{m-1}, \,g_m, \,g_{m+1}, \,\ldots,\, g_{m+n}]_{\stF_{m+n}})= [g_1, \,\ldots,\, g_{m-1}, \,g_m g_{m+1} \cdots g_{m+n}]_{\tF_m},\\
&\mu_-([g_1, \ldots,  \,g_m, \,g_{m+1}, \,\ldots, \,g_{m+n}]_{\stF_{-(m+n)}}) = [g_1 \cdots  g_m g_{m+1}, \, g_{m+2}, \,\ldots,\, g_{m+n}]_{\tF_{-n}}.
\end{align*}
By \exref{:XGY}, one sees that
\[
\phi:= \mu|_{\sZ_m} \times \mu_{-}|_{\sZ_{-n}}: \; (Z_m \times Z_{-n}, \tpi_{m, n})
 \lrw (\wF_m \times \wF_{-n}, \,\widetilde{\pi}_{m, n})
\]
is a $T$-equivariant Poisson isomorphism.
The statement on the $T$-leaf decomposition of $(\tF_m \times \tF_{-n}, \tpi_{m, n})$ now follows from that for 
$(\tF_{m+n} \times \tF_{-(m+n)}, \tpi_{m+n, m+n})$
 given in \cite[Theorem 1.4]{Lu-Mou:flags}.
\end{proof}

Specializing to the case when $C = \{e\}$ is the trivial conjugacy class, we have the Poisson sub-manifold $G_{m, n}$ of $(\tF_m \times \tF_{-n}, \tpi_{m, n})$.,
where, we have already introduced in $\S$\ref{ss:Guu}, 
\[
G_{m, n} = G_{m, n, \{e\}} = \{([g_1, g_2, \ldots, g_m]_{\stF_m}, [h_1, h_2, \ldots, h_n]_{\stF_n}): g_1g_2 \cdots g_n = h_1h_2 \cdots h_n\}.
\]
The restriction of $\tpi_{m, n}$ to $G_{m, n}$ will still be denoted as $\pi_{m, n}$. 
Again as we have done in $\S$\ref{ss:Guu}, for $\bfu \in W^m$ and $\bfv \in W^n$, 
\[
G^{\bfu, \bfv} = G^{\bfu, \bfv}_{\{e\}} \subset G_{m, n}
\]
is called a {\it generalized double Bruhat cell}. 
We now have the following
\coref{:Tleaf-Gmn}
which  generalizes the corresponding statements  in \cite{FZ:double} on for double Bruhat cells in $G$.

\bco{:Tleaf-Gmn}
1) For any $\bfu \in W^m, \bfv \in W^n$, the generalized double Bruhat cell
$G^{\bfu, \bfv}$ is
a non-empty connected sub-manifold of $G_{m, n}$ of dimension $l(\bfu) + l(\bfv) + \dim T$;

2) the $T$-leaves of $(G_{m, n}, \tpi_{m, n})$ are
precisely the generalized double Bruhat cells $G^{\bfu, \bfv}$, where $\bfu \in W^m$ and $\bfv \in W^n$.
\eco

\subsection{Piece-wise Poisson isomorphisms from $(G_{m, n}, \tpi_{m,n})$ to $(\Gamma_{m+n}, \tpi_{m+n})$}\label{ss:piece}
Fix integers $m, n \geq 1$.
Recall  that the maximal torus $T$ acts on both $G_{m, n}$ and $\Gamma_{m+n}$, respectively via \eqref{eq:T-tFmn} and \eqref{eq:T-tFn}.
Recall also that we have set
\[
\bfv^{-1} = (v_n^{-1}, \ldots, v_2^{-1}, v_1^{-1}) \hs \mbox{if} \hs \bfv = (v_1, v_2, \ldots, v_n) \in W^n.
\]
Writing an arbitrary $\bfw \in W^{m+n}$ as $\bfw = (\bfu, \bfv^{-1})$ with
$\bfu \in W^m$ and $\bfv \in W^n$, by \prref{:T-leaves-all}, the $T$-leaf decomposition of $(\Gamma_{m+n}, \tpi_{m+n})$ can be re-written as
\[
\Gamma_{m+n} = \bigsqcup_{\bfu \in W^m, \bfv \in W^n}\Gamma^{(\bfu, \bfv^{-1})},
\]
where recall that $\Gamma^{(\bfu, \bfv^{-1})} = \Gamma_{m+n} \cap (B (\bfu, \bfv^{-1})B)$. 
For $\bfv \in W^n$, we set
\[
G_{m,n}^{\bfv} = G_{m, n} \cap (\tF_m \times (B_- \bfv B_-) \hs \mbox{and} \hs
\Gamma^{\bfv}_{m+n} = \Gamma_{m+n} \cap (\tF_m \times_B (B \bfv^{-1} B)).
\]
Here $\tF_m \times_B (B \bfv^{-1} B)$ is the quotient of $\tF_m \times (B \bfv^{-1} B)$ by the diagonal $B$-action
\[
([g_1, \ldots, g_{m-1}, g_m]_{\stF_m}, \; g) \cdot b = ([g_1, \ldots, g_{m-1}, g_mb]_{\stF_m}, \; b^{-1}g).
\]
See \exref{:quot} for the notation $X \times_B Y$ for manifolds $X$ with a right $B$-action and $Y$ with a left $B$-action.
Both $G_{m, n}^{\bfv}$ and $\Gamma_{m+n}^{\bfv}$ are then unions of $T$-leaves in the respective Poisson manifolds
$(G_{m, n}, \tpi_{m, n})$
and  $(\Gamma_{m+n}, \tpi_{m+n})$, and one has the disjoint unions
\[
G_{m, n} = \bigsqcup_{\bfv \in W^n} G_{m, n}^{\bfv} \hs \mbox{and} \hs
\Gamma_{m+n} = \bigsqcup_{\bfv \in W^n} \Gamma_{m+n}^{\bfv}.
\]

For $\bfv = (v_1, \ldots, v_n) \in W^n$ and any
representative $\dbfv = (\dv_1, \ldots, \dv_n) \in \bfv T^n$, recall that
 \[
C_{\dbfv} = C_{\dv_1}  \times \cdots \times C_{\dv_n}, \hs \mbox{where}\;\; C_{\dv_j} = N \dv_j \cap \dv_j N_- \;\; \mbox{for}\;\; 1 \leq j \leq n.
\]
Denoting an arbitrary element in $C_{\dbfv}$ as
$\dot{c} = (\dot{c}_1, \dot{c}_2, \ldots, \dot{c}_n)$ and recalling from \eqref{eq:BBv0} and \eqref{eq:BDu}  the isomorphisms,
\begin{align*}
&B_- \times C_{\dbfv} \lrw B_- \bfv B_-, \;\; (b_-, \,\dot{c}_1, \, \dot{c}_2, \, \ldots, \, \dot{c}_n) 
\Maps [b_-\dot{c}_1, \, \dot{c}_2, \, \ldots, \, \dot{c}_n]_{\stF_{-n}},\\
&B \times C_{\dbfv} \lrw B \bfv^{-1} B, \;\; (b, \,\dot{c}_1, \, \ldots, \, \dot{c}_{n-1}, \, \dot{c}_n) \Maps
[b \dot{c}_n^{-1}, \dot{c}_{n-1}^{-1}, \, \ldots,  \, \dot{c}_1^{-1}]_{\stF_n},
\end{align*}
we have the well-defined map
\begin{align}\label{eq:Ev}
E_{m,\dbfv}: & \;\;\tF_m \times (B_- \bfv B_-) \lrw \tF_m \times_B (B \bfv^{-1} B),\\
\nonumber
&\;\; ([g_1, \ldots, g_m]_{\stF_m}, [b_-\dot{c}_1, \dot{c}_2, \ldots, \dot{c}_n]_{\stF_{-n}})
\Maps [g_1, \ldots, g_m,\, \dot{c}_n^{-1}, \ldots, \,\dot{c}_2^{-1}, \,\dot{c}_1^{-1}]_{\stF_{m+n}}.
\end{align}
Note that $\tF_m \times (B_- \bfv B_-)$ is a Poisson sub-manifold of $(\tF_m \times \tF_{-n}, \,\tpi_{m, n})$, and
$\tF_m \times_B (B \bfv^{-1} B)$ is a Poisson sub-manifold of $(\tF_{m+n}, \,\tpi_{m+n})$. 

\bthm{:Ev-Poi} For any $m\geq 0$, any $\bfv  \in W^n$ and any
$\dbfv  \in \bfv T^n$,  the map
\[
E_{m,\dbfv}: \; (\tF_m \times (B_- \bfv B_-), \, \tpi_{m, n})  \lrw (\tF_m \times_B (B \bfv^{-1} B),\, \tpi_{m+n})
\]
is Poisson and restricts to a Poisson isomorphism from $(G_{m,n}^{\bfv}, \,\tpi_{m, n})$ to $(\Gamma_{m+n}^{\bfv}, \, \tpi_{m+n})$.
\ethm

\thmref{:Ev-Poi} will be proved in $\S$\ref{ss:proof-Emv}. We first give some consequences of \thmref{:Ev-Poi}.

\bre{:piece}
{\rm
The isomorphism $E_{m,\dbfv}$ depends  on the representative $\dbfv$ of $\bfv$, 
as indicated in the notation.
For a different choice $\dbfv'$, one has $E_{m, \dbfv'} = r_t \circ E_{m,\dbfv}$ for some $t \in T$, where
\[
r_t: \; \Gamma_{m+n} \lrw \Gamma_{m+n}, \; [g_1, \,\ldots, \,g_{m+n-1}, \, g_{m+n}]_{\tF_{m+n}} 
\Maps [g_1, \,\ldots, \,g_{m+n-1}, \, g_{m+n}t]_{\tF_{m+n}}.
\]

For $m = n=1$, $G_{m, n} \cong G$ and $\Gamma_2 \cong \calB \times B_-$ via $[g_1, g_2]_{\sF_2} \to (g_1{}_\cdot B, g_1g_2)$.
The piece-wise isomorphisms from $G$ to $\calB \times B_-$ have already been observed in \cite[Remark 10]{Lu-Mou:groupoids}.
See
\cite[Example 2]{Lu-Mou:groupoids} for concrete calculations for the case of $G = SL(2, \CC)$.
\hfill $\diamond$
}
\ere

For $\bfu \in W^m, \bfv \in W^n$ and  $\dbfv \in \bfv T^n$, set now
\begin{equation}\label{eq:E-uv}
E_{\bfu, \dbfv} \stackrel{{\rm def}}{=} E_{m,\dbfv}|_{\Guv}: \;\; G^{\bfu, \bfv}\lrw \Gamma^{(\bfu, \bfv^{-1})}.
\end{equation}

\bco{:Guv-emb} 
The map $E_{\bfu, \dbfv}$ gives a $T$-equivariant
Poisson isomorphism 
\[
E_{\bfu, \dbfv}: \;\; (G^{\bfu, \bfv}, \tpi_{m, n}) \lrw (\Gamma^{(\bfu, \bfv^{-1})}, \tpi_{m+n}) \subset (\Gamma_{m+n}, \tpi_{m+n}).
\]
\eco

\begin{proof} 
We have already seen that $E_{\bfu, \dbfv} : G^{\bfu, \bfv} \to \Gamma_{m+n}$ is $T$-equivariant,
and it is clear that $E_{\bfu, \dbfv} (G^{\bfu, \bfv}) \subset \Gamma^{(\bfu, \bfv^{-1})}$. We show that
$E_{\bfu, \dbfv} $ is an isomorphism by 
writing down its inverse. By \eqref{eq:CuB} and \eqref{eq:BDu},  one has the isomorphism
$B \bfu B \times C_{\dbfv} \to B(\bfu, \bfv^{-1}) B$ given by
\[
([g_1, \ldots, g_{m-1}, g_m]_{\stF_n},  \dot{c}_1,  \ldots,  \dot{c}_{n-1},\, \dot{c}_n)
\Maps [g_1, \ldots, g_{m-1}, g_m, \, \dot{c}_n^{-1},  \,\dot{c}_{n-1}^{-1}, \, \ldots,  \,\dot{c}_1^{-1} ]_{\stF_{m+n}}.
\]
We thus have
\begin{align*}
(E_{\bfu, \dbfv})^{-1}:  \Gamma^{(\bfu, \bfv^{-1})} \to G^{\bfu, \bfv}, \;\;&
[g_1, \ldots, g_{m-1}, g_m, \, \dot{c}_n^{-1}, \, \dot{c}_{n-1}^{-1}, \, \ldots, \, \dot{c}_1^{-1}]_{\stF_{m+n}} \\
&\Maps  ([g_1, \ldots, g_{m-1}, g_m]_{\tF_m}, \; [b_- \dot{c}_1, \, \ldots, \, \dot{c}_{n-1},\, \dot{c}_n]_{\tF_{-n}}),
\end{align*}
where $[g_1, \ldots, g_{m-1}, g_m]_{\tF_m} \in B \bfu B$ and $(\dot{c}_1, \, \ldots, \, \dot{c}_{n-1},\, \dot{c}_n) \in C_{\dbfv}$ are such that 
\[
 g_1 \cdots g_{m-1}g_m \,\dot{c}_n^{-1} \,\dot{c}_{n-1}^{-1}\cdots \dot{c}_1^{-1} = b_- \in B_-.
\]
\end{proof}

Recall now from $\S$\ref{s:app-Jn} that for any $n \geq 1$ we have the $T$-equivariant Poisson isomorphism 
\[
J_{n}:\; (\Gamma_n, \tpi_n) \longrightarrow (F_n^o \times T, \, \pi_n \bowtie 0), \;
J_n([g_1, g_2,\ldots, g_n]_{\stF_n}) = ([g_1, g_2,\ldots, g_n]_{\sF_n}, \, [g_1g_2\cdots g_n]_0).
\]
For $m, n \geq 1$, composing $J_{m+n}$ with the piece-wise Poisson isomorphisms from $(G_{m, n}, \tpi_{m, n})$ to 
$(\Gamma_{m+n}, \tpi_{m+n})$ in \thmref{:Ev-Poi}, we obtain piece-wise Poisson isomorphism from $(G_{m, n}, \tpi_{m, n})$ to
$(F_{m+n}^o \times T, \pi_{m+n} \bowtie 0)$, which we state in the next \coref{:JEv-Poi}.

\bco{:JEv-Poi} For any $m \geq 1$, $\bfv = (v_1, \ldots, v_n) \in W^n$, and 
$\dbfv = (\dv_1, \ldots, \dv_n) \in \bfv T^n$,  one has the $T$-equivariant Poisson isomorphism
\[
J_{m+n} \circ \left(E_{m,\dbfv}|_{G_{m, n}^{\bfv}}\right): \; (G_{m, n}^{\bfv}, \tpi_{m, n}) \longrightarrow 
\left(F_{m+n}^o \cap \left(\tF_m \times_B (B \bfv^{-1} B/B) \right)\times T,\, \pi_{m+n} \bowtie 0\right),
\]
explicitly given by 
\begin{equation}\label{eq:JE}
\left([g_1, \ldots, g_m]_{\stF_m}, \, [b_- \dot{c}_1, \,\dot{c}_2, \ldots, \,\dot{c}_n]_{\stF_{-n}}\right) \longmapsto
\left([g_1, \ldots, g_m, \,\dot{c}_n^{-1}, \ldots, \dot{c}_2^{-1}, \, \dot{c}_1^{-1}]_{\sF_{m+n}}, \, [b_-]_0\right),
\end{equation}
where $[g_1, \ldots, g_m]_{\stF_m} \in \tF_m, \, b_- \in B_-$, and $(\dot{c}_1, \, \ldots, \,\dot{c}_n) \in C_{\dbfv}$ are such that
$g_1\cdots g_m = b_-\dot{c}_1 \cdots \dot{c}_n$. 
In particular, for any $\bfu \in W^m$, one has the $T$-equivariant isomorphism 
\begin{equation}\label{eq:K-uv}
K_{\bfu, \dbfv} \stackrel{{\rm def}}{=} J_{m+n} \circ \left(E_{m,\dbfv}|_{G^{\bfu, \bfv}}\right): \; (G^{\bfu, \bfv}, \tpi_{m, n}) \longrightarrow 
\left(\O^{(\bfu, \bfv^{-1})}_e \times T, \, \pi_{m+n} \bowtie 0\right)
\end{equation}
of single $T$-leaves, explicitly given by \eqref{eq:JE} by restricting $[g_1, \ldots, g_m]_{\stF_m}$ to $B \bfu B \subset \tF_m$. 
\eco

Note now that for any $n \geq 1$, the map 
\begin{equation}\label{eq:Fnn1}
(\tF_n, \tpi_n) \longrightarrow (G_{n, 1}, \tpi_{n, 1}), \; [g_1, g_2, \ldots, g_n]_{\stF_n} \longmapsto ([g_1, g_2, \ldots, g_n]_{\stF_n}, \, g_1g_2 \cdots g_n),
\end{equation}
is a Poisson isomorphism due to the definition of the Poisson structure $\tpi_{n, 1}$ on $G_{n, 1}$. Applying \coref{:JEv-Poi} to $G_{n, 1}$, we 
also obtain piece-wise $T$-equivariant Poisson isomorphisms from $(\tF_n, \tpi_n)$ to $(F_{n+1}^o \times T, \pi_{2n+1} \bowtie 0)$ which carry
$T$-leaves  to $T$-leaves. To give the precise statement, recall from \prref{:tFn-0} that $T$-leaves of 
$(\tF_n, \tpi_n)$  are precisely of the form 
\begin{equation}\label{eq:tFn-T-leaves}
\tF_n^{\bfu, v} \stackrel{{\rm def}}{=} (B\bfu B) \cap \mu_{\stF_n}^{-1}(B_-vB_-),
\end{equation}
where $\bfu \in W^n$ and $v \in W$, and that 
\[
\mu_{\stF_n}: \; (\tF_n, \tpi_n) \lrw (G, \pist), \;\; [g_1, g_2, \ldots, g_n]_{\stF_n} \Maps g_1g_2 \cdots g_n.
\]

\bco{:tFn-1}
For any $v  \in W$ and $\dv \in vT$,  one has the $T$-equivariant Poisson isomorphism
\[
J_{n, \dv}: \; 
\mu_{\stF_n}^{-1}(B_-vB_-) \longrightarrow \left(\left(F_{n+1}^o \cap \left(\tF_n \times_B (B v^{-1} B)/B\right)\right) \times T, \, \pi_{n+1} \bowtie 0\right)
\]
explicitly given by 
\begin{equation}\label{eq:Jnv}
J_{n, \dv}([g_1, \ldots, g_n]_{\stF_n}) = ([g_1, \ldots, g_n, \, \dv^{-1}]_{\sF_{n+1}}, \, [g_1\cdots g_n \dv^{-1}]_0).
\end{equation}
In particular, for any $\bfu \in W^n$,  one has the $T$-equivariant isomorphism 
\begin{equation}\label{eq:J-uv}
J_{\bfu, \dv} \stackrel{{\rm def}}{=} J_{n, \dv}|_{\stF_n^{\bfu, v}}: \; \tF_n^{\bfu, v} \longrightarrow (\O^{(\bfu, v^{-1})}_e \times T, \, \pi_{n+1} \bowtie 0)
\end{equation}
of single $T$-leaves.
\eco

\begin{proof}
By \coref{:JEv-Poi} and using the isomorphism in \eqref{eq:Fnn1}, we have the Poisson isomorphism
\[
J_{n, \dv}: \; 
\mu_{\stF_n}^{-1}(B_-vB_-) \longrightarrow \left(\left(F_{n+1}^o \cap \left(\tF_n \times_B (B v^{-1} B)/B\right)\right) \times T, \, \pi_{n+1} \bowtie 0\right)
\]
given by 
$J_{n, \dv}([g_1, \ldots, g_n]_{\stF_n}) = ([g_1, \ldots, g_n, \, \dot{c}^{-1}]_{\sF_{n+1}}, \, [b_-]_0)$, 
where $[g_1, \ldots, g_n]_{\stF_n} \in \mu_{\stF_n}^{-1}(B_-vB_-)$, and we write $g_1g_2 \cdots g_n =b_- \dot{c}$ with unique $b_-\in B_-$
and $\dot{c} \in C_{\dv}$. Since $\dot{c}^{-1} \in \dv^{-1}N$ and $[b_-]_0 = [g_1\cdots g_n \dv^{-1}]_0$,  $J_{n, \dv}$ is also given as in 
\eqref{eq:Jnv}.
\end{proof}

\subsection{Proof of \thmref{:Ev-Poi}}\label{ss:proof-Emv}
We first prove two auxiliary lemmas.

Consider the quotient space
\[
F_n^\prime =  B \backslash \tF_n = B\backslash \overbrace{G \times_B \cdots \times_B G}^n,
\]
where $B$ acts on $\tF_n$ as a subgroup of $G$ via the action of $G$ on $\tF_n$ by \eqref{eq:G-tFn}, and denote by
$\pi_n^\prime$ the Poisson structure on $F_n^\prime$ that is the quotient Poisson structure of $\tpi_n$ on $\tF_n$.
 Define similarly the quotient
Poisson manifold $(F_{-n}^\prime,  \pi_{-n}^\prime)$ of $(\tF_{-n}, \tpi_{-n})$, where
\[
F_{-n}^\prime = B_-\backslash \tF_{-n} = B_-\backslash \overbrace{G \times_{B_-} \cdots \times_{B_-} G}^n.
\]
For every
$\bfw \in W^n$, $B\backslash B\bfw B$ is then a Poisson sub-manifold of $(F_n^\prime, \pi_n^\prime)$ and
$B_-\backslash B_- \bfw B_-$ a Poisson sub-manifold of $(F_{-n}^\prime,  \pi_{-n}^\prime)$.

Fix now $\bfv =(v_1, \ldots, v_n) \in W^n$ and any $\dbfv =(\dv_1, \ldots, \dv_n) \in \bfv T^n$. 
By \eqref{eq:BBv0} and \eqref{eq:BDu}, we have parametrizations
\begin{align*}
&C_{\dbfv} \,\stackrel{\sim}{\lrw} \,B_-\backslash B_- \bfv B_-: \;\; (\dot{c}_1, \dot{c}_2, \ldots, \dot{c}_n) \Maps [\dot{c}_1, \, \dot{c}_2, \, \ldots, \, \dot{c}_n]_{\sF_{-n}'},\\
&C_{\dbfv} \,\stackrel{\sim}{\lrw} \,B\backslash B \bfv^{-1} B: \;\; (\dot{c}_1, \dot{c}_2, \ldots, \dot{c}_n) \Maps
[\dot{c}_n^{-1}, \, \ldots, \, \dot{c}_2^{-1}, \, \dot{c}_1^{-1}]_{\sF_n^\prime}.
\end{align*}
One thus has the isomorphism
\[
\psi_{\dbfv}:\; B_-\backslash B_- \bfv B_-\lrw B \backslash B \bfv^{-1} B,\;\;
[\dot{c}_1, \, \dot{c}_2, \, \ldots, \, \dot{c}_n]_{\sF_{-n}'} \Maps [\dot{c}_n^{-1}, \, \ldots, \, \dot{c}_2^{-1}, \, \dot{c}_1^{-1}]_{\sF_n^\prime},
\]
where $(\dot{c}_1, \, \dot{c}_2, \, \ldots, \, \dot{c}_n) \in C_{\dbfv}$.

\ble{:phi-dbfv}
The map 
\begin{equation}\label{eq:psiv}
\psi_{\dbfv}: (B_-\backslash B_- \bfv B_-, \, \pi_{-n}^\prime) \to (B \backslash B \bfv^{-1} B, \, \pi_n^\prime)
\end{equation}
 is a Poisson isomorphism.
\ele

\begin{proof} Define $I_{\dbfv}:  B_-\backslash B_- \bfv B_-\to B\bfv B/B$ by
\[
I_{\dbfv}([\dot{c}_1, \, \dot{c}_2, \, \ldots, \, \dot{c}_n]_{\sF_{-n}^\prime}) = [\dot{c}_1, \, \dot{c}_2, \, \ldots, \, \dot{c}_n]_{\sF_{n}},
\hs (\dot{c}_1, \, \dot{c}_2, \, \ldots, \, \dot{c}_n) \in C_{\dbfv}.
\]
By \cite[Lemma 7.1]{Mou:Guu}, $I_{\dbfv}:  (B_-\backslash B_- \bfv B_-, \, \pi_{-n}^\prime) \to  (B\bfv B/B, \, -\pi_n)$ is Poisson.
One also has the Poisson isomorphism
\begin{equation}\label{eq:In}
I_n: \; (F_n, \pi_n) \lrw (F_n^\prime, -\pi_n^\prime), \;  [g_1, g_2, \ldots, g_n]_{\sF_n} \Maps [g_n^{-1}, \ldots, g_2^{-1}, g_1^{-1}]_{\sF_n^\prime}.
\end{equation}
Thus $\phi_{\dbfv} = I_n \circ I_{\dbfv}$ is Poisson as stated.
\end{proof}

Recall from \exref{:pair-AB} the pair $((B_-^{\rm op}, \pist), (B, \pist))$ of dual Poisson Lie groups. Note that one has the left Poisson action
\begin{align}\label{eq:lam1}
\lam_{\bfv}^-:&\;\;(B_-^{\rm op}, \pist) \times (B_-\backslash B_-\bfv B_-, \, \pi_{-n}^\prime) \lrw (B_-\backslash B_-\bfv B_-, \, \pi_{-n}^\prime),\\
&\;\;(b_-, \; [g_1, \,\ldots, \,g_{n-1}, \,g_n]_{\sF_{-n}^\prime}) \Maps[g_1, \,\ldots, \,g_{n-1}, \,g_nb_-]_{\sF_{-n}^\prime},
\end{align}
of the Poisson Lie group $(B_-^{\rm op}, \pist)$. Let
\begin{equation}\label{eq:lam-v}
\lambda_{\dbfv^{-1}}: \; (B_-^{\rm op}, \pist) \times (B\backslash B\bfv^{-1} B, \, \pi_{n}^\prime) \lrw (B\backslash B\bfv^{-1} B, \, \pi_{n}^\prime)
\end{equation}
be the unique left Poisson action of $(B_-^{\rm op}, \pist)$ on $(B\backslash B\bfv^{-1} B, \, \pi_{n}^\prime)$ such that
$\psi_{\dbfv}$ in \eqref{eq:psiv}
becomes an isomorphism of Poisson manifolds with left Poisson actions by the Poisson Lie group $(B_-^{\rm op}, \pist)$.
Let $\rho_+$ be the right Poisson action of $(B, \pist)$ on itself by right translation. Using the pair $(\rho_+, \lam_{\dbfv^{-1}})$ of  respectively right and left Poisson 
actions of
the Poisson Lie groups $(B, \pist)$ and $(B_-^{\rm op}, \pist)$, one then has the mixed product Poisson structure
$\pist \times_{(\rho_+, \lam_{\dbfv^{-1}})}\pi_n^\prime$ on $B \times (B\backslash B \bfv^{-1} B)$. Define the isomorphism
\begin{align*}
K_{\dbfv}: &\;\; B \bfv^{-1} B \lrw B \times (B\backslash B \bfv^{-1} B), \;\\
&\;\;[b\dot{c}_n^{-1}, \dot{c}_{n-1}^{-1}, \ldots, \dot{c}_1^{-1}]_{\stF_n} \Maps (b, \, [\dot{c}_n^{-1}, \dot{c}_{n-1}^{-1}, \ldots, \dot{c}_1^{-1}]_{\sF_n^\prime}),
\end{align*}
where again $(\dot{c}_1, \ldots, \dot{c}_{n-1}, \dot{c}_n) \in C_{\dbfv}$.

\ble{:BvB-1} The map $K_{\dbfv}: (B \bfv^{-1} B, \, \tpi_n) \to (B \times (B\backslash B \bfv^{-1} B), \, \pist \times_{(\rho_+, \lam_{\dbfv^{-1}})}\pi_n^\prime)$
is a Poisson isomorphism.
\ele

\begin{proof} Note first that $K_{\dbfv} = K_n \circ J_{\dbfv} \circ \tilde{I}_n$ with
\[
B \bfv^{-1} B \, \stackrel{\tilde{I}_n}{\lrw}\, B\bfv B \, \stackrel{J_{\dbfv}}{\lrw}\, (B\bfv B/B) \times B \, \stackrel{K_n}{\lrw}\,
B \times (B\backslash B \bfv^{-1} B),
\]
where $\tilde{I}_n: B \bfv^{-1} B \to B \bfv B,
 [g_1, g_2, \ldots, g_n]_{\stF_n} \to [g_n^{-1}, \ldots, g_2^{-1}, g_1^{-1}]_{\stF_n}$,
\[
J_{\dbfv}: \; B\bfv B \lrw (B\bfv B/B) \times B, \; [\dot{c}_1, \ldots, \dot{c}_{n-1}, \dot{c}_nb]_{\stF_n} \Maps ([\dot{c}_1, \ldots, \dot{c}_{n-1}, \dot{c}_n]_{\sF_n}, \, b),
\]
where $(\dot{c}_1, \ldots, \dot{c}_{n-1}, \dot{c}_n) \in C_{\dbfv}$, and
$K_n: (B\bfv B/B) \times B \to B \times (B\backslash B \bfv^{-1} B)$ is given by
\[
K_n([\dot{c}_1, \ldots, \dot{c}_{n-1}, \dot{c}_n]_{\sF_n}, \, b) \Maps (b^{-1}, \; [\dot{c}_n^{-1}, \, \dot{c}_{n-1}^{-1},\, \ldots, \, \dot{c}_1^{-1}]_{\sF_n^\prime}),
\]
where again $(\dot{c}_1, \ldots, \dot{c}_{n-1}, \dot{c}_n) \in C_{\dbfv}$. It is clear that $\tilde{I}_n: (B \bfv^{-1} B, \, \tpi_n) \to (B \bfv B, \,-\tpi_n)$ is Poisson.
Denote $\pi = \pist \times_{(\rho_+, \lam_{\dbfv^{-1}})}\pi_n^\prime$. To show $K_{\dbfv}(\tpi_n) = \pi$, one needs to show
\[
J_{\dbfv}(\tpi_n) = -K_n^{-1}(\pi)
\]
as Poisson structures on $(B\bfv B/B) \times B$. The Poisson structure $J_{\dbfv}(\tpi_n)$ has been shown in
\cite[Proposition 7.3]{Mou:Guu} to be a mixed product. We now compute $-K_n^{-1}(\pi)$ using the definition of $\pi$ and then compare
$-K_n^{-1}(\pi)$ with the formula for $J_{\dbfv}(\tpi_n)$ given in \cite[Proposition 7.3]{Mou:Guu}.

Let $\{x_i\}_{i = 1}^{\dim \b_-}$ be any basis of $\b_-$ with $\{x^i\}_{i=1}^{\dim \b_-}$ the dual basis of $\b$ under the pairing between
$\langle \,, \, \rangle_{(\b_-, \b)}$ in \eqref{eq:bb-pair}. Recall that for $x \in \b$, $x^L$ (resp. $x^R$)
denotes the left (resp. right) invariant vector field of $B$ with value $x$ at the
identity element of $B$. By the definition of $\pi$, we have
\begin{align}\nonumber
-K_n^{-1}(\pi) & = K_n^{-1}\left((-\pist, \, -\pi_n^\prime) + \sum_{i=1}^{\dim \b_-}((x^i)^L, \,0) \wedge (0, \,\lam_{\dbfv^{-1}}(x_i))\right)\\
\label{eq:Jpi}
& = (\pi_n, \, \pist) + \sum_{i=1}^{\dim \b_-} (I_n^{-1}(\lam_{\dbfv^{-1}}(x_i)), \, 0) \wedge (0, \, (x^i)^R),
\end{align}
where $I_n: F_n \to F_n^\prime$ is given in \eqref{eq:In}.

Consider now the Poisson Lie groups $(B_-, -\pist)$ and $(B, \pist)$ which become dual Poisson
Lie groups under the pairing $-\langle \,, \, \rangle_{(\b_-, \b)}$. Let now the right Poisson action
\begin{align*}
&\;\;(B_-\backslash B_-\bfv B_-, \, -\pi_{-n}^\prime) \times (B_-, -\pist)\lrw (B_-\backslash B_-\bfv B_-, \, -\pi_{-n}^\prime),\\
&\;\;([g_1, \,\ldots, \,g_{n-1}, \,g_n]_{\sF_{-n}^\prime}, \, b_-) \Maps[g_1, \,\ldots, \,g_{n-1}, \,g_nb_-]_{\sF_{-n}^\prime},
\end{align*}
of the Poisson Lie group $(B_-, -\pist)$, and let
\[
\rho_{\dbfv}: \; (B\bfv B/B, \, \pi_n) \times (B_-, -\pist) \lrw (B\bfv B/B, \, \pi_n)
\]
be the unique right Poisson action of $(B_-, -\pist)$ on $(B\bfv B/B, \, \pi_n) $ such that
\[
I_{\dbfv}: \;\;(B_-\backslash B_- \bfv B_-, \, -\pi_{-n}^\prime) \lrw (B \bfv^{-1} B/B, \, \pi_n)
\]
becomes an isomorphism of Poisson manifolds with right Poisson actions by the Poisson Lie group $(B_-, -\pist)$.
Let $\lam_+$ be the left Poisson action of $(B, \pist)$ on itself by left translation. Using the pair $(\rho_{\dbfv}, \lam_+)$ of 
respectively right and left Poisson actions 
of
the Poisson Lie groups $(B_-, -\pist)$ and $(B, \pist)$, one then has the mixed product Poisson structure
$\pi_n \times_{(\rho_{\dbfv}, \lam_+)} \pist$ on   $(B \bfv B/B) \times B$.  By \cite[Proposition 7.3]{Mou:Guu},
\[
J_{\dbfv}(\tpi_n) =\pi_n \times_{(\rho_{\dbfv}, \lam_+)} \pist =  (\pi_n, \, \pist) + \sum_{i=1}^{\dim \b_-} (\rho_{\dbfv}(x_i), \, 0) \wedge (0, \, (x^i)^R).
\]
One now checks from the definitions of the actions $\rho_{\dbfv}$ and $\lam_{\dbfv^{-1}}$ that
\[
\rho_{\dbfv}(x) = I_n^{-1}(\lam_{\dbfv^{-1}}(x)) \in \mX^1(B\bfv B/B), \hs \forall \; x \in \b_-.
\]
Comparing with \eqref{eq:Jpi}, one shows \eqref{eq:Jpi}. This finishes the proof of \leref{:BvB-1}.
\end{proof}

\bre{:BvB-decom}
{\rm
The proof of \cite[Proposition 7.3]{Mou:Guu} uses a quotient space of the Drinfeld double of $(G, \pist)$, and one can 
also use similar arguments to prove \leref{:BvB-1} directly. 
\hfill $\diamond$
}
\ere

\noindent
{\it Proof of \thmref{:Ev-Poi}}. 
Let $p_{-n}^\prime: \tF_{-n} \to F_{-n}^\prime$ be the projection map, and let
\[
P_n = {\rm Id}_{\stF_m} \times p_{-n}^\prime: \;\; \tF_m \times B_- \bfv B_- \lrw \tF_m \times (B_-\backslash B_- \bfv B_-).
\]
By the definition of the Poisson structure $\tpi_{m,n}$,
\[
P_n: \; (\tF_m \times B_- \bfv B_- , \; \tpi_{m,n}) \lrw (\tF_m \times (B_-\backslash B_- \bfv B_-), \; \tpi_m \times_{(\tilde{\rho}_m, \lam_{\bfv}^-)} \pi_{-n}^\prime)
\]
is Poisson, where $\lam_{\bfv}^-$ is given in \eqref{eq:lam1}, and $\tilde{\rho}_m$ is the right Poisson action
\[
\tilde{\rho}_m: \; (\tF_m, \, \tpi_m) \times (B, \pist) \lrw (\tF_m, \, \tpi_m), \; ([g_1, g_2, \ldots, g_m]_{\stF_m}, b) \Maps
[g_1, g_2, \ldots, g_m b]_{\stF_m}.
\]
Let $P_n^\prime = {\rm Id}_{\stF_m} \times \phi_{\dbfv}:
\tF_m \times (B_-\backslash B_- \bfv B_-)\to
\tF_m \times (B \bfv^{-1} B/B)$.
By \leref{:phi-dbfv} and by the definition of the Poisson action
$\lambda_{\dbfv^{-1}}$ in \eqref{eq:lam-v},
\[
P_n^\prime: \;
(\tF_m \times (B_-\backslash B_- \bfv B_-), \; \tpi_m \times_{(\tilde{\rho}_m, \lam_{\bfv}^-)} \pi_{-n}^\prime) \lrw
(\tF_m \times (B \bfv^{-1} B/B), \, \tpi_m \times_{(\tilde{\rho}_m, \lambda_{\dbfv^{-1}})} \pi_n)
\]
is Poisson. On the other hand, by \leref{:BvB-1} and \exref{:XGY},  the map
\begin{align*}
Q_n: &\; (\tF_m \times (B \bfv^{-1} B/B), \, \tpi_m \times_{(\tilde{\rho}_m, \lambda_{\dbfv^{-1}})} \pi_n) \lrw
(\tF_m \times_B B \bfv^{-1} B, \, \tpi_{m+n}),\\
& \;\;([g_1, g_2, \ldots, g_m]_{\stF_m}, \; [\dot{c}_n^{-1}, \,\ldots, \,\dot{c}_2^{-1}, \, \dot{c}_1^{-1}]_{\sF_n}) \Maps
 [g_1, \,g_2, \,\ldots, \,g_m,\, \dot{c}_n^{-1}, \,\ldots, \,\dot{c}_2^{-1}, \, \dot{c}_1^{-1}]_{\stF_{m+n}},
\end{align*}
is a Poisson isomorphism, where $g_1, \ldots, g_m \in G$ and $(\dot{c}_1, \dot{c}_2, \ldots, \dot{c}_n) \in C_{\dbfv}$.  As
\[
E_{m, \dbfv} = Q_n \circ P_n^\prime \circ P_n,
\]
one concludes that $E_{m, \dbfv}$ is Poisson as stated. This finishes the proof of \thmref{:Ev-Poi}.

\section{Symplectic leaves of $(\O^\bfw_e \times T, \pi_n \bowtie 0)$}\label{s:app-leaves}
In this appendix, we assume that $G$ is connected and simply connected, and we describe in \thmref{:leaves-Ow} 
the symplectic leaves of
$(\O^\bfw_e \times T, \pi_n \bowtie 0)$ for arbitrary $\bfw \in W^n$. We then apply \thmref{:leaves-Ow} 
to obtain explicit descriptions of all the symplectic
leaves in the three series
\[
(F_n^o \times T, \,\pi_n \bowtie 0), \hs (G_{m+n}, \,\tpi_{m, n}), \hs (\tF_n, \,\tpi_n), \hs m, n \geq 1.
\]
In particular, we describe all symplectic leaves in all generalized double Bruhat cells $G^{\bfu, \bfv}$, 
generalizing the result of \cite{k-z:leaves} for the case
of $u, v \in W$.

\subsection{Notation}\label{ss:nota-D}
Assume that $G$ is connected and simply connected. Recall from $\S$\ref{ss:pist} that $\h$ denotes the Lie algebra of 
the maximal torus $T  = B \cap B_-$ of $G$.
Let $X_*(T)$ and $X^*(T)$ be respectively the co-character and the character lattices of $T$. 
Then
\[
X_*(T) \otimes_{\ZZ} \CC = \h \hs \mbox{and} \hs X^*(T) \otimes_{\ZZ} \CC = \h^*,
\]
and we also regard $X_*(T)$ and $X^*(T)$ as respective subsets of $\h$ and $\h^*$. 
For $\lambda \in X^*(T)$ and $t \in T$, write $t^\lambda \in \CC^\times$ for the value of $\lambda$ on $t$. 
Recall that for $w \in W$, we 
denote by $wT$ the set of all representative of $w$ in the normalizer subgroup $N_\sG(T)$ of $T$.
One has the right action
\[
T \times W \lrw T, \;\; (t, w) \longmapsto t^w \stackrel{{\rm def}}{=} \dot{w}^{-1} t \dot{w},
\]
where $\dot{w} \in wT$ for $w \in W$, and  $(t^w)^\lambda = t^{w\lambda}$ for $t \in T, w \in W$, and $\lambda \in X^*(T)$.

Let $\Phi_0\subset X^*(T)$ and $\{\alpha^\vee: \alpha \in \Phi_0\} \subset X_*(T)$ be respectively the set of simple roots and the set of simple
co-roots determined by $B$.
For $\al \in \Phi_0$, fix root vectors $e_\al$ for $\al$ and  
$e_{-\al}$  for $-\alpha$  such that $[e_\al, e_{-\al}] =\alpha^\vee\in \h$.
Let
$x_{\pm \al}: \CC \to G$ be the one-parameter subgroups given by
\[
x_\al(z) = \exp(z \,e_\al), \hs x_{-\al}(z) = \exp(z\,e_{-\al}), \hs z \in \CC.
\]
For $\alpha \in \Phi_0$, let $s_\alpha \in W$ be the corresponding simple reflection,
and choose $\overline{s}_\alpha \in s_\alpha T$ 
by $\overline{s}_\alpha = x_\alpha(-1) x_{-\alpha}(1) x_\alpha(-1)$. 
For future use, we also note that 
\begin{equation}\label{eq:SL2}
x_{-\alpha}(z) = x_\alpha(z^{-1}) \overline{s}_\alpha \alpha^\vee(z) x_\alpha(z^{-1}), \hs \alpha \in \Phi_0, \; z \in \CC^\times.
\end{equation}
By \cite[$\S$1.4]{FZ:double}, for any reduced decomposition $w = s_{\al_1}s_{\al_2} \cdots s_{\al_l}$ of $w$, the element
\begin{equation}\label{eq:wbar}
\overline{w} \stackrel{\rm def}{=} \overline{s}_{\alpha_1} \;\overline{s}_{\alpha_2} \;\cdots \; \overline{s}_{\alpha_l} \in N_\sG(T)
\end{equation}
represents $w$ and is independent of the choice of the reduced decomposition. 

Let $\{\omega_\alpha: \alpha \in \Phi_0\}$ be the set of fundamental weights. For $\alpha \in \Phi_0$, let $\Delta^{\omega_\al} \in \CC[G]$ 
be the corresponding principal minor on $G$, uniquely determined by
\[
\Delta^{\omega_\alpha}(g_-g_0g_+) = (g_0)^{\omega_\alpha}, \hs g_- \in N_-,\, g_0 \in T, \,g_+ \in N.
\]
By \cite[Proposition 2.3]{FZ:double}, when $\alpha, \alpha' \in \Phi_0$ and $\alpha \neq \alpha'$, one has
\begin{equation}\label{eq:gig}
\Delta^{\omega_\al}(g x_{\alpha'}(z)\overline{s}_{\alpha'}) = \Delta^{\omega_\alpha}(g), \hs \forall\;\; g \in G, \, z \in \CC.
\end{equation}
More generally, let $u \in W$, let ${\rm supp}(u)$ be the set of all $\alpha \in \Phi_0$ such that $s_\alpha$ appears in some, equivalently every,
reduced decomposition of $u$, and let 
\[
{\rm supp}^o(u)  =\{\alpha \in \Phi_0: u\omega_\alpha = \omega_\alpha\}.
\]
Recall that $C_{\overline{u}} = N \overline{u} \cap \overline{u} N_-$. 
If $u = s_{\al_1}s_{\al_2} \cdots s_{\al_l}$ is a reduced decomposition  of $u$, setting $\bfu = (s_{\alpha_1}, \ldots, s_{\alpha_l})$, 
one then  has the isomorphism \cite[Proposition 2.11]{FZ:double}
\begin{equation}\label{eq:g-bfu}
g_{\bfu}: \;\; \CC^l \lrw C_{\overline{u}}, \; g_\bfu(z_1, z_2, \ldots, z_l) = x_{\alpha_1}(z_1) \overline{s}_{\alpha_1}
\, x_{\alpha_2}(z_2) \overline{s}_{\alpha_2}\, \cdots \, x_{\alpha_l}(z_l) \overline{s}_{\alpha_l}.
\end{equation}
It then follows from \eqref{eq:gig} that 
\begin{equation}\label{eq:c-gig}
\Delta^{\omega_\al}(g c) = \Delta^{\omega_\alpha}(cg) =\Delta^{\omega_\alpha}(g), \hs \forall\;\; \alpha \in {\rm supp}^o(u), \;  g \in G, \,  \; c \in C_{\overline{u}}.
\end{equation}
We also note the open embedding \cite[Proposition 2.7]{Lusz:pos} (see also \cite[Proposition 2.18]{FZ:double}) 
\begin{equation}\label{eq:m-bfu}
m_\bfu: \;\ (\CC^\times)^l \longrightarrow N_- \cap B u B, \;\;
m_\bfu(\varepsilon_1, \varepsilon_2, \ldots, \varepsilon_l)= 
x_{-\alpha_1} (\varepsilon_1) x_{-\alpha_2} (\varepsilon_2)  \cdots x_{-\alpha_l} (\varepsilon_l).
\end{equation}

\subsection{Description of symplectic leaves of $(\O^\bfw_e \times T, \pi_n \bowtie 0)$}\label{ss:leaves-D}
Let $\bfw = (w_1, \ldots, w_n) \in W^n$, and let $\dbfw = (\dw_1, \ldots, \dw_n) \in \bfw T^n$ be arbitrary. 
Recall from \eqref{eq:CD-bfu} that
$C_{\dbfw} = C_{\dw_1} \times \cdots \times C_{\dw_n}$.
By \eqref{eq:C-Ou}, one has the isomorphism
\begin{equation}\label{eq:sigma-w}
\rho_{{\dbfw}}:\;\; C_{{\dbfw}} \lrw \O^\bfw, \;\; (\dot{c}_1, \dot{c}_2, \ldots, \dot{c}_n) \longmapsto [\dot{c}_1, \dot{c}_2, \ldots, \dot{c}_n]_{\sF_n}.
\end{equation}
Recall that $\O^\bfw_e = \{[g_1, \ldots, g_n]_{\sF_n} \in \O^\bfw: g_1g_2 \cdots g_n \in B_-B\}$. For 
$(\dot{c}_1, \dot{c}_2, \ldots, \dot{c}_n) \in C_{{\dbfw}}$, then
\[
 [\dot{c}_1, \dot{c}_2, \ldots, \dot{c}_n]_{\sF_n} \in \O^\bfw_e \hs \mbox{iff} \hs \dot{c}_1\dot{c}_2 \cdots \dot{c}_n \in B_-B.
\]
Recall also that for $x \in B_-B$, we write $x = [x]_-[x]_0[x]_+$, where $[x]_- \in N_-, [x]_0 \in T$, and $[x]_+ \in N$.
One thus has the well-defined map
\begin{equation}\label{eq:tauw}
\tau_{\dbfw}: \; \O^{\bfw}_e \lrw T, \; [\dot{c}_1, \dot{c}_2, \ldots, \dot{c}_n]_{\sF_n} \longmapsto [\dot{c}_1\dot{c}_2 \cdots \dot{c}_n]_0, 
\end{equation}
where $(\dot{c}_1, \dot{c}_2, \ldots, \dot{c}_n) \in C_{{\dbfw}}$ and  $\dot{c}_1\dot{c}_2 \cdots \dot{c}_n \in B_-B$.  
Recall the $T$-action on $\O^\bfw \subset F_n$ in \eqref{eq:T-Fn}.

\ble{:ccc-T}  For any $\bfw = (w_1, \ldots, w_n) \in W^n$ and $\dbfw = (\dw_1, \ldots, \dw_n) \in \bfw T^n$, 
 one has 
\begin{equation}\label{eq:tau-aq0}
\tau_{\dbfw}(a \cdot q) = h (h^{-1})^w \tau_{\dbfw} (q), \hs h \in T, \, q \in \O^\bfw_e,
\end{equation}
where $w = w_1w_2 \cdots w_n \in W$.
\ele

\begin{proof} 
Let $(\dot{c}_1, \ldots, \dot{c}_n) \in C_{\overline{\bfw}}$ and write $\dot{c}_i = x_i \dot{w}_i$, where 
$x_i \in N \cap \dot{w}_i N_- \dot{w}_i^{\, -1}$ for $i \in [1, n]$. For $h \in T$, one then has
\[
h \cdot [\dot{c}_1, \dot{c}_2, \ldots, \dot{c}_n]_{\sF_n} = [h\dot{c}_1, \dot{c}_2, \ldots, \dot{c}_n]_{\sF_n} =
 [x_1^\prime \dot{w}_1, x_2^\prime \dot{w}_2, \ldots, x_n^\prime \dot{w}_n]_{\sF_n},
\]
where $x_1^\prime = hx_1h^{-1}$ and  $x_i^\prime = h^{w_1\cdots w_{i-1}} x_i   (h^{-1})^{w_1\cdots w_{i-1}}$ for $i\in [2,n]$.
Now  \eqref{eq:tau-aq0} follows from 
\[
(x_1^\prime \dot{w}_1) (x_2^\prime \dot{w}_2) \cdots (x_n^\prime \dot{w}_n) = h \dot{c}_1\dot{c}_2 \cdots \dot{c}_n (h^{-1})^w.
\]
\end{proof}

\bre{:dots}
{\rm
The definition of the map $\tau_{\dbfw}: \O^\bfw_e \to T$ depends on the choice of $\dbfw \in \bfw T^n$ as indicated in the notation.
If $\hat{\bfw} = (\hat{w}_1, \ldots, \hat{w}_n) \in \bfw T^n$ is another choice, and if  
\[
\dot{w}_1 \dot{w}_2 \cdots \dot{w}_n = \hat{w}_1 \hat{w}_2 \cdots \hat{w}_n t \in w T^n,
\]
where $t \in T$,
then each $(\dot{c}_1, \ldots, \dot{c}_n) \in C_{\dbfw}$ corresponds to a unique $(\hat{c}_1, \ldots, \hat{c}_n) \in C_{\hat{\bfw}}$
such that
$[\dot{c}_1, \,\ldots, \,\dot{c}_{n-1}, \,\dot{c}_n]_{\stF_n} = [\hat{c}_1, \,\ldots, \,\hat{c}_{n-1}, \,\hat{c}_n t']_{\stF_n}$.
It follows that 
\begin{equation}\label{eq:tauw-tauw}
\tau_{\dbfw}(q) = t \tau_{\hat{\bfw}}(q), \hs q \in \O^\bfw_e.
\end{equation}
\hfill $\diamond$
}
\ere

Let again $\bfw = (w_1, w_2, \ldots, w_n) \in W^n$ and $w = w_1w_2 \cdots w_n \in W$. Set
\begin{equation}\label{eq:Tw-00}
T^w = \{a (a^{-1})^w: a \in T\}.
\end{equation}
Then $\dim T^w = {\rm dim} ({\rm Im}(1-w))$, 
where $1-w: \h \to \h$. For
$\dbfw \in \bfw T^n$, let 
\begin{align}\label{eq:muw-1dot}
&\mu_{\dbfw}: \;\; \O^\bfw_e \times T \longrightarrow T/T^w, \;\; (q, t) \longmapsto t^{-2} \tau_{\dbfw}(q)_\cdot T^w.
\end{align}
Recall that  $l(\bfw) = l(w_1) + l(w_2) + \cdots + l(w_n)$.

\bpr{:leaves-conn}
For any $\bfw \in W^n$ and $\dbfw \in \bfw T^n$, symplectic leaves of $(\O^\bfw_e \times T, \pi_n \bowtie 0)$ are precisely all 
the connected components of the level sets of the map $\mu_{\dbfw}$. 
In particular, all symplectic leaves of $(\O^\bfw_e \times T, \pi_n \bowtie 0)$ have dimension 
$l(\bfw)+{\rm dim} {(\rm Im}(1-w))$. 
\epr

\begin{proof} As 
$\mu_{{\dbfw}}$ is a surjective submersion, all of its level sets are smooth and have dimension equal to 
$l(\bfw)+{\rm dim} {(\rm Im}(1-w))$.
By \reref{:dots}, the collection of level sets of the map $\mu_{\dbfw}$ is
independent of the choice of the representative $\dbfw \in \bfw T^n$.
We may thus choose $\overline{\bfw} = (\overline{w}_1, \ldots, \overline{w}_n) \in \bfw T^n$.

For notational simplicity, we set $\pi = \pi_n \bowtie 0$, and let 
\[
\mu: \; \O^\bfw_e \times T \longrightarrow  T, \; \mu(q, t) = t^{-2} \tau_{\overline{\bfw}}(q).
\]
For $(q, t) \in \O^\bfw_e \times T$, let 
\[
\pi_{(q, t)}^\#: \; \; T^*_{(q, t)}(\O^\bfw_e \times T) \longrightarrow T_{(q, t)}(\O^\bfw_e \times T), \;  
\left(\pi_{(q, t)}^\#(\beta_1), \beta_2\right)
= \pi(q, t)(\beta_1, \beta_2), 
\]
where $\beta_1, \beta_2 \in T^*_{(q, t)}(\O^\bfw_e \times T)$,
and we use translation in $T$ to identify the tangent space of $T$ at $\mu(q, t) \in T$ with $\h = {\rm Lie}(T)$. 
Note that the Lie algebra of $T^w$ is ${\rm Im}(1-w) \subset \h$. As $\mu$ is a surjective submersion, 
it is enough to show that for every $(q, t) \in \O^\bfw_e \times T$, 
\begin{equation}\label{eq:need}
{\rm Im}\left(\pi_{(q, t)}^\#\right) = \mu_*^{-1}({\rm Im}(1-w)),
\end{equation}
where $\mu_*: T_{(q, t)}(\O^\bfw_e \times T) \to T_{\mu(q, t)} T \cong \h$ is the differential of $\mu$ at $(q, t)$.
Recall that $\O^\bfw_e \times T$ is a single $T$-leaf of $\pi$, and that 
by \cite[Proposition 2.24]{Lu-Mi:Kostant}, the co-rank of $\pi$ in $\O^\bfw_e \times T$ 
is equal to the co-dimension of $T^w$ in $T$. Thus the two vector spaces on the two sides of \eqref{eq:need} have the same dimension, and it is
enough to show that
${\rm Im}\left(\pi_{(q, t)}^\#\right) \subset \mu_*^{-1}({\rm Im}(1-w))$. 

For $\chi \in X^*(T)$, denote by $\mu^\chi$ the regular function on $\O^\bfw_e \times T$ defined by 
\[
\mu^\chi(q, t) = (\mu(q, t))^\chi = t^{-2 \chi} \tau_{\overline{\bfw}}(q)^{\chi}.
\]
 For any regular function $f$ on $\O^\bfw_e \times T$ that is a weight
vector with weight $\chi_f \in X^*(T)$ for the diagonal action of $T$ on $\O^\bfw_e \times T$, one has 
by \cite[Corollary 2.15]{Lu-Mi:Kostant} and the definition $\pi$ that
\begin{equation}\label{eq:mu-f}
\{\mu^\chi, \, f\}_\pi = \langle \chi - w\chi, \, \chi_f\rangle \mu^\chi f,
\end{equation}
where $\{\, , \, \}_\pi$ is the Poisson bracket on the coordinate ring of $\O^\bfw_e \times T$ defined by $\pi$. For $\chi \in \h^*$, let $\tilde{\chi}$ 
be the left invariant $1$-form on $T$ with value $\chi$ at the identity element.  By \eqref{eq:mu-f},
\[
\pi^\#(\mu^*(\tilde{\chi})) = \sigma(\chi^\# - w \chi^\#) \in {\mathfrak{X}}^1(\O^\bfw_e \times T),
\]
where $\sigma: \h \to {\mathfrak{X}}^1(\O^w_e \times T)$ is the Lie algebra homomorphism defined by the diagonal $T$-action on
$\O^\bfw_e \times T$, and $\chi^\# \in \h$ is such that $\chi^\prime(\chi^\#) = \langle \chi, \chi'\rangle$ for $\chi' \in \h^*$. 
Note that for $\chi \in \h^*$, $\chi|_{{\rm Im}(1-w)} = 0$ if and only if $\chi = w\chi$. 
Now for $\chi \in \h^*$ with $\chi = w\chi$ and $\beta \in T^*_{(q, t)} (\O^\bfw \times T)$,
\[
\left(\mu_* \left(\pi^\#_{(q, t)} (\beta), \, \chi\right)\right) = -\left(\beta, \; \pi^\#_{(q, t)}(\mu^*(\tilde{\chi}))\right) = -(\beta, \,  \sigma(\chi^\# - w \chi^\#))=0.
\]
This shows that ${\rm Im}\left(\pi_{(q, t)}^\#\right) \subset \mu_*^{-1}({\rm Im}(1-w))$ and thus \eqref{eq:need}.
\end{proof}

For $\dbfw =(\dw_1, \dw_2, \ldots, \dw_n) \in \bfw T^n$ and for $a \in T$, we thus need to determine the connected components of the level set 
\begin{equation}\label{eq:tilde-simga-a}
\widetilde{\Sigma}^{\dbfw}_a \, \stackrel{\rm def}{=}\, \mu_{\dbfw}^{-1}(a_\cdot T^w) = \{(q, t) \in \O^\bfw_e \times T: \, t^{-2} \tau_{\dbfw}(q) \in a T^w\}.
\end{equation}
To this end, let ${\rm supp}(\bfw) = \bigcup_{i=1}^n {\rm supp}(w_i)$, and let 
\[
{\rm supp}^o(\bfw)  =\Phi_0\backslash {\rm supp}(\bfw) = 
\bigcap_{i=1}^n {\rm supp}^o(w_i)  =\{\alpha \in \Phi_0: w_i \omega_{\alpha} = \omega_\alpha, \, \forall \, i \in [1, n]\}.
\]
Introduce the sub-torus 
\begin{equation}\label{eq:Tbfw}
\Tbfw = \{t \in T: t^{\omega_\alpha} = 1, \,\forall  \alpha \in {\rm supp}^o(\bfw)\}
\end{equation}
of $T$. Note that $T^w \subset \Tbfw$ and 
that $\dim (\Tbfw) = |{\rm supp}(\bfw)|$.  Let
\begin{equation}\label{eq:delta-00}
\delta_\bfw: \; \O^\bfw_e \times T \longrightarrow T/\Tbfw:\; (q, t) \longmapsto t_\cdot \Tbfw,
\end{equation}
and note that $\delta_\bfw$ has the $T$-equivariance
\begin{equation}\label{eq:delta-equi}
\delta_\bfw (h \cdot (q, t)) = h \delta_\bfw(q, t), \hs h \in T, \, (q, t) \in \O^\bfw_e \times T.
\end{equation}
Let $t_{\dbfw} \in T$ be such that 
$\dw_1 \dw_2 \cdots \dw_n = \overline{w}_1 \overline{w}_2 \cdots \overline{w}_n \, t_{\dbfw} \in w T$. 

\ble{:conn-00}
For any $\dbfw =(\dw_1, \dw_2, \ldots, \dw_n) \in \bfw T^n$ and $a \in T$, the restriction of the map 
\[
\delta_\bfw^2: \; \O^\bfw_e \times T \longrightarrow  T/\Tbfw, \; (q, t) \longmapsto {t^2}_\cdot \Tbfw,
\]
 to $\widetilde{\Sigma}^{\dbfw}_a \subset \O^\bfw_e \times T$ is a constant map:
one has $t^2 \in a^{-1} t_{\dbfw} \Tbfw$ for all $(q, t) \in \widetilde{\Sigma}^{\dbfw}_a$.
\ele

\begin{proof} By \eqref{eq:c-gig},  $\tau_{\overline{\bfw}}(q) \in \Tbfw$ for all $q \in \O^\bfw_e$, so by \reref{:dots},
\begin{equation}\label{eq:tau-tau}
\tau_{\dbfw}(q) = 
t_{\dbfw} \, \tau_{\overline{\bfw}}(q) \in t_{\dbfw} \, \Tbfw, \hs \forall \; q \in \O^\bfw_e.
\end{equation}
It follows that for every $(q, t) \in \widetilde{\Sigma}^{\dbfw}_a$, one has
$t^{-2}  \in  a t_{\dbfw}^{-1}\Tbfw$, and thus 
$t^2 \in a^{-1} t_{\dbfw}\, \Tbfw$.
\end{proof}

\bre{:delta-bfw}
{\rm
By \eqref{eq:c-gig}, $\Delta^{\omega_\alpha}(\overline{w}_1 \overline{w}_2 \cdots \overline{w}_n)= 1$ for every $\alpha \in {\rm supp}^o(\bfw)$. Thus
\[
\Delta^{\omega_\alpha} (\dw_1 \dw_2 \cdots \dw_n) = \Delta^{\omega_\alpha}(t_{\dbfw}), \hs \forall \, \alpha \in {\rm supp}^o(\bfw).
\]
\leref{:conn-00} is then equivalent to saying for every $(q, t) \in \widetilde{\Sigma}^{\dbfw}_a$ one has 
\begin{equation}\label{eq:t2-00}
(t^{\omega_\alpha})^2 =a^{-\omega_\alpha} \Delta^{\omega_\alpha}(\dw_1 \dw_2 \cdots \dw_n),
\hs \forall \; \alpha \in {\rm supp}^o(\bfw).
\end{equation}
\hfill $\diamond$
}
\ere

For $a \in T$, define a {\it level set of $\delta_\bfw$ in $\widetilde{\Sigma}^{\dbfw}_a$} to be any non-empty level set of the map
\[
\delta_\bfw|_{\widetilde{\Sigma}^{\dbfw}_a}: \; 
\widetilde{\Sigma}^{\dbfw}_a \longrightarrow T/\Tbfw.
\]
By \eqref{eq:t2-00}, $\delta_\bfw$ has at most 
$2^{|{\rm supp}^o(\bfw)|}$ level sets 
in $\widetilde{\Sigma}^{\dbfw}_a$, and every connected component of $\widetilde{\Sigma}^{\dbfw}_a$ is contained in one such 
level set. 
Consider the order $2^{|{\rm supp}^o(\bfw)|}$ sub-group 
\[
T^{(2)}_{{\rm supp}^o(\bfw)} = \{\alpha^\vee (\pm 1): \alpha \in {\rm supp}^o(\bfw)\}
\]
of $T^{(2)} = \{t \in T: t^2 = e\}$.
Note that for each $a \in T$, $\widetilde{\Sigma}^{\dbfw}_a$ is $T^{(2)}_{{\rm supp}^o(\bfw)}$-invariant for the
diagonal $T$-action on $\O^\bfw_e \times T$.

\ble{:level-sets-11} For any $a \in T$, there are precisely $2^{|{\rm supp}^o(\bfw)|}$  level sets of
$\delta_\bfw$ in $\widetilde{\Sigma}^{\dbfw}_a$, each pair mutually  
isomorphic by the action of a unique element in $T^{(2)}_{{\rm supp}^o(\bfw)}$. 
\ele

\begin{proof}
By the $T$-equivariance of 
$\delta_\bfw: \O^\bfw_e \times T \to T/\Tbfw$ in \eqref{eq:delta-equi}, 
the group $T^{(2)}_{{\rm supp}^o(\bfw)}$ 
acts freely and transitively on the collection of all level sets of $\delta_\bfw$ in $\widetilde{\Sigma}^{\dbfw}_a$.
\end{proof}

Take now the special representative $\overline{\bfw} = (\overline{w}_1, \overline{w}_2, \ldots, \overline{w}_n) \in \bfw T^n$, and let
\begin{equation}\label{eq:Sigma-bfw-bar}
\Sigma^{\overline{\bfw}} = 
\{(q, t) \in \O^\bfw_e \times T: \, t^{-2} \tau_{\overline{\bfw}}(q) \in T^w, \, t \in \Tbfw\} \subset \widetilde{\Sigma}^{\overline{\bfw}}_e.
\end{equation} 
Pick any $q_0 \in \O^\bfw_e$. Then $\tau_{\overline{\bfw}}(q_0) \in \Tbfw$ by \eqref{eq:c-gig}.  Pick any $t_0 \in \Tbfw$ such that $t_0^2 =  
\tau_{\overline{\bfw}}(q_0)$. Then $(q_0, t_0) \in \Sigma^{\overline{\bfw}}$, showing that $\Sigma^{\overline{\bfw}}\neq \emptyset$.
Thus $\Sigma^{\overline{\bfw}}$ is a 
level set of $\delta_\bfw$ in  $\widetilde{\Sigma}^{\overline{\bfw}}_e$.

\bthm{:conn-0}
For any $\bfw \in W^n$, the sub-variety $\Sigma^{\overline{\bfw}}$ of $\O^\bfw_e \times T$ is connected.
\ethm

\thmref{:conn-0} will be proved in $\S$\ref{ss:proof-conn}. In the rest of $\S$\ref{ss:leaves-D}, we use 
\thmref{:conn-0} to describe all the symplectic leaves of $(\O^\bfw_e \times T, \pi_n\bowtie 0)$
for every $\bfw \in W^n$.

\bthm{:leaves-0}
For any $\dbfw \in \bfw T^n$ and any $a \in T$, the $2^{|{\rm supp}^o(\bfw)|}$ level sets of $\delta_\bfw$ in
$\widetilde{\Sigma}^{\dbfw}_{a}$ are precisely all the connected components of $\widetilde{\Sigma}^{\dbfw}_{a}$ and are thus also all the
symplectic leaves of $(\O^\bfw_e \times T, \pi_n\bowtie 0)$ contained in $\widetilde{\Sigma}^{\dbfw}_{a}$. 
\ethm

\begin{proof} Consider first when $\dbfw = \overline{\bfw}$.  The level sets of $\delta_\bfw$ in
$\widetilde{\Sigma}^{\overline{\bfw}}_{e}$, being isomorphic to $\Sigma^{\overline{\bfw}}
\subset \widetilde{\Sigma}^{\overline{\bfw}}_{e}$ by \leref{:level-sets-11}, are connected 
by \thmref{:conn-0}, and since they are both open and closed, they are all the connected components of $\widetilde{\Sigma}^{\overline{\bfw}}_{e}$.
For any $a \in T$, choose any $h \in T$ such that $h^{-2} = a$. Then $h \cdot \widetilde{\Sigma}^{\overline{\bfw}}_{e} 
=\widetilde{\Sigma}^{\overline{\bfw}}_{a}$. By the $T$-equivariance of $\delta_\bfw$ in \eqref{eq:delta-equi}, the level sets of $\delta_\bfw$ in 
$\widetilde{\Sigma}^{\overline{\bfw}}_{a}$ are in bijection with the level sets of $\delta_\bfw$ in $\widetilde{\Sigma}^{\overline{\bfw}}_{e}$ by the
action of $h$ and are thus all the connected components of $\widetilde{\Sigma}^{\overline{\bfw}}_{a}$.
For an arbitrary $\dbfw \in \bfw T^n$, since $\{\widetilde{\Sigma}^{\dbfw}_a: a \in T\} = \{\widetilde{\Sigma}^{\overline{\bfw}}_{a}: a \in T\}$
by \reref{:dots}, the level sets of $\delta_\bfw$ in
$\widetilde{\Sigma}^{\dbfw}_{a}$ are also all the connected components of $\widetilde{\Sigma}^{\dbfw}_{a}$. 

By \prref{:leaves-conn}, the level sets of $\delta_\bfw$ in
$\widetilde{\Sigma}^{\dbfw}_{a}$ are precisely all the
symplectic leaves of $(\O^\bfw_e \times T, \pi_n\bowtie 0)$ contained in $\widetilde{\Sigma}^{\dbfw}_{a}$. 
\end{proof}

We have the following immediate consequence of \thmref{:leaves-0}.

\bco{:level-Ow}
Let $\bfw \in W^n$ and $\dbfw \in \bfw T^n$. The symplectic leaf of $(\O^\bfw_e \times T, \pi_n \bowtie 0)$ 
through any $(q_0, t_0) \in \O^\bfw_e \times T$ consists precisely of all $(q, t) \in \O^\bfw_e \times T$ satisfying 
\[
\mu_{\dbfw} (q, t) = \mu_{\dbfw} (q_0, t_0) \hs \mbox{and} \hs \delta_\bfw(q, t) = \delta_\bfw(q_0, t_0).
\]
\eco

In view of \coref{:level-Ow}, for $\bfw \in W^n$ and $\dbfw \in \bfw T^n$, it is natural to consider the map
\begin{equation}\label{eq:mu-delta}
\mu_{\dbfw} \times \delta_\bfw: \; \O^\bfw_e \times T \longrightarrow T/T^w \times T/\Tbfw, \;
(q, t) \longmapsto (t^{-2} \tau_{\dbfw}(q)_\cdot T^w, \, t_\cdot \Tbfw).
\end{equation}
Let $X^{\dbfw} =  (\mu_{\dbfw} \times \delta_\bfw)(\O^\bfw_e \times T) \subset T/T^w \times T/\Tbfw$ be the image. 
Symplectic leaves of $(\O^\bfw_e \times T, \, \pi_n \bowtie 0)$ are then 
precisely all the level sets of the now surjective map
\[
\mu_{\dbfw} \times \delta_\bfw:\; \O^\bfw_e \times T \longrightarrow X^{\dbfw}.
\]
To characterize $X^{\dbfw}$, note that $\mu_{\dbfw} \times \delta_\bfw$ is $T$-equivariant, where
$T$ acts on $T/T^w \times T/\Tbfw$ by 
\[
h \cdot (a_\cdot T^w, \, a'{}_\cdot \Tbfw) = (h^{-2}a_\cdot T^w, \, ha'{}_\cdot \Tbfw), \hs a, a' \in T.
\]
Note that $T$ has the same the stabilizer subgroup every point in $T/T^w \times T/\Tbfw$, which is
\begin{equation}\label{eq:Stab}
{\rm Stab}^{\bfw} =\{h \in T: \, h \in \Tbfw, \, h^2 \in T^w\} \subset T.
\end{equation}
Let $p_1: T/T^w \times T/\Tbfw \to T/T^w$ be the projection to the first factor, and note that $p_1$ is $T$-equivariant, where $T$ acts on
$T/T^w$ by $h \cdot (a_\cdot T^w) = h^{-2}a_\cdot T^w$ for $h, a \in T$.
Recall again that $t_{\dbfw} \in T$ is such that
$\dw_1 \dw_2 \cdots \dw_n = \overline{w}_1 \overline{w}_2 \cdots \overline{w}_n \, t_{\dbfw}$. 

\bthm{:X-w}
For any $\bfw \in W^n$ and $\dbfw \in \bfw T^n$, one has
\begin{equation}\label{eq:X-w}
X^{\dbfw} = \{(a_\cdot T^w, a'{}_\cdot \Tbfw) \in T/T^w \times T/\Tbfw: \; a (a')^2 \in t_{\dbfw} \, \Tbfw\}.
\end{equation}
Moreover, $X^{\dbfw}$ is a single $T$-orbit in $T/T^w \times T/\Tbfw$ and is thus smooth and isomorphic to $T/{\rm Stab^\bfw}$. 
The subgroup $T^{(2)}_{{\rm supp}^o(\bfw)}$ of $T$ acts freely on $X^{\dbfw}$, and
the restriction of $p_1$ to $X^{\dbfw}$ gives a covering map $p_1: X^{\dbfw} \to T/T^w$
whose fibers are orbits of
$T^{(2)}_{{\rm supp}^o(\bfw)}$ in $X^{\dbfw}$.
\ethm

\begin{proof}
If $a, a' \in T$ are such that $(a_\cdot T^w, a'{}_\cdot \Tbfw) \in X^{\dbfw}$, it follows from \eqref{eq:mu-delta} that 
 $a (a')^2 \in t_{\dbfw} \, \Tbfw$. Conversely,
suppose that $a, a' \in T$ are such that $a (a')^2 \in t_{\dbfw} \, \Tbfw$. Let $x \in \Tbfw$ be such that $a (a')^2 = t_{\dbfw}x^2$. Let
$\sqrt{t_{\dbfw}}$ be any element in $T$ such that $\sqrt{t_{\dbfw}}^2 = t_{\dbfw}$. Then 
\[
(a_\cdot T^w, a'{}_\cdot \Tbfw) = h\cdot (e_\cdot T^w, \; \sqrt{t_{\dbfw}}{}_\cdot \Tbfw) \in X^{\dbfw},
\]
where $h = a' (\sqrt{t_{\dbfw}})^{-1} x^{-1} \in T$. Thus $X^{\dbfw}$ is given as in \eqref{eq:X-w}.

As $T$ acts transitively on the set of all symplectic leaves of $(\O^\bfw_e \times T, \, \pi_n \bowtie 0)$, the subset
$X^{\dbfw}$ of  $T/T^w \times T/\Tbfw$ is a single $T$-orbit and is thus also smooth. One checks directly from the definitions
that the stabilizer subgroup of $T$ at every point in $X^{\dbfw}$ is ${\rm Stab}^\bfw$. The map $p_1: X^{\dbfw} \to T/T^w$, being is $T$-equivariant,
is thus surjective. For any fixed $a_\cdot T^w \in T/T^w$, one has 
\[
p_1^{-1}(a_\cdot T^w) = \{(a_\cdot T^w, \, a'_\cdot \Tbfw): \; a' \in T, \, (a'_\cdot \Tbfw)^2 = a^{-1}t_{\dbfw} {}_\cdot\Tbfw\}
\]
which is precisely an orbit of $T^{(2)}_{{\rm supp}^o(\bfw)}$ in $X^{\dbfw}$. 
\end{proof}

\bco{:X-w} 
For any $\bfw \in W^n$ and $\dbfw \in \bfw T^n$, the map 
$\mu_{\dbfw} \times \delta_\bfw:  \O^\bfw_e \times T \to X^{\dbfw}$ 
 is a surjective submersion whose level sets are precisely all the symplectic leaves of $(\O^\bfw_e \times T, \pi_n\bowtie 0)$.
\eco

\begin{proof}
As $p_1 \circ (\mu_{\dbfw} \times \delta_\bfw): \O^\bfw_e \times T \to
T/T^w$ is a submersion, and as $p_1: X^{\dbfw} \to T/T^w$ is a covering map, 
$\mu_{\dbfw} \times \delta_\bfw: \O^\bfw_e \times T \to X^{\dbfw}$ is also a submersion.
\end{proof}

For $\dbfw \in \bfw T^n$, let again $\sqrt{t_{\dbfw}}$ be any element in $T$ such that
$(\sqrt{t_{\dbfw}})^2 = t_{\dbfw}$, and set
\begin{equation}\label{eq:E-dbfw-0}
\Sigma^{{\dbfw}} \, \stackrel{{\rm def}}{=}\,  
\{(q, t) \in \O^w_e \times T: \; \; \,t^{-2} \tau_{{\dbfw}}(q) \in T^w, \; t  \in \sqrt{t_{\dbfw}} \,\Tbfw\},
\end{equation}
which is a level set of $\delta_\bfw$ in $\widetilde{\Sigma}^{\dbfw}_{a}$ and thus a symplectic leaf of 
$(\O^\bfw_e \times T, \, \pi_n \bowtie 0)$. Note that 
\begin{equation}\label{eq:SSigma}
\Sigma^{{\dbfw}} = \sqrt{t_{\dbfw}}  \cdot \Sigma^{\overline{\bfw}}.
\end{equation}
We now have the following alternative reformulation of \thmref{:leaves-0}.

\bthm{:leaves-Ow}
For any $\bfw \in W^n$ and $\dbfw \in \bfw T^n$, symplectic leaves of $(\O^\bfw_e \times T, \, \pi_n \bowtie 0)$ are precisely all the sub-varieties of
$\O^\bfw_e \times T$ of the form
\[
h \cdot \Sigma^{{\dbfw}}=  
\{(q, t) \in \O^w_e \times T: \; \; \,t^{-2} \tau_{{\dbfw}}(q) \in h^{-2}T^w, \; t  \in h \sqrt{t_{\dbfw}} \,\Tbfw\},
\]
where $h \in T$.  For  $h_1, h_2 \in T$, 
$h_1 \cdot \Sigma^{{\dbfw}} = h_2 \cdot \Sigma^{{\dbfw}}$ if and only if $h_1^{-1}h_2 \in {\rm Stab}^{\bfw}$ given in \eqref{eq:Stab}.
\ethm

\bre{:Stab}
{\rm
Note that for any  symplectic leaf $\Sigma$ 
of $(\O^\bfw_e \times T, \, \pi_n \bowtie 0)$ and for any $h \in T$, $h \cdot \Sigma = \Sigma$ if and only if $h \in {\rm Stab}^{\bfw}$.
For this reason, we call ${\rm Stab}^{\bfw}$ is the {\it leaf-stabilizer of $T$} in $(\O^\bfw_e \times T, \pi_b\bowtie 0)$.
\hfill $\diamond$
}
\ere

We have already seen in \prref{:leaves-conn} that every symplectic leaf of $(\O^\bfw_e \times T, \pi_n\bowtie 0)$ has dimension equal to
$l(\bfw) + \dim ({\rm Im} (1-w)) = l(\bfw) + \dim(T^w)$. We now show that $\Sigma^{{\dbfw}}$ is a 
$2^{|{\rm supp}(\bfw)|}$-to-$1$ cover of $\O^\bfw_e \times T^w$.  
To this end, let $T$ act on $\O^\bfw_e \times T$ by
\[
h \circ (q, t) = (q, ht), \hs h \in T, \, (q, t) \in \O^\bfw \times T.
\]
Consider $T^{(2)}_{{\rm supp}(\bfw)} =\{\alpha^\vee(\pm 1): \alpha \in {\rm supp}(\bfw)\} = \{h \in \Tbfw: h^2 = e\}$,
a group of order $2^{|{\rm supp}(\bfw)|}$. It follows from the definitions that $\Sigma^{{\dbfw}}$ is invariant under the action of 
$T^{(2)}_{{\rm supp}(\bfw)}$. 

\bpr{:cover}
For any $\dbfw \in \bfw T^n$, the map
\[
\Sigma^{\dbfw} \longrightarrow \O^\bfw_e \times T^w, \; (q, \,t) \longmapsto (q, \, t^{-2} \tau_{\dbfw}(q)),
\]
 is a covering map whose fibers are the orbits of $T^{(2)}_{{\rm supp}(\bfw)}$ in $\Sigma^{\dbfw}$.
\epr

\begin{proof} Let $(q, t') \in \O^\bfw_e \times T^w$. Writing $t \in \sqrt{t_{\dbfw}} \, \Tbfw$ as $t = \sqrt{t_{\dbfw}} \, x$ for $x \in \Tbfw$, then
$t^{-2}\tau_{\dbfw}(q) = t'$ if and only if
$x^2 = (t')^{-1} \tau_{\overline{\bfw}}(q)$. 
By \eqref{eq:c-gig}, $\tau_{\overline{\bfw}}(q) \in \Tbfw$. The equation  $x^2 = (t')^{-1} \tau_{\overline{\bfw}}(q)$,
regarded as one in $\Tbfw$,
then has exactly
$2^{|{\rm supp}(\bfw)|}$ solutions, consisting of a single $T^{(2)}_{{\rm supp}(\bfw)}$ coset in $\Tbfw$. 
\end{proof}

\bex{:w-e}
{\rm
Consider the special case 
when $\bfw = (w_1, \ldots, w_n) \in W^n$ is such that $w = w_1w_2 \cdots w_n = e \in W$. Then $T^w = \{e\}$. Assume also that 
$\bfw = (\dw_1, \ldots, \dw_n) \in \bfw T^n$ satisfies $\dw_1 \dw_2 \cdots \dw_n = e \in G$. Then element
$t_{\dbfw} \in T$ given by 
$e = \overline{w}_1 \cdots \overline{w}_n \, t_{\dbfw}$ lies in $\Tbfw$, and by \eqref{eq:tau-tau}
the image of $\tau_{\dbfw}: \O^\bfw_e \to T$ also lies in $\Tbfw$. 
Thus the symplectic leaf of
$(\O^\bfw_e \times T, \pi_n \bowtie 0)$ through the point $([\dw_1, \ldots, \dw_n]_{\sF_n}, e) \in \O^\bfw_e \times T$ is  
given by 
\[
\Sigma^{\dbfw} = \{(q, t) \in \O^\bfw_e \times \Tbfw: \; t^2= \tau_{\dbfw}(q)\}. 
\]
By \prref{:cover}, the projection $\Sigma^{\dbfw} \to \O^\bfw_e, (q, t) \mapsto q$, is a $2^{|{\rm supp}(\bfw)|}$-to-$1$
covering map.
\hfill $\diamond$
}
\eex

Recall that $\Gamma_n =\{[g_1, g_2, \ldots, g_n]_{\stF_n}: g_1g_2 \cdots g_n \in B_-\}$, with the $T$-action 
\[
t \cdot [g_1, g_2, \ldots, g_n]_{\stF_n} = [tg_1, g_2, \ldots, g_n]_{\stF_n}. 
\]
By \prref{:T-leaves-all}, $T$-leaves of $(\Gamma_n, \tpi_n)$ are precisely all the $\Gamma^\bfw$'s as $\bfw$ runs over $W^n$, where
\[
\Gamma^\bfw = (B \bfw B) \cap \Gamma_n.
\]
In the remainder of $\S$\ref{ss:leaves-D}, we determine symplectic leaves in $(\Gamma^\bfw, \tpi_n)$ for every $\bfw \in W^n$, which is enough for
the discussion in $\S$\ref{s:leaves} on configuration symplectic groupoids. In fact only the cases of 
$\bfw = (\bfu, \bfu^{-1})$ for $\bfu \in W^m$ are needed in $\S$\ref{s:leaves}. Symplectic leaves in 
$(G_{m, n}, \tpi_{m, n})$ and $(\tF_n, \tpi_n)$ for all integers $m, n \geq 1$ are determined in $\S$\ref{ss:Guv-leaves} and
$\S$\ref{ss:leaves-tFn}. 

Let $\bfw = (w_1, \ldots, w_n) \in W^n$ and choose any $\dbfw = (\dot{w}_1, \ldots, \dot{w}_n)\in \bfw T^n$. 
With an arbitrary $\gamma \in \Gamma^\bfw$ uniquely written as (see \eqref{eq:CuB})
\begin{equation}\label{eq:gamma-w}
\gamma = [{c}_1, \ldots, {c}_{n-1}, {c}_n b]_{\stF_n}, \hs \mbox{where}\;\;
(c_1, \ldots, c_n) \in C_{\dbfw}, \; b \in B, \; c_1\cdots c_{n-1}c_n b \in B_-,
\end{equation}
and with $b_- = 
c_1\cdots c_{n-1}c_n b \in B_-$, define
\begin{equation}\label{eq:beta-w}
\beta_{{\dbfw}}: \; \Gamma^\bfw \longrightarrow T/T^w \times T/\Tbfw, \;
\gamma \longmapsto \left([b]_0[b_-]_0 {}_\cdot T^w, \; [b_-]_0{}_\cdot \Tbfw\right).
\end{equation}
Let $\sqrt{t_{\dbfw}} \in T$ be as in \thmref{:leaves-Ow}, and let 
\begin{equation}\label{eq:Lambda-w}
\Lambda^{\dbfw} = \{\gamma \in \Gamma^\bfw: \; [b]_0[b_-]_0 \in T^w, \, [b_-]_0 \in \sqrt{t_{\dbfw}} \,\Tbfw\}
= \beta_{{\dbfw}}^{-1} \left(e_\cdot T^w, \, \sqrt{t_{\dbfw}}{}_\cdot \Tbfw\right).
\end{equation}

\bthm{:Gaw-leaves}
Let $\bfw \in W^n$ and $\dbfw \in \bfw T^n$. 

1)  Symplectic leaves of 
$(\Gamma^\bfw, \tpi_n)$ are precisely all the non-empty level sets of the map 
$\beta_{\dbfw}$ and all have dimension equal to $l(\bfw) + \dim ({\rm Im}(1-w))$. 

2) Alternatively, symplectic leaves of $(\Gamma^\bfw, \tpi_n)$ are all the sub-varieties of $\Gamma^\bfw$ of the form
\[
h \cdot \Lambda^{\dbfw} = \{\gamma \in \Gamma^\bfw: \; [b]_0[b_-]_0 \in h^2 T^w, \, [b_-]_0 \in h \sqrt{t_{\dbfw}} \,\Tbfw\},
\]
where $h \in T$ and $\gamma \in \Gamma^\bfw$ is written as in \eqref{eq:gamma-w}. 
\ethm

\begin{proof}
Under the Poisson isomorphism $J_n: (\Gamma^\bfw, \tpi_n) \to (\O^\bfw_e \times T, \pi_n \bowtie 0)$ in \eqref{eq:Jn}, one 
has
\[
J_n(\gamma) = ([{c}_1, \ldots, {c}_{n-1}, {c}_n]_{\sF_n}, \; [b_-]_0),
\]
where $\gamma$ is as in \eqref{eq:gamma-w}. Thus
\begin{equation}\label{eq:J-bb}
(\mu_{\dbfw} \times \delta_\bfw)(J_n(\gamma)) = ({[b]_0[b_-]_0)^{-1}}_\cdot T^w, \; [b_-]_0{}_\cdot \Tbfw),
\end{equation}
and  $J_n(\Lambda^{\dbfw}) = \Sigma^{\dbfw}$.
\thmref{:Gaw-leaves} now follows from \coref{:level-Ow} and \thmref{:leaves-Ow}.
\end{proof}

\bre{:beta-w}
{\rm
By  \eqref{eq:J-bb} and \thmref{:X-w}, 
the image of $\beta_{\dbfw}: \Gamma^\bfw \to T/T^w \times T/\Tbfw$ is
\[
Y^{\dbfw} = \{(a_\cdot T^w, a'{}_\cdot \Tbfw) \in T/T^w \times T/\Tbfw: \; a^{-1} (a')^2 \in t_{\dbfw} \, \Tbfw\},
\]
and 
symplectic leaves of $(\Gamma^\bfw, \tpi_n)$ are precisely all the level sets of the surjective submersion
\[
\beta_{{\dbfw}}: \; \Gamma^\bfw \longrightarrow Y^{\dbfw}.
\]
Moreover, $\beta_{\dbfw}$ is $T$-equivariant, where $T$ acts on $T/T^w \times T/\Tbfw$ by 
\[
h \cdot (a_\cdot T^w, \; a'{}_\cdot \Tbfw) = (h^2 a_\cdot T^w, \; ha'{}_\cdot \Tbfw), \hs a, a' \in T,
\]
and $Y^{\dbfw} \subset T/T^w \times T/\Tbfw$ is a single $T$-orbit with ${\rm Stab}^\bfw$ in \eqref{eq:Stab} as the stabilizer subgroup.
\hfill $\diamond$
}
\ere

\bex{:uu-leaf}
{\rm
Consider the special case when $\bfw = (\bfu, \bfu^{-1}) \in W^{2n}$, where $\bfu = (u_1, \ldots, u_n) \in W^n$ and
$\bfu^{-1} = (u_n^{-1}, \ldots, u_1^{-1})$. Choose any $\dbfu = (\dot{u}_1, \ldots, \dot{u}_n) \in \bfu T^n$, so we have the point
\[
[\bfu, \bfu^{-1}]_{\stF_{2n}} \, \stackrel{\rm def}{=}\, [\dot{u}_1, \ldots, \dot{u}_n, \, \dot{u}_n^{-1}, \ldots, 
\dot{u}_1^{-1}]_{\stF_{2n}} \in \Gauu
\]
which is in fact independent of the choice of the representative $\dbfu$ for $\bfu$.  
Let 
\[
\dbfw = (\dot{u}_1, \,\ldots, \,\dot{u}_n, \,\dot{u}_n^{-1}, \,\ldots, \,\dot{u}_1^{-1}) \in \bfw T^{2n}. 
\]
By \thmref{:Gaw-leaves} and by \exref{:w-e}, the 
symplectic leaf of $(\Gauu, \tpi_{2n})$ through 
$[\bfu, \bfu^{-1}]_{\stF_{2n}}$ 
is $\Lambda^{\dbfw}$ in \eqref{eq:Lambda-w}, i.e. the sub-variety of $\Gauu$ consisting of all 
$\gamma =[c_1, \ldots, c_{2n-1}, c_{2n}b]_{\stF_{2n}}$, 
where $(c_1, \ldots, c_{2n-1}, c_{2n}) \in C_{\dbfw}$ and $b \in B$ such that 
\[
c_1\cdots c_{2n-1} c_{2n}b = b_- \in B_- \hs \mbox{and} \hs
[b]_0[b_-]_0=e, \; [b_-]_0 \in \Tbfu.
\]
This example is used in $\S$\ref{ss:Lambda-poid}.
\hfill $\diamond$
}
\eex

\subsection{Bott-Samelson coordinates and Lusztig toric charts on $\O^\bfw$}\label{ss:Lusz-charts}
To prepare for the proof of  \thmref{:conn-0} in  $\S$\ref{ss:proof-conn}, we recall in this section some toric charts on 
generalized Schubert cells  which are of interests of their own (see \reref{:cluster}).

Consider again an arbitrary 
$\bfw = (w_1, w_2, \ldots, w_n) \in W^n$. For each $i \in [1, n]$, choose a reduced word $\bfw_i$ of $w_i$ and regard
$\bfw_i$ as in $W^{l(w_i)}$. One then has the concatenation
\begin{equation}\label{eq:tilde-bfw}
\widetilde{\bfw} =(\bfw_1, \, \bfw_2, \, \ldots, \, \bfw_n) =  (s_{\alpha_1}, \; s_{\alpha_2}, \; \ldots, \; s_{\alpha_{n'}}) \in W^{n'},
\end{equation}
where $n' =l(\bfw) =  l(w_1) + l(w_2) + \cdots + l(w_n)$, and $\alpha_j \in \Phi_0$ for each $j \in [1, n']$. 
Using the parametrization of $\O^\bfw$ by
$C_{\overline{\bfw}} = C_{\overline{w}_1}  \times \cdots \times C_{\overline{w}_n}$ in \eqref{eq:sigma-w} and 
the parametrization $g_{\bfw_i}: \CC^{l(w_i)} \to C_{\overline{w}_i}$  in \eqref{eq:g-bfu} for each 
$i \in [1, n]$, one obtains the isomorphism
\begin{equation}\label{eq:q-w}
q_{\widetilde{\bfw}}: \; \CC^{n'} \longrightarrow \O^\bfw, \; 
q_{\widetilde{\bfw}} (z_1, z_2, \ldots, z_{n'}) = [g_{\bfw_1}(z_1, \ldots, z_{l_1}), \, \ldots, \, 
g_{\bfw_n}(z_{l_{n-1}+1}, \ldots, z_{n'})]_{\tF_{n}},
\end{equation}
where $l_i = l(w_1) + \cdots l(w_i)$ for $i \in [1, n]$ and $n' = l_n$.
Following \cite{Elek-Lu:BS}, we call $q_{\widetilde{\bfw}}: \CC^{n'} \to \O^\bfw$ the {\it Bott-Samelson parametrization of $\O^\bfw$ defined by $\wbfw$}, and
the resulting coordinates $(z_1, z_2, \ldots, z_{n'})$ on $\O^\bfw$ the
{\it Bott-Samelson coordinates} on $\O^\bfw$ defined by $\widetilde{\bfw}$.

We now use $\widetilde{\bfw}$ to define an open toric chart on $\O^\bfw$. 

\ble{:m-bfw}
For any $\bfw = (w_1, \ldots,w_n) \in W^n$, the map
\begin{align}\label{eq:delta-bfw-00}
\varsigma_\bfw: \; (N_- \cap Bw_1B) \times \cdots \times (N_- \cap Bw_n B) \longrightarrow \O^\bfw,\;\;
(m_1, \ldots, m_n) \longmapsto [m_1, \, \ldots, m_n]_{\sF_n}
\end{align}
induces a biregular isomorphism between $(N^- \cap Bu_1B) \times \cdots \times (N_- \cap Bw_n B) $ and the Zariski open subset
$(\O^\bfw)_0$ of $\O^\bfw$ given by
\[
(\O^\bfw)_0 = \{[c_1, c_2, \ldots, c_n]_{\sF_n}: c_i \in C_{\overline{w}_i}, \, c_1\cdots c_i \in B^-B,\, \forall \,\; i \in [1, n]\}.
\]
\ele

\begin{proof}
Let $(m_1, \ldots, m_n) \in (N_- \cap Bw_1B) \times \cdots \times (N_- \cap Bw_n B)$, and let $c_i \in C_{\overline{w}_i}$ and $b_i \in B$
for $i \in [1, n]$ be defined by 
\[
m_1 = c_1 b_1, \;\;\; b_1m_2 = c_2 b_2, \;\;\; \ldots, \;\; \;b_{n-1}m_n = c_n b_n.
\]
Then $\varsigma_\bfw(m_1, \ldots, m_n) = [c_1, \, c_2, \, \ldots, \, c_n]_{\sF_n}$ by definition. As
$c_1 \cdots c_i= m_1\cdots m_i b_i^{-1} \in B^-B$ for each $i \in [1, n]$,
one has $\varsigma_\bfw(m_1, \ldots, m_n) \in (\O^\bfw)_0$. Furthermore, for any $(c_1, \ldots, c_n) \in C_{\dbfw}$
such that $[c_1, \ldots, c_n]_{\sF_n} \in (\O^\bfw)_0$, 
let $m_i \in N_-$, $i \in [1, n]$,  be given by
\[
m_1 = [c_1]_-, \hs m_1m_2 = [c_1c_2]_-, \hs \ldots, \hs
m_1 m_2\cdots m_n = [c_1c_2 \cdots c_n]_-.
\]
Then $(m_1, \ldots, m_n)$ is the unique element in $(N_- \cap Bw_1B) \times \cdots \times (N_- \cap Bw_n B)$ such that
$\varsigma_\bfw(m_1, m_2, \ldots, m_r)  = [c_1, c_2, \ldots, c_n]_{\sF_n}$.  Thus $\varsigma_\bfw$ induces the biregular isomorphism as described.
\end{proof}

Consider now the concatenation $\widetilde{\bfw}$ in \eqref{eq:tilde-bfw}. Combining $\varsigma_{\bfw}$ in \eqref{eq:delta-bfw-00}
with  the open embeddings
$m_{\bfw_i}: (\CC^\times)^{l(w_i)} \to N_- \cap Bw_iB$ in \eqref{eq:m-bfu} for $i \in [1, n]$, one obtains the open embedding
$\sigma_{\widetilde{\bfw}}: (\CC^\times)^{n'} \to \O^\bfw_e$ given by
\begin{equation}\label{eq:sigma-w-0}
\sigma_{\widetilde{\bfw}}(\varepsilon_1, \varepsilon_2, \ldots, \varepsilon_{n'}) 
= [m_{\bfw_1}(\varepsilon_1, \ldots, \varepsilon_{l_1}),  \ldots, 
m_{\bfw_n}(\varepsilon_{l_{n-1}+1}, \ldots, \varepsilon_{n'})]_{\sF_n}.
\end{equation}
We call $\sigma_{\widetilde{\bfw}}:  (\CC^\times)^{n'} \to \O^\bfw_e$ the {\it Lusztig toric chart on $\O^\bfw$ defined by $\widetilde{\bfw}$}.

We now solve the {\it inverse parameter problem} for the open embedding $\sigma_{\widetilde{\bfw}}: (\CC^\times)^{n'} \to\O^\bfw_e$.
More specifically, for each $j \in [1, n']$, we will express $\varepsilon_j$, regarded as a rational function on $\O^\bfw$ through 
\eqref{eq:sigma-w-0}, as a monomial of certain regular functions on $\O^\bfw$. To this end, using 
the Bott-Samelson parametrization
$q_{\widetilde{\bfw}}: \CC^{n'} \to \O^\bfw$ in \eqref{eq:q-w}, define  
regular functions $\phi_{\wbfw, j}$ on $\O^\bfw$ by 
\begin{equation}\label{eq:phi-w-k}
\phi_{\wbfw, j}(q_\wbfw(z_1, z_2, \ldots, z_{n'})) =
\Delta^{\omega_{\alpha_j}}(x_{\alpha_1}(z_1) \overline{s}_{\alpha_1} 
x_{\alpha_2}(z_2) \overline{s}_{\alpha_2} \cdots  x_{\alpha_j}(z_j) \overline{s}_{\alpha_j}), \hs j \in [1, n'].
\end{equation}
 Set  $\O^\bfw_{\phi_{\wbfw} \neq 0} = \{q \in \O^\bfw: \phi_{\wbfw, 1}(q)  \phi_{\wbfw, 2}(q) \cdots \phi_{\wbfw, n'}(q)\neq 0\} \subset \O^\bfw$.

\ble{:toric} One has  $\sigma_\wbfw((\CC^\times)^{n'}) = \O^\bfw_{\phi_\wbfw \neq 0}$. 
Moreover, for $\varepsilon = (\varepsilon_1, \ldots, \varepsilon_{n'}) \in (\CC^\times)^{n'}$, 
\begin{equation}\label{eq:epsilon-phi}
\varepsilon_j =  
\left(\phi_{\wbfw, 1}^{r_{1, j}} \phi_{\wbfw, 2}^{r_{2, j}} \cdots \phi_{\wbfw,  j-1}^{r_{j-1, j}}\phi_{\wbfw, j}^{-1}\right)(\sigma_\wbfw(\varepsilon)), 
\hs j \in [1, n'],
\end{equation}
where for $j \in [1, n']$ and $i \in [1, j-1]$, 
\[
r_{i, j} = \begin{cases} 0, & \hs \mbox{if} \;\; \alpha_i \in \{\alpha_{i+1}, \ldots, \alpha_{j-1}\},\\
-\alpha_j(\alpha_i^\vee), & \hs \mbox{if} \;\; \alpha_i \notin \{\alpha_{i+1}, \ldots, \alpha_{j-1}\},\; \alpha_i \neq \alpha_j,\\
-1, & \hs \mbox{if} \;\; \alpha_i \notin \{\alpha_{i+1}, \ldots, \alpha_{j-1}\},\; \alpha_i = \alpha_j.\end{cases}
\]
\ele

\begin{proof} 
For notational simplicity, we set $\phi_j = \phi_{\wbfw, j}$ for $j \in [1, n']$.

Consider the generalized Schubert cell $\O^{\wbfw} \subset F_{n'}$ defined by $\wbfw \in W^{n'}$.
 As each $\bfw_i$, for $i \in [1, n]$, is a reduced word of $w_i$,  we have the  isomorphism
${\bf m}:  \O^{\wbfw} \to \O^\bfw$ given by
\[
{\bf m}\left([g_1, \ldots, g_{n'}]_{\sF_{n'}}\right) = [g_1\cdots g_{l_1}, \, g_{l_1+1} \cdots g_{l_2}, \, \ldots, \, g_{l_{n-1}+1} \cdots g_{n'}]_{\sF_{n}},
\]
where $g_j \in Bs_{\al_j} B$ for $j \in [1, n']$.
Let $\varepsilon = (\varepsilon_1,\ldots, \varepsilon_{n'}) \in (\CC^\times)^{n'}$ and $z = (z_1, \ldots, z_{n'}) \in \CC^{n'}$ 
be such that $\sigma_\wbfw(\varepsilon) = q_\wbfw(z)$, i.e., 
\[
[x_{-\al_1}(\varepsilon_1), \, \ldots, \, x_{-\al_{n'}}(\varepsilon_{n'})]_{\sF_{n'}} 
= [x_{\al_1}(z_1)\overline{s}_{\alpha_1}, \,\ldots, \, x_{\al_{n'}}(z_{n'})\overline{s}_{\alpha_{n'}}]_{\sF_{n'}}.
\]
Let $j \in [1, n']$. Then $[x_{-\al_1}(\varepsilon_1), \, \ldots, \, x_{-\al_j}(\varepsilon_j)]_{\sF_j} 
= [x_{\al_1}(z_1)\overline{s}_{\alpha_1}, \,\ldots, \, x_{\al_j}(z_j)\overline{s}_{\alpha_j}]_{\sF_j}$ and thus
\begin{equation}\label{eq-cc-gg}
x_{-\al_1}(\varepsilon_1)\cdots x_{-\al_j}(\varepsilon_j)t_jn_j = x_{\al_1}(z_1)\overline{s}_{\alpha_1} \cdots  x_{\al_j}(z_j)\overline{s}_{\alpha_j}
\end{equation}
for some unique $t_j \in T$ and $n_j \in N$. It follows that
\begin{equation}\label{eq:t-phi}
\phi_j(q_\wbfw(z)) = t_j^{\omega_{\alpha_j}}.
\end{equation}
This shows in particular that $\sigma_\wbfw(\varepsilon) \in \O^\bfw_{\phi_\wbfw \neq 0}$. 
Comparing \eqref{eq-cc-gg} for $j$ and $j-1$, one then has
$x_{-\al_j}(\varepsilon_j) t_jn_j = t_{j-1}n_{j-1}x_{\al_j}(z_j)\overline{s}_{\alpha_j}$.
By \eqref{eq:SL2}, one  has
\[
x_{\al_j} (\varepsilon_j^{-1}) \overline{s}_{\al_j} \al_j^\vee(\varepsilon_j) x_{\al_j}(\varepsilon_j^{-1}) t_jn_j=t_{j-1}n_{j-1}x_{\al_j}(z_j)\overline{s}_{\al_j} 
\in Bs_{\al_j}B.
\]
It follows that $\al_j^\vee(\varepsilon_j) t_j = t_{j-1}^{s_{\al_j}}$, and thus by \eqref{eq:t-phi},
\[
\varepsilon_j  \phi_j(q_\wbfw(z)) = (\al_j^\vee(\varepsilon_j) t_j )^{\omega_{\al_j}} =  t_{j-1}^{s_{\al_j}\omega_{\al_j}}.
\]
Writing $s_{\al_j}\omega_{\al_j} = \sum_{\al \in \Phi_0} (s_{\al_j}\omega_{\al_j}, \al^\vee) \omega_\al$, one then has
\[
\varepsilon_j  \phi_j(q_\wbfw(z))   =  \prod_{\al \in \Phi_0} (t_{j-1}^{\omega_\al})^{(s_{\al_j}\omega_{\al_j}, \,\al^\vee)}
= \prod_{\al \in \Phi_0} (t_{j-1}^{\omega_\al})^{(\omega_{\al_j}-\alpha_j, \al^\vee)}.
\]
If $\al \notin \{\al_1, \ldots,\al_{j-1}\}$, then 
$t_{j-1}^{\omega_\al} = \Delta^{\omega_\al}(x_{\al_1}(z_1)\overline{s}_{\alpha_1} \cdots  x_{\al_{j-1}}(z_{j-1})\overline{s}_{\alpha_{j-1}}) = 1$ by
\eqref{eq:gig}. If $\al \in \{\al_1, \ldots,\al_{j-1}\}$, then $t_{j-1}^{\omega_\al} = \phi_{j_\al}(q_\wbfw(z))$, where
$j_\al = {\rm max}\{i \in [1, j-1]: \al = \al_i\}$. With  $r_{1, j}, \ldots, r_{j-1, j}$ as described, one thus has
\[
\varepsilon_j = \frac{1}{\phi_j(\sigma_\wbfw(\varepsilon))}
(\phi_1(\sigma_\wbfw(\varepsilon)))^{r_{1,j}} (\phi_2(\sigma_\wbfw(\varepsilon)))^{r_{2,j}}
 \cdots (\phi_{j-1}(\sigma_\wbfw(\varepsilon)))^{r_{j-1, j}}.
\]
Furthermore, given $q = q_\wbfw(z) \in \O^\bfw_{\phi_{\wbfw} \neq 0}$, one sees by induction on $j$ that there exist
unique $\varepsilon_j \in \CC^\times$ and
unique $t_j \in T$ and $n_j \in N$ for $j \in [1, n']$ such that \eqref{eq-cc-gg} holds. Thus $\sigma_\wbfw(\varepsilon) = q_\wbfw(z)$. 
This shows that $\sigma_\wbfw((\CC^\times)^n) = \O^\bfw_{\phi_\wbfw \neq 0}$.
This finishes the proof of \leref{:toric}.
\end{proof}

\begin{remark}\label{re-c-phi}
{\rm For $u \in W$ let $N_u = N \cap \overline{u} N_- \overline{u}^{\, -1}$ so that $C_{\overline{u}} = N_u \overline{u}$. 
For $\bfw = (w_1, w_2, \ldots, w_n) \in W^n$, one then has the isomorphism
\begin{equation}\label{eq:nnn}
N_{w_1} \times N_{w_2} \times \cdots \times N_{w_n} \lrw \O^\bfw, \; (g_1, g_2, \ldots, g_n) \longmapsto [g_1\overline{w}_1, g_2 \overline{w}_2, \ldots, 
g_n\overline{w}_n]_{\sF_n}.
\end{equation}
Recall  also from \cite{FZ:double} that one has 
generalized minors 
\[
\Delta_{u\omega_\alpha, v\omega_\alpha}(g) = \Delta^{\omega_\alpha}(\overline{u}^{\, -1} g \overline{v}), \hs g \in G,
\]
where $u, v \in W$ and $\alpha \in \Phi_0$.
Using the parametrization of $\O^\bfw$ in \eqref{eq:nnn}, one can also express the functions $\phi_{\wbfw, j} \in \CC[\O^\bfw]$ as follows:
for $i \in [1, n]$ and $j \in [l_{i-1}+1, l_i]$, 
\[
\phi_{\wbfw, j}([g_1\overline{w}_1, g_2 \overline{w}_2, \ldots, 
g_n\overline{w}_n]_{\sF_n}) = 
\Delta_{\omega_{\alpha_j}, \,s_{l_{i-1}+1} \cdots s_j \omega_{\alpha_j}}([g_1\overline{w}_1 \cdots g_{i-1} \overline{w}_{i-1} g_i),
\]
where we have set $s_k = s_{\alpha_k}$ for $k \in [1, n']$.
In the special case when $n = 1$ so $\bfw = w \in W$ and $\wbfw = (s_{\alpha_1}, \cdots, s_{\alpha_l})$ is a reduced word for $w$,
we can parametrize $\O^w = BwB/B \subset G/B$ by $N_w \to \O^w, g \mapsto g \overline{w}_\cdot B$, 
and \eqref{eq:epsilon-phi} can be rewritten as 
\begin{equation}\label{eq-c-phi-1}
\varepsilon_j = \frac{\prod_{\al\neq \al_j} \Delta_{\omega_{\al}, \, s_{\al_1} \cdots s_{\al_j} \omega_{\al}}(g)^{-\al_j(\al^\vee)}}{
\Delta_{\omega_{\al_j}, \, s_{\al_1} \cdots s_{\al_j} \omega_{\al_j}}(g)
\Delta_{\omega_{\al_j}, \, s_{\al_1} \cdots s_{\al_{j-1}} \omega_{\al_j}}(g)}, \hs j \in [1, l], \, g \in N_w.
\end{equation}
In this case,  the fact that $\sigma_\wbfw((\CC^\times)^l) =
\O^\bfw_{\phi_\wbfw \neq 0}$ also follows from \cite[Proposition 5.2, Corollary 6.6]{MR:para}, and 
\eqref{eq-c-phi-1} has been proved in 
\cite[Theorem 7.1]{MR:para}. Equivalent formulations of \eqref{eq-c-phi-1}
can be found in \cite[Theorem 1.4]{BZ:total} and \cite[Theorem 2.19]{FZ:double}.
\hfill $\diamond$
}
\end{remark}

\bco{:toric}
In the context of \leref{:toric}, one has the isomorphism
\begin{equation}\label{eq:phi-wbfw}
\phi_{\wbfw} = (\phi_{\wbfw, 1}, \ldots, \phi_{\wbfw, n'}): \; \; \O^\bfw_{\phi_{\wbfw} \neq 0} \longrightarrow (\CC^\times)^{n'}.
\end{equation}
\eco

 By \leref{:toric},  $\phi_{\wbfw}: \O^\bfw_{\phi_{\wbfw} \neq 0} \rightarrow (\CC^\times)^{n'}$ is the
inverse of  $\sigma_\wbfw: (\CC^\times)^{n'} \to \O^\bfw_{\phi_\wbfw \neq 0}$ up to an invertible monomial transformation
on $(\CC^\times)^{n'}$.

\bre{:cluster}
{\rm 
The functions $\{\phi_{\wbfw, 1}, \ldots, \phi_{\wbfw, n'}\}$ form an initial cluster for a cluster structure on $\O^\bfw$ defined by 
K. Goodearl and M. Yakimov  using the theory of symmetric Poisson CGLs \cite{GY:Poi-CGL}, an aspect of generalized Schubert cells that
will be explored elsewhere.
\hfill $\diamond$
}
\ere

\subsection{Proof of \thmref{:conn-0}}\label{ss:proof-conn} Let again $\bfw = (w_1,\ldots, w_n) \in W^n$ and 
we  return to the (non-empty) sub-variety $\Sigma^{\overline{\bfw}}$ of $\O^\bfw \times T$ given
in \eqref{eq:Sigma-bfw-bar}, i.e., 
\[
\Sigma^{\overline{\bfw}} = \{(q,\, t) \in \O^\bfw_e \times T: \; t \in \Tbfw, \, t^{-2} \tau_\bfw(q) \in T^w\}.
\]
We now prove \thmref{:conn-0} which states that $\Sigma^{\overline{\bfw}}$ is connected. 

Let $A^\bfw = \{(q, t, t') \in \O^\bfw_e \times \Tbfw \times T^w: t^{-2} \tau_{\overline{\bfw}}(q) = t'\}$. 
Then $A^\bfw \to \Sigma^{\overline{\bfw}}, (q, t, t') \mapsto (q, t)$, is an isomorphism. We thus only need to show that $A^\bfw$ is connected.
Using the fact that for $a_1, a_2 \in T$, $a_1 = a_2$ if and only if $a_1^{\omega_\alpha} = 
a_2^{\omega_\alpha}$ for all $\alpha \in \Phi_0$, one has
\[
A^\bfw = \{(q, t, t') \in  \O^\bfw_e \times \Tbfw \times T^w: \; (\tau_{\overline{\bfw}}(q))^{\omega_\alpha} = (t^2t')^{\omega_\alpha}, \,\forall \, \alpha \in \Phi_0\}.
\]
Since $T^w \subset \Tbfw = \{t \in T: t^{\omega_\alpha} = 1, \, \forall \, \alpha \in {\rm supp}^o(\bfw)\}$ and
$\tau_{\overline{\bfw}}(\O^\bfw_e) \subset \Tbfw$, one has 
\[
A^\bfw = \{(q, t, t') \in \O^\bfw_e \times \Tbfw \times T^w:  \; (\tau_{\overline{\bfw}}(q))^{\omega_\alpha} = 
(t^2t')^{\omega_\alpha}, \,\forall \, \alpha \in {\rm supp}(\bfw)\}.
\]
Choose a reduced word for each $w_i$ and consider again the 
sequence $\wbfw \in W^{n'}$ of simple reflections in \eqref{eq:tilde-bfw},
where $n' = l(\bfw)$, and the sequence $\{\phi_{\wbfw, j}: j \in [1, n']\}$ of
regular functions on $\O^\bfw$ given in \eqref{eq:phi-w-k}. 
For $\alpha \in {\rm supp}(\bfw)$, let $j_\bullet(\alpha) \in [1, n']$ be the maximal $j \in [1, n']$ such that $\alpha = \alpha_j$. By \eqref{eq:c-gig},  
one has $(\tau_{\overline{\bfw}}(q))^{\omega_\alpha} = \phi_{\wbfw, j_\bullet(\alpha)}(q)$ for all $\alpha \in 
{\rm supp}(\bfw)$ and $q \in \O^\bfw_e$. 
On the other hand, since an element $g \in G$ lies in $B_-B$ if and only if $\Delta^{\omega_\alpha}(g) \neq 0$ for all $\alpha \in \Phi_0$, 
one has by \eqref{eq:c-gig} again that 
$\O^\bfw_e = \{q \in \O^\bfw: \; \phi_{\wbfw, j_\bullet(\alpha)}(q) \neq 0,  \,\forall \, \alpha \in {\rm supp}(\bfw)\}\}$. 
It follows that 
\[
A^\bfw = \{(q, t, t') \in \O^\bfw \times \Tbfw \times T^w: \; 
\phi_{\wbfw, j_\bullet(\alpha)}(q) =( t^{2}t')^{\omega_\alpha},\, \forall \,\alpha \in {\rm supp}(\bfw)\}.
\]
Consider the Zariski open subset $A^\bfw_0$ of $A^\bfw$ given by
\[
A^\bfw_0 =A^\bfw \cap (\O^\bfw_{\phi_\wbfw\neq 0} \times \widetilde{T}^\bfw \times T^w).
\]
Let $J =\{j_\bullet(\alpha): \alpha \in {\rm supp}(\bfw)\}$ and $J^o = [1, n'] \backslash J$. Under the isomorphism
\[
\phi_{\wbfw} \times {\rm Id}_{\Tbfw \times T^w}:\;\;
\O^\bfw_{\phi_\wbfw\neq 0} \times \widetilde{T}^\bfw \times T^w \longrightarrow (\CC^\times)^{n'} \times \widetilde{T}^\bfw \times T^w,
\]
the sub-variety $A^\bfw_0 \subset \O^\bfw_{\phi_\wbfw\neq 0} \times \widetilde{T}^\bfw \times T^w$ is then defined by the equations 
\[
\phi_{\wbfw, j_\bullet(\alpha)}(q) =( t^{2}t')^{\omega_\alpha},\, \hs \forall \,\alpha \in {\rm supp}(\bfw)
\]
on $(q, t, t') \in \O^\bfw_{\phi_\wbfw\neq 0} \times \widetilde{T}^\bfw \times T^w$. 
Thus $A^\bfw_0 \cong (\CC^\times)^{|J^o|} \times \widetilde{T}^\bfw \times T^w \cong (\CC^\times)^d$, where
$d:=n' + \dim ({\rm Im}(1-w))$. 
Now the complement of $A^\bfw_0$ in $A^\bfw$ is 
$\bigcup_{j \in J^o} Z_j$, where
\[
Z_j = \{(q, t, t') \in A^\bfw: \; \phi_{\wbfw, j}(q) = 0\}
\]
for $j \in J^o$. Thus $\dim Z_j < d$ for each $j \in J^o$.  As $A^\bfw$
is smooth of dimension $d$, $A^\bfw_0$ is dense in $A^\bfw$.
Since $A^\bfw_0 \cong (\CC^\times)^d$ is connected, it follows that $A^\bfw$ is connected.

This finishes the proof of \thmref{:conn-0}.

\subsection{Symplectic leaves in generalized double Bruhat cells}\label{ss:Guv-leaves}
Assume again that $G$ is connected and simply connected.
For any integers $m, n \geq 1$, we describe in this section all the symplectic leaves of the $T$-Poisson variety $(G_{m, n}, \tpi_{m, n})$ introduced in 
$\S$\ref{ss:Guu} and $\S$\ref{ss:Guv-app}. By \coref{:Tleaf-Gmn}, the $T$-leaves of $(G_{m, n}, \tpi_{m, n})$ are precisely the
generalized double Bruhat cells in $G_{m, n}$. It is thus enough to determine symplectic leaves in all generalized Bruhat cells. We remark that 
symplectic leaves in $(G_{1, 1}, \tpi_{1, 1}) \cong (G, \pist)$ are determined by M. Kogan and A. Zelevinsky in
\cite{k-z:leaves}.

Let thus $\bfu = (u_1, \ldots, u_m) \in W^m$ and $\bfv = (v_1, \ldots, v_n) \in W^n$ be arbitrary. 
By \coref{:JEv-Poi}, one has a  $T$-equivalent
isomorphism $K_{\bfu, \dot{\bfv}}: (G^{\bfu, \bfv}, \tpi_{m, n}) \to (\O^{(\bfu, \bfv^{-1})}_e \times T, \pi_{m+n} \bowtie 0)$ 
for each choice of $\dot{\bfv} \in \bfv T^n$. As symplectic leaves of $(\O^{(\bfu, \bfv^{-1})}_e \times T, \pi_{m+n} \bowtie 0)$
are determined by  \thmref{:leaves-Ow},
we will use  $K_{\bfu, \dot{\bfv}}$ and  \thmref{:leaves-Ow} for a special choice of $\dot{\bfv}$ to determine all the
symplectic leaves of $(G^{\bfu, \bfv}, \tpi_{m, n})$. 

For $w \in W$, let $\overline{\overline{w}} = (\overline{w^{-1}})^{\, -1} \in wT$ (see \cite[$\S$1.4]{FZ:double}).
For representatives of $\bfu$ and $\bfv$, we choose
\[
\overline{\bfu} = (\overline{u}_1, \ldots, \overline{u}_m) \in \bfu T^m \hs \mbox{and} \hs
\overline{\overline{\bfv}} = (\overline{\overline{v}}_1, \, \ldots, \,  \overline{\overline{v}}_n) \in \bfv T^n.
\]
 Recall that
$C_{\overline{{\bfu}}} = C_{{\overline{u}}_1} \times \cdots \times C_{\overline{{u}}_m}$ and
$C_{\overline{\overline{\bfv}}} = C_{\overline{\overline{v}}_1} \times \cdots \times C_{\overline{\overline{v}}_n}$.
Introduce
\begin{equation}\label{eq:G-bar-u-v}
\calG^{\overline{\bfu}, \overline{\overline{{\bfv}}}} = 
\{({c}_1, \ldots, {c}_m, b, b_-, {c}_1^\prime, \ldots, {c}^\prime_n) \in 
C_{\overline{\bfu}} \times B \times B_- \times C_{\overline{\overline{\bfv}}}: \,
{c}_1\cdots {c}_mb = b_- {c}^\prime_1 \cdots {c}^\prime_n\}.
\end{equation}
One can then parametrize $G^{\bfu, \bfv}$ by $\calG^{\overline{\bfu}, \overline{\overline{{\bfv}}}}$ by sending 
$({c}_1, \ldots, {c}_m, b, b_-, {c}_1^\prime, \ldots, {c}^\prime_n) \in \calG^{\overline{\bfu}, \overline{\overline{{\bfv}}}}$
 to 
\begin{equation}\label{eq:bf-g-hat-dot}
{\bf g} = \left([{c}_1, \ldots, {c}_{m-1}, {c}_m b]_{\stF_m}, \, [b_- {c}^\prime_1, {c}^\prime_2, \ldots, {c}^\prime_n]_{\stF_{-n}}\right) \in \Guv.
\end{equation}
Let $u = u_1\cdots u_m \in W$, $v = v_1\cdots v_n \in W$, and recall that $T^{uv^{-1}} = \{a(a^{-1})^{uv^{-1}}: a \in T\}$.
Set 
\[
T^{u, v} = \{(a^{-1})^u a^v: a \in T\} = \{a^v: a \in T^{uv^{-1}}\}. 
\] 
Let again ${\rm supp}^o(\bfu, \bfv) = \{\alpha \in \Phi_0: u_i \omega_\alpha = v_j \omega_\alpha = \omega_\alpha, \forall \; i \in [1, m], j \in [1, n]\}$, and
 set 
\[
\widetilde{T}^{\bfu, \bfv} = \{t \in T: t^{\omega_\alpha} = 1, \; \forall \; \alpha \in {\rm supp}^o(\bfu, \bfv)\}.
\]
With ${\bf g} \in G^{\bfu, \bfv}$ expressed as in \eqref{eq:bf-g-hat-dot}, define now
\begin{equation}\label{eq:chi-1}
\chi_{{\bfu}, {\bfv}}: \;\;  G^{\bfu, \bfv} \longrightarrow (T/T^{u, v}) \times (T/\widetilde{T}^{\bfu, \bfv}), \;\; 
{\bf g} \longmapsto \left([b]_0 [b_-]_0^v {}_\cdot T^{u, v}, \; [b]_0{}_\cdot \widetilde{T}^{\bfu, \bfv}\right).
\end{equation}
Recall that $T$ acts on $G^{\bfu, \bfv} \subset G_{m, n}$ by \eqref{eq:T-tFmn}. 
Recall also that $l(\bfu) = l(u_1) + \cdots + l(u_m)$ and $l(\bfv) = l(v_1) + \cdots + l(v_n)$.

\bthm{:Guv-leaves}
For any $\bfu \in W^m$ and $\bfv \in W^n$, 
symplectic leaves of $(G^{\bfu, \bfv}, \tpi_{m, n})$ are all the (non-empty)  level sets of the map  
$\chi_{{\bfu}, {\bfv}}$ in \eqref{eq:chi-1} and all have dimension equal to 
$l(\bfu) + l(\bfv) + \dim ({\rm Im}(1-uv^{-1}))$.
 Alternatively, 
\[
S^{{\bfu},{\bfv}} \stackrel{\rm def}{=}\{{\bf g} \in G^{\bfu, \bfv}\, \mbox{as in \eqref{eq:bf-g-hat-dot}}: \;  \; [b]_0 
\in \widetilde{T}^{\bfu, \bfv}, \, [b]_0 [b_-]_0^v \in T^{u, v}\}
\]
is a symplectic leaf of $(G^{\bfu, \bfv}, \tpi_{m, n})$, and
every symplectic leaf of $(G^{\bfu, \bfv}, \tpi_{m, n})$ is of the form 
\[
a \cdot S^{{\bfu}, {\bfv}}= \{{\bf g} \in G^{\bfu, \bfv}\, \mbox{as in \eqref{eq:bf-g-hat-dot}}: \;  \; [b]_0 
\in a^u\widetilde{T}^{\bfu, \bfv}, \, [b]_0 [b_-]_0^v \in (a^2)^uT^{u, v}\}
\]
for some $a \in T$. 
Furthermore, for $a_1, a_2 \in T$, $a_1 \cdot S^{{\bfu}, {\bfv}} = a_2 \cdot S^{{\bfu}, {\bfv}}$ 
if and only if $a_1^{-1}a_2 \in \widetilde{T}^{\bfu, \bfv}$ and
$(a_1^{-1}a_2)^2 \in T^{uv^{-1}}$. 
\ethm

\begin{proof}
For ${\bf g} \in G^{\bfu, \bfv}$ as in \eqref{eq:bf-g-hat-dot}, the $T$-equivariant Poisson isomorphism 
$K_{\bfu, \overline{\overline{\bfv}}}: (G^{\bfu, \bfv}, \tpi_{m, n}) \to (\O^{(\bfu, \bfv^{-1})}_e \times T, \pi_{m+n} \bowtie 0)$
in \eqref{eq:K-uv} is given by
\[
K_{\bfu, \overline{\overline{\bfv}}}({\bf g}) = ([c_1, \ldots, c_{m-1}, c_mb, \, (c_n^\prime)^{-1}, \cdots, 
(c_2^\prime)^{-1}, (c_1^\prime)^{-1}]_{\sF_{m+n}}, \, [b_-]_0).
\]
By \eqref{eq:CBBC}, for ${\bf g} \in G^{\bfu, \bfv}$ as in \eqref{eq:bf-g-hat-dot},
there are unique $({c}_{m+1}, \ldots, {c}_{m+n}) \in C_{\overline{v_n^{-1}}} \times \cdots \times 
C_{\overline{v_1^{-1}}}$ and $b_{m+1}, \ldots, b_{m+n} \in B$ such that
\[
b ({c}_n^\prime)^{-1} = {c}_{m+1} b_{m+1}, \;\; b_{m+1} (c_{n-1}^\prime)^{-1} = {c}_{m+2} b_{m+2}, \;\ldots, \;\;
b_{n+m-1} (c_1^\prime)^{-1} = {c}_{m+n} b_{m+n}.
\]
Then 
$[{c}_1, \ldots, {c}_{m-1}, {c}_mb, \, (c_n^\prime)^{-1}, \cdots, 
(c_2^\prime)^{-1}, (c_1^\prime)^{-1}] _{\sF_{m+n}} =
[{c}_1, \ldots, {c}_m, \, {c}_{m+1}, \ldots, {c}_{m+n}]_{\sF_{m+n}}$. Since
\[
[b(c_n^\prime)^{-1}, \cdots, 
(c_2^\prime)^{-1}, (c_1^\prime)^{-1}]_{\stF_n} = [{c}_{m+1}, \ldots, {c}_{m+n-1}, {c}_{m+n}b_{m+n}]_{\stF_n},
\]
one has ${c}_1 \cdots {c}_{m-1} {c}_m {c}_{m+1} \cdots {c}_{m+n-1} {c}_{m+n} = b_-b_{m+n}^{-1}$, and 
$[b_{m+n}]_0 = [b]_0^{v^{-1}}$.
Thus 
\[
[{c}_1 \cdots {c}_{m-1} {c}_m {c}_{m+1} \cdots {c}_{m+n-1} {c}_{m+n} ]_0 = [b_-b_{m+n}^{-1}]_0 = [b_-]_0 ([b]_0^{-1})^{v^{-1}}. 
\]
Let $\bfw = (\bfu, \bfv^{-1}) \in W^{m+n}$, so  $\overline{\bfw} = (\overline{u}_1, \ldots, \overline{u}_m, \overline{v_n^{-1}}, \ldots, 
\overline{v_1^{-1}},) \in \bfw T^{m+n}$. Let $\Sigma^{\overline{\bfw}} \subset \O^\bfw_e \times T$ be as in \eqref{eq:Sigma-bfw-bar}. It then follows from 
definitions that 
$K_{\bfu, \overline{\overline{\bfv}}}(S^{\bfu, \bfv})  = \Sigma^{\overline{\bfw}}$. 

 By  \thmref{:leaves-Ow},
 $\Sigma^{\overline{\bfw}}$ is a symplectic leaf of 
$(\O^{(\bfu, \bfv^{-1})}_e \times T, \pi_{m+n} \bowtie 0)$ of dimension equal to $l(\bfu) + l(\bfv) + {\rm dim} ({\rm Im}(1-uv^{-1}))$. 
It follows that  $S^{\bfu, \bfv}$ is a symplectic leaf of $(G^{\bfu, \bfv}, \tpi_{m, n})$ of the same dimension. 

2) follows from \thmref{:leaves-Ow} and the fact that $K_{\bfu, \overline{\overline{\bfv}}}$ is a $T$-equivariant Poisson isomorphism. 
\end{proof}

\bre{re:nunu}
{\rm 
The Poisson structure $\tpi_{m, n}$ is also invariant under the $T$-action 
\[
([g_1, \ldots, g_m]_{\stF_m}, [k_1, \ldots, k_n]_{\stF_{-n}}) \cdot t=
([g_1, \ldots, g_{m-1}, g_m t]_{\stF_m}, [k_1, \ldots, k_{n-1}, k_nt]_{\stF_{-n}}). 
\]
One checks directly that for any $\bfu \in W^m, \bfv \in W^n$, and $a \in T$, one has
\[
S^{\bfu, \bfv} \cdot a = \{{\bf g} \in G^{\bfu, \bfv}\, \mbox{as in \eqref{eq:bf-g-hat-dot}}: \;  \; [b]_0 
\in a\widetilde{T}^{\bfu, \bfv}, \, [b]_0 [b_-]_0^v \in a^2T^{u, v}\} = (a^{-1})^u \cdot S^{{\bfu}, {\bfv}}.
\]
\hfill $\diamond$
}
\ere

\subsection{Symplectic leaves of $(\tF_n, \tpi_n)$}\label{ss:leaves-tFn}
Consider now the $T$-Poisson variety $(\tF_n, \tpi_n)$ for $n \geq 1$, and recall that from \eqref{eq:tFn-T-leaves} that
$T$-leaves of $(\tF_n, \tpi_n)$
are precisely of the form 
\[
\widetilde{F}_n^{\bfu, v} = \{[g_1, g_2, \ldots, g_n]_{\stF_n} \in B\bfu B: g_1g_2 \cdots g_n \in B_-vB_-\},
\]
where $\bfu \in W^n$ and $v \in W$. The $T$-equivariant Poisson isomorphism
$(\tF_n, \tpi_n) \to (G_{n, 1}, \pi_{n, 1})$ in \eqref{eq:Fnn1} gives a $T$-equivariant Poisson isomorphism from 
$(\widetilde{F}_n^{\bfu, v}, \tpi_n)$ to $(G^{\bfu, v}, \tpi_{n, 1})$. We can thus use  \thmref{:Guv-leaves} to get a description of all
symplectic leaves of $(\widetilde{F}_n^{\bfu, v}, \tpi_n)$.

More precisely, let $\bfu = (u_1, u_2, \ldots, u_n) \in W^n$ and $v \in W$, and let $u = u_1u_2 \cdots u_n \in W$. Write an element in
$B\bfu B$ uniquely as $[c_1, c_2, \ldots, c_nb]_{\stF_n}$, where $(c_1, c_2, \ldots, c_n) \in C_{\overline{u}_1} \times C_{\overline{u}_2} \times \cdots \times
C_{\overline{u}_n}$ and $b \in B$, and let $\underline{c} = c_1c_2 \cdots c_n$. Let
\[
\Lambda^{\bfu, v}=\{[c_1, c_2, \ldots, c_nb]_{\stF_n} \in B\bfu B:  
\underline{c}\,b \in B_-vB_-, \,
[b]_0 \in \widetilde{T}^{\bfu, v}, [b]_0[\underline{c}\,b\,\overline{v^{-1}}]_0^v \in T^{u, v}\}.
\]

\bthm{:leaves-tFn} For any $\bfu \in W^n$ and $v \in W$,

1) $\Lambda^{\bfu, v}$ is a symplectic leaf of $(\tF_n^{\bfu, v}, \tpi_n)$ of dimension
$l(\bfu) + l(v) + \dim ({\rm Im}(1-uv^{-1}))$; 

2) every symplectic leaf of $(\tF_n^{\bfu, v}, \tpi_n)$ is of the form $a \cdot \Lambda^{\bfu, v}$ for some $a \in T$. Moreover, for $a_1, a_2 \in T$,
$a_1 \cdot \Lambda^{\bfu, v} = a_2 \cdot \Lambda^{\bfu, v}$ if and only if $a_1^{-1}a_2 \in \widetilde{T}^{\bfu, v}$ and $(a_1^{-1}a_2)^2 \in T^{uv^{-1}}$.
\ethm

\end{document}